\documentclass[11pt]{article}
\usepackage{latexsym}
\usepackage[all]{xy}
\usepackage{amsmath}
\usepackage{amssymb}
\usepackage{amsfonts}

\newtheorem{thm}{Th\'eor\`eme}[section]
\newtheorem{prop}[thm]{Proposition}
\newtheorem{lem}[thm]{Lemme}

\newtheorem{df}[thm]{D\'efinition}
\newtheorem{cor}[thm]{Corollaire}

\newtheorem{rmk}[thm]{Remarque}

\newtheorem{hyp}[thm]{Hypoth\`ese}

\begin{document}

\title{\textbf{Au-dessous de $Spec\, \mathbb{Z}$}}
\bigskip
\bigskip

\author{\bigskip\\
Bertrand To\"en et Michel Vaqui\'e\\
\small{Laboratoire Emile Picard}\\
\small{UMR CNRS 5580} \\
\small{Universit\'{e} Paul Sabatier, Toulouse}\\
\small{France}}

\date{Octobre 2007}

\maketitle

\begin{abstract}
Dans ce travail nous utilisons les th\'eories de g\'eom\'etrie alg\'ebrique relative 
et de g\'eom\'etrie alg\'ebrique homotopique (voir \cite{hagII}) afin de construire 
plusieurs cat\'egories de sch\'emas d\'efinis \emph{au-dessous de $Spec\, \mathbb{Z}$}. 
Nous d\'efinissons ainsi les cat\'egories de $\mathbb{N}$-sch\'emas, $\mathbb{F}_{1}$-sch\'emas, 
$\mathbb{S}$-sch\'emas, $\mathbb{S}_{+}$-sch\'emas, et $\mathbb{S}_{1}$-sch\'emas, 
o\`u (d'un point de vue tr\`es intuitif) $\mathbb{N}$ est le semi-anneau des entiers naturels, 
$\mathbb{F}_{1}$ est le corps \`a un \'el\'ement, $\mathbb{S}$ est l'anneau en spectres des entiers,
$\mathbb{S}_{+}$ est le semi-anneau en spectres des entiers naturels et 
$\mathbb{S}_{1}$ est l'anneau en spectres \`a un \'el\'ement. 
Ces cat\'egories de sch\'emas sont reli\'ees entre elles \`a l'aide de foncteurs
de changement de bases, et poss\`edent toutes un foncteur de changement 
de bases vers les $\mathbb{Z}$-sch\'emas. Nous montrons comment 
les groupes lin\'eaires $Gl_{n}$ et les vari\'et\'es toriques peuvent-\^etre d\'efinis
comme certains objets de ces cat\'egories. 
\end{abstract}

\medskip

\tableofcontents

\bigskip

\section{Introduction}

Le but de ce travail est de construire plusieurs cat\'egories de \emph{sch\'emas}Ê
qui sont d\'efinis sur des bases se trouvant \emph{au-dessous de $Spec\, \mathbb{Z}$}. 
Bien entendu, comme $\mathbb{Z}$ est l'anneau commutatif initial, il est indispensable
de sortir du cadre usuel des anneaux et de s'autoriser \`a utiliser des
objets plus g\'en\'eraux mais qui ressemblent suffisamment \`a des anneaux 
commutatifs afin que la notion de sch\'ema puisse \^etre d\'efinie. Notre approche 
\`a cette question est bas\'ee sur la th\'eorie de la g\'eom\'etrie alg\'ebrique relative, 
largement insipir\'ee de \cite{ha}. 
Elle consiste \`a remarquer qu'un anneau commutatif n'est rien d'autre qu'un 
mono\"\i de commutatif dans la cat\'egorie mono\"\i dale des $\mathbb{Z}$-modules, et qu'en g\'en\'eral pour
une cat\'egorie mono\"\i dale sym\'etrique $(C,\otimes,\mathbf{1})$ les mono\"\i des commutatifs
dans $C$ peuvent \^etre pens\'es comme des mod\`eles pour les \emph{sch\'emas affines relatifs \`a $C$}. Il est remarquable qu'une
approche si g\'en\'erale (voire simpliste) permette effectivement de d\'efinir une notion de sch\'ema, et de plus de fa\c{c}on
fonctorielle en $C$. Ainsi, en choisissant $C$ munie d'un foncteur mono\"\i dal sym\'etrique
$C \longrightarrow \mathbb{Z}-Mod$ raisonnable, on trouve une notion de sch\'ema
relatif \`a $C$ et un foncteur de changement de bases vers les $\mathbb{Z}$-sch\'emas, et donc une
notion de sch\'ema au-dessous de $Spec\, \mathbb{Z}$.  Dans cet article nous montrerons comment cette 
approche, ainsi que sa g\'en\'eralisation homotopique o\`u $C$ est de plus munie d'une structure
de cat\'egorie de mod\`eles de Quillen, 
permettent de d\'efinir cinq cat\'egories de sch\'emas se trouvant en dessous de $Spec\, \mathbb{Z}$. \\

\begin{center} \textit{G\'eom\'etrie alg\'ebrique relative} \end{center}

Les id\'ees g\'en\'erales de la g\'eom\'etrie alg\'ebrique relative remontent
\`a \cite{ha}, o\`u des sch\'emas relatifs au-dessus d'un topos annel\'e
sont d\'efinis. Dans \cite{del} le cas des sch\'emas au-dessus 
d'une cat\'egorie Tannakienne est aussi consid\'er\'e. La th\'eorie
de g\'eom\'etrie alg\'ebrique relative que nous allons pr\'esenter est largement
inspir\'ee de ses deux r\'ef\'erences, bien que consid\'erer des
cat\'egories de bases qui ne sont pas ab\'eliennes, ni m\^eme 
additives, semble une nouveaut\'e. 

Donnons-nous une cat\'egorie mono\"\i dale sym\'etrique $(C,\otimes,\mathbf{1})$ que l'on supposera
compl\`ete, cocompl\`ete et ferm\'ee (i.e. poss\`ede des Hom internes relatifs \`a la structure
mono\"\i dale $\otimes$). Il est bien connu que l'on dispose dans $C$ d'une notion 
de mono\"\i de, pour un tel mono\"\i de $A$ d'une notion de module et pour 
un morphisme de mono\"\i des $A\longrightarrow B$ d'un foncteur $-\otimes_{A}B$ de changement de bases
des $A$-modules vers les $B$-modules (voir par exemple \cite{sav}).
En particulier il existe une notion de mono\"\i de commutatif (associatif et unitaire) dans $C$, 
et ils forment une cat\'egorie que l'on note $Comm(C)$. On d\'efinit formellement la cat\'egorie
des sch\'emas affines relatifs \`a $C$ par $Aff_{C}:=Comm(C)^{op}$. Tout ceci est pour le moment tr\`es formel, 
mais il se produit alors plusieurs miracles. 
\begin{itemize}
\item Il existe une topologie de Grothendieck naturelle sur $Aff_{C}$ appel\'ee la \emph{topologie plate}. 
Les familles couvrantes $\{X_{i} \longrightarrow X\}$ pour cette topologie correspondent aux familles
finies de morphismes $\{A \longrightarrow A_{i}\}$ dans $Comm(C)$ telles que le foncteur de changement 
de bases sur les cat\'egories de modules
$$\prod_{i} -\otimes_{A}A_{i} : A-Mod \longrightarrow \prod_{i}A_{i}-Mod$$
soit exact et conservatif. 
\item La topologie plate sur $Aff_{C}$ ainsi d\'efinie est sous-canonique (i.e. les pr\'efaisceaux repr\'esentables
sont des faisceaux). 
\item  Il existe une notion d'ouvert de Zariski dans $Aff_{C}$, qui par d\'efinition sont les
morphismes $f : X \longrightarrow Y$ dont le morphisme $A \longrightarrow B$ correspondant dans $Comm(C)$
satisfait aux trois conditions suivantes. 
\begin{enumerate}
\item (\emph{f est un monomorphisme}) Pour tout $A'\in Comm(C)$ le morphisme induit 
$Hom(B,A') \longrightarrow Hom(A,A')$ est injectif. 
\item (\emph{f est plat}) 
Le foncteur de changement de bases 
$$-\otimes_{A}B : A-Mod \longrightarrow B-Mod$$
est exact. 
\item (\emph{f est de pr\'esentation finie}) Pour tout diagramme filtrant d'objets $C_{i} \in A/Comm(C)$, le morphisme naturel
$$Colim Hom_{A/Comm(C)}(B,C_{i}) \longrightarrow Hom_{A/Comm(C)}(B,Colim C_{i})$$
est bijectif.

\end{enumerate}
\item La notion d'ouvert de Zariski s'\'etend de fa\c{c}on naturelle aux morphismes
entre faisceaux quelconques (voir d\'efinition \ref{d3}). 
\item Les ouverts de Zariski sont stables par compositions, isomorphismes et changements de bases. 
\item Les ouverts de Zariski donnent lieu \`a une notion de topologie de Zariski, et celle-ci est
encore sous-canonique.

\end{itemize}

Les propri\'et\'es ci-dessus sont tout ce dont on a besoin pour d\'efinir une cat\'egorie
de sch\'emas relatifs \`a $(C,\otimes,\mathbf{1})$. En effet, un sch\'ema
relatif est par d\'efinition un faisceau sur le site $Aff_{C}$ muni de la topologie de Zariski, et qui 
poss\`ede un recouvrement Zariski par des sch\'emas affines
(voir d\'efinition \ref{sch}). 
La cat\'egorie des sch\'emas ainsi obtenue est not\'ee $Sch(C)$. C'est une sous-cat\'egorie
pleine, stable par produits fibr\'es et sommes disjointes, de la cat\'egorie
des faisceaux sur $Aff_{C}$. De plus, elle contient une sous-cat\'egorie pleine
de sch\'emas affines, qui sont exactement les faisceaux repr\'esentables, et 
qui est naturellement \'equivalente \`a la cat\'egorie $Comm(C)^{op}$, oppos\'ee
de la cat\'egorie des mono\"\i des commutatifs dans $C$ (voir \S 2.2). Enfin, la nature
purement cat\'egorique de la construction rend la cat\'egorie
$Sch(C)$ fonctorielle en $C$, tout au moins pour des adjoints \`a gauches
sym\'etriques mono\"\i daux $(C,\otimes,\mathbf{1})\longrightarrow (D,\otimes,\mathbf{1})$
satisfaisant \`a quelques conditions faciles \`a v\'erifier dans la pratique (voir $\S 2.3$). \\

\begin{center} \textit{Trois exemples de g\'eom\'etries alg\'ebriques relatives} \end{center}

Nous consid\`ererons trois exemples de g\'eom\'etries alg\'ebriques 
relatives, correspondant \`a trois choix pour $(C,\otimes,\mathbf{1})$. 
Tout d'abord on posera  
$(C,\otimes,\mathbf{1})=(\mathbb{Z}-Mod,\otimes,\mathbb{Z})$, la cat\'egorie mono\"\i dale sym\'etrique
des groupes ab\'eliens (pour le produit tensoriel). La cat\'egorie des sch\'emas
ainsi obtenue $\mathbb{Z}-Sch$ se trouve \^etre \'equivalente \`a la cat\'egorie
des sch\'emas au sens usuel. Ce fait justifie 
notre terminologie de \emph{sch\'emas relatifs}. 

Notre deuxi\`eme exemple sera $(C,\otimes,\mathbf{1})=(\mathbb{N}-Mod,\otimes,\mathbb{N})$
la cat\'egorie des mono\"\i des commutatifs, ou encore des semi-groupes ab\'eliens (pour
la notion naturelle de produit tensoriel), que l'on pourrait aussi appeler
des \emph{$\mathbb{N}$-modules}. La cat\'egorie des sch\'emas
sera dans ce cas not\'ee $\mathbb{N}-Sch$, et la sous-cat\'egorie des sch\'emas affines
est \'equivalente \`a la cat\'egorie oppos\'ee de celle des 
semi-anneaux commutatifs. Elle est munie d'une adjonction
$$i : \mathbb{Z}-Sch \longrightarrow \mathbb{N}-Sch
\qquad 
 \mathbb{Z}-Sch \longleftarrow \mathbb{N}-Sch : -\otimes_{\mathbb{N}}\mathbb{Z},$$
et l'adjoint \`a gauche $i$ est pleinement fid\`ele. Le foncteur
$-\otimes_{\mathbb{N}}\mathbb{Z}$ quant \`a lui est une globalisation 
du foncteur de compl\'etion en groupes des semi-anneaux commutatifs
vers les anneaux commutatifs. 

Enfin, notre troisi\`eme exemple est $(C,\otimes,\mathbf{1})=(Ens,\times,*)$, la cat\'egorie mono\"\i dale
sym\'etrique des ensembles (pour le produit direct). La cat\'egorie des sch\'emas relatifs sera
not\'ee $\mathbb{F}_{1}-Sch$, et l'on y pense comme des mod\`eles pour
des \emph{vari\'et\'es d\'efinies sur le corps \`a un \'el\'ement}, au sens 
o\`u cette notion apparait dans \cite{so}. Par d\'efinition la sous-cat\'egorie des
$\mathbb{F}_{1}$-sch\'emas affines est \'equivalente \`a l'oppos\'ee de la cat\'egorie
des mono\"\i des commutatifs. 
On dispose d'un foncteur de changement de bases
$$-\otimes_{\mathbb{F}_{1}}\mathbb{N} : \mathbb{F}_{1}-Sch \longrightarrow \mathbb{N}-Sch,$$
qui est une globalisation du foncteur qui envoit un mono\"\i de commutatif $M$ sur 
son semi-anneau en mono\"\i des $\mathbb{N}[M]$, analogue \emph{semi} 
des anneaux en groupes. En composant avec le changement de base pr\'ec\'edent, on obtient 
deux foncteurs
$$\xymatrix{
 \mathbb{F}_{1}-Sch \ar[r]^-{-\otimes_{\mathbb{F}_{1}}\mathbb{N}} & \mathbb{N}-Sch
 \ar[r]^-{-\otimes_{\mathbb{N}}\mathbb{Z}} &\mathbb{Z}-Sch,}$$
que l'on peut repr\'esenter sch\'ematiquement par le diagramme
$$\xymatrix{
Spec\, \mathbb{Z} \ar[r] & Spec\, \mathbb{N} \ar[r] & Spec\, \mathbb{F}_{1}.}$$

Comme exemple de sch\'emas relatifs, nous montrerons que les vari\'et\'es toriques
sont naturellement d\'efinies sur $\mathbb{F}_{1}$ (voir $\S 4.2$). Ceci est tr\`es naturel car elles sont
obtenues par recollements  formels de mono\"\i des commutatifs. Nous montrons aussi qu'il existe des 
sch\'emas $Gl_{n,\mathbb{N}}\in \mathbb{N}-Sch$ et 
$Gl_{n,\mathbb{F}_{1}}\in \mathbb{F}_{1}-Sch$ qui sont des versions du sch\'ema en groupes lin\'eaires. 
Cependant, bien que $Gl_{n,\mathbb{N}}\otimes_{\mathbb{N}}\mathbb{Z}\simeq Gl_{n,\mathbb{Z}}$, 
il n'est pas vrai que $Gl_{n,\mathbb{F}_{1}}\otimes_{\mathbb{F}_{1}}\mathbb{Z}$ soit isomorphe
\`a $Gl_{n,\mathbb{Z}}$ (contrairement \`a ce que l'on pourrait attendre). 
Les sch\'emas 
$Gl_{n,\mathbb{F}_{1}}$ et $Gl_{n,\mathbb{N}}$ ont ceci de remarquable que
$$Gl_{n,\mathbb{F}_{1}}(\mathbb{F}_{1})\simeq Gl_{n,\mathbb{N}}(\mathbb{N})\simeq \Sigma_{n},$$
ce qui montre que $Gl_{n,\mathbb{N}}(\mathbb{N})$ est une extension non-triviale
de $Gl_{n,\mathbb{Z}}$ au-dessus de $Spec\, \mathbb{N}$ (i.e. 
$i(Gl_{n,\mathbb{Z}})\ncong Gl_{n,\mathbb{N}}$). Cet exemple montre 
que la propri\'et\'e \emph{d'\^etre d\'efini sur} n'est pas tr\`es pertinente (car 
tout $\mathbb{Z}$-sch\'ema est d\'efini sur $\mathbb{N}$), et le fait 
int\'eressant est en r\'ealit\'e que certains $\mathbb{Z}$-sch\'emas 
poss\`edent des mod\`eles naturels d\'efinis sur $\mathbb{N}$ ou sur
$\mathbb{F}_{1}$. La situation est donc tout \`a fait comparable
avec ce qu'il se passe en g\'eom\'etrie alg\'ebrique d\'eriv\'ee
dont la pertinence r\'eside dans le fait que les espaces
de modules poss\`edent des extensions naturelles et non-triviales
en des espaces de modules d\'eriv\'es (voir \cite{hagdag}). \\

\begin{center} \textit{G\'eom\'etrie alg\'ebrique relative homotopique} \end{center}

Dans une derni\`ere section nous introduirons trois nouvelles cat\'egories
de sch\'emas d\'efinies \`a l'aide du formalisme de la g\'eom\'etrie alg\'ebrique
homotopique de \cite{hagII}. L'id\'ee g\'en\'erale est que l'on peut g\'en\'eraliser
la g\'eom\'etrie alg\'ebrique relative en supposant de plus que 
la cat\'egorie $C$ est munie d'une structure de mod\`eles compatible
avec sa structure mono\"\i dale (le cas particulier de la g\'eom\'etrie alg\'ebrique
relative non-homotopique se retrouve en prenant la structure de mod\`eles
triviale pour laquelle les \'equivalences sont les isomorphismes). Il se trouve que les
notions de topologie plate et d'ouverts de Zariski gardent un sens, bien que dans le cas
g\'en\'eral nous ne savons pas montrer que la topologie plate est sous-canonique
(nous le montrerons cependant pour les exemples qui nous int\'eressent). 
On dispose donc d'une notion de sch\'emas relatifs \`a $C$, dont la cat\'egorie
sera not\'ee $Sch(C)$ (voir $\S 5$ pour plus de d\'etails). 

Ceci nous permettra de trouver trois nouvelles notions de sch\'emas au-dessous de 
$Spec\, \mathbb{Z}$ en trouvant
trois exemples de cat\'egories de mod\`eles mono\"\i dales sym\'etriques $C$ munies de foncteurs
de Quillen \`a gauche mono\"\i daux $C \longrightarrow \mathbb{Z}-Mod$. Le premier de ces exemples
est lorsque l'on pose $(C,\otimes,\mathbf{1})=(\mathcal{GS},\wedge,\mathbb{S})$, la cat\'egorie
de mod\`eles mono\"\i dale des $\Gamma$-espaces tr\`es sp\'eciaux. Cette cat\'egorie
de mod\`eles est un mod\`ele pour la th\'eorie homotopique des spectres connectifs (i.e. sans
homotopie n\'egative), que l'on peut voir
comme des analogues homotopiques des groupes ab\'eliens. On dispose ainsi d'une
cat\'egorie de sch\'emas relatifs \`a $\mathcal{GS}$, que l'on notera
$\mathbb{S}-Sch$, o\`u la notation $\mathbb{S}$ rappelle le spectre en sph\`ere. Les 
$\mathbb{S}$-sch\'emas affines sont en correspondance avec les
$E_{\infty}$-anneaux en spectres connectifs (i.e. sans homotopie n\'egative), qui
souvent portent le nom de \emph{(commutative) brave new rings} dans la litt\'erature. 
Les objets de $\mathbb{S}-Sch$ peuvent donc s'appeler 
\emph{brave new schemes}, ou encore 
d'apr\`es une suggestion amusante de J. Tapia des 
\emph{nouveaux sch\'emas courageux}. 

Notre second exemple est $(C,\otimes,\mathbf{1})=(SEns,\times,*)$, la cat\'egorie 
de mod\`eles des ensembles simpliciaux munie de son produit direct. Les sch\'emas que l'on obtient
sont des versions homotopiques des $\mathbb{F}_{1}$-sch\'emas, et seront appel\'es
des $\mathbb{S}_{1}$-sch\'emas, la notation $\mathbb{S}_{1}$ signifiant intuitivement
\emph{l'anneau en spectres \`a un \'el\'ement}, bien que $\mathbb{S}_{1}$ ne soit ni un anneau 
ni un spectre (mais ceci est compatible avec la terminologie de \emph{corps \`a un \'el\'ement} qui
ne d\'esigne pas un corps). 

Enfin, notre dernier exemple est $(C,\otimes,\mathbf{1})=(\mathcal{MS},\wedge,\mathbb{S}_{+})$, 
la cat\'egorie de mod\`eles des $\Gamma$-espaces sp\'eciaux. La cat\'egorie de mod\`eles
$\mathcal{MS}$ est un mod\`ele pour la th\'eorie homotopique des $E_{\infty}$-mono\"\i des 
simpliciaux. La cat\'egorie des sch\'emas relatifs sera not\'ee $\mathbb{S}_{+}-Sch$, et 
ses objets affines sont en correspondance avec les \emph{semi-anneaux en spectres commutatifs}
(bien que ceux-ci ne soient pas des spectres), ou encore des \emph{brave new commutative semi-rings}, 
analogues
homotopiques des semi-anneaux commutatifs. La notation $\mathbb{S}_{+}$ signifie intuitivement
\emph{le semi-anneau en spectres des entiers positifs}, et est une version homotopique
du semi-anneau $\mathbb{N}$.

Pour terminer, on dispose de foncteurs de changement de bases naturels entre toutes ces notions, et
d'un diagramme commutatif 
$$\xymatrix{
\mathbb{S}_{1}-Sch \ar[rr]^-{-\otimes_{\mathbb{S}_{1}}\mathbb{S}_{+}} \ar[dd]_-{-\otimes_{\mathbb{S}_{1}}\mathbb{F}_{1}} & & \mathbb{S}_{+}-Sch 
\ar[dd]^-{-\otimes_{\mathbb{S}_{+}}\mathbb{N}} \ar[rr]^-{-\otimes_{\mathbb{S}_{+}}\mathbb{S}} & & \mathbb{S}-Sch
\ar[dd]^-{-\otimes_{\mathbb{S}}\mathbb{Z}} \\
 &  & & & \\
\mathbb{F}_{1}-Sch \ar[rr]_-{-\otimes_{\mathbb{F}_{1}}\mathbb{N}} & & 
\mathbb{N}-Sch \ar[rr]_-{-\otimes_{\mathbb{N}}\mathbb{Z}} & & \mathbb{Z}-Sch.}$$
Sch\'ematiquement on repr\'esente ce diagramme par le diagramme suivant
$$\xymatrix{
Spec\, \mathbb{Z} \ar[r] \ar[d] & Spec\, \mathbb{N} \ar[r] \ar[d] & Spec\, \mathbb{F}_{1} \ar[d] \\
Spec\, \mathbb{S} \ar[r] & Spec\, \mathbb{S}_{+} \ar[r] & Spec\, \mathbb{S}_{1}.}$$
On montrera que $Gl_{n}$ peut \^etre d\'efini sur $\mathbb{S}_{+}$, et que les vari\'et\'es toriques 
le sont sur $\mathbb{S}_{1}$, et donc par changement de bases sur $\mathbb{S}_{+}$ et 
$\mathbb{S}$. \\

\begin{center} \textit{Ce que nous n'avons pas fait} \end{center}

Notre but principal dans ce travail \'etait d'amorcer une \'etude syst\'ematique de la notion
de sch\'emas d\'efinis au-dessous de $Spec\, \mathbb{Z}$, et de montrer que les
techniques de g\'eom\'etrie alg\'ebrique relative et de g\'eom\'etrie alg\'ebrique homotopique
semblent bien adapt\'ees pour aborder cette question. Cependant, nous convenons que ce
texte ne contient que quelques pr\'emisses de g\'eom\'etrie alg\'ebrique au-dessous
de $Spec\, \mathbb{Z}$, et n\'ecessiterait de nombreux compl\'ements afin d'aboutir \`a une
th\'eorie riche et int\'eressante. Dans cet ordre d'id\'ees, signalons quelques 
questions que nous pensons importantes et qui ne sont pas trait\'ees dans 
ce travail.

Pour commencer, nous ne nous sommes pas attaqu\'es \`a la recherche d'une description plus
explicite de la topologie de Zariski d'un sch\'ema relatif (tel que pr\'esent\'ee \`a la fin
du paragraphe \S 2.2). Il devrait toute fois \^etre possible de d\'ecrire l'espace topologique
$|Spec\, A|$, sous-jacent \`a un sch\'ema relatif affine, en des termes
plus standards utilisant par exemple une notion d'id\'eal dans $A$ (i.e. 
des sous-objets du $A$-module $A$), tout au moins
si la cat\'egorie de base $C$ satisfait \`a quelques hypoth\`eses
suppl\'ementaires (mais raisonables). Dans le cadre de la g\'eom\'etrie alg\'ebrique
homotopique cette question semble reli\'ee de tr\`es pr\`es \`a la notion 
de topologie enrichie de \cite{ve}. Une telle description ne nous semble pas vraiment
indispensable, mais elle permettrait de comparer plus facilement notre notion 
de sch\'emas \`a celles d\'ej\`a existantes, comme par exemple 
la notion de $\mathbb{F}_{1}$-sch\'emas  de \cite{de}. 

De nombreux espaces de modules ne sont pas des sch\'emas, mais 
des espaces alg\'ebriques voire des champs alg\'ebriques. Il nous semble
donc capital de d\'evelopper aussi une notion de \emph{champs alg\'ebriques
relatifs} afin de pouvoir constuire de nouveaux exemples d'objets
g\'eom\'etriques d\'efinis au-dessous de $Spec\, \mathbb{Z}$.
Ceci n\'ecessite bien \'evidemment une notion de morphismes \'etales et 
de morphismes lisses. Dans le cas o\`u la cat\'egorie de mod\`eles $C$ est additive (ou m\^eme
additive \`a homotopie pr\`es), de telles notions ont \'et\'e introduites et \'etudi\'ees
dans \cite{hagII}, et dans ce cas on arrive \`a d\'efinir une notion int\'eressante
de champs g\'eom\'etriques relatifs (voir \cite{hagII} pour des exemples de tels objets). 
Les deux seuls mod\`eles pr\'esent\'es dans ce travail dont les
cat\'egories de bases sont additives sont les $\mathbb{Z}$-sch\'emas et 
les $\mathbb{S}$-sch\'emas. Ainsi, on dispose d'une notion raisonable
de \emph{$\mathbb{S}$-champs g\'eom\'etriques}, et 
on peut montrer par exemple que les champs de modules de repr\'esentations
de carquois poss\`edent des mod\`eles naturels qui sont des 
$\mathbb{S}$-champs g\'eom\'etriques. Plus g\'en\'eralement, on peut
d\'emontrer un th\'eor\`eme de repr\'esentabilit\'e analogue \`a celui 
d\'emontr\'e dans \cite{tv} qui fournit des exemples 
de $\mathbb{S}$-champs g\'eom\'etriques associ\'es \`a des
\emph{cat\'egories spectrales} (cat\'egories enrichies dans les spectres).

En dehors des deux mod\`eles $Spec\, \mathbb{Z}$ et 
$Spec\, \mathbb{S}$, les notions de morphismes
\'etales et de morphismes lisses d\'evelopp\'ees
dans \cite{hagII} ne sont plus disponibles.  Il nous semble crucial
de chercher \`a d\'evelopper de telles notions dans le cas
non-additif, ce qui permettrait d'avoir acc\`es \`a la notion
de champs g\'eom\'etriques relatifs, et par exemple de montrer que
les champs de modules de repr\'esentations de carquois sont 
en fait d\'efinis sur $Spec\, \mathbb{S}_{+}$. D'autres part, 
la notion de lissit\'e semble aussi tr\`es importante si l'on souhaite
d\'efinir une notion de motifs pour les sch\'emas relatifs, qui 
comme il est expliqu\'e dans \cite{so} est largement souhait\'ee
pour les $\mathbb{F}_{1}$-sch\'emas. 

Nous n'avons pas \'etudi\'e les notions de faisceaux 
quasi-coh\'erents et coh\'erents sur les sch\'emas
relatifs. La notion de faisceau quasi-coh\'erent est 
claire et facile \`a deviner. En contre partie, 
la notion de coh\'erence semble plus d\'elicate car elle fait
intervenir des conditions de finitudes qui ne parraissent pas
\'evidentes \`a g\'en\'eraliser au cadre des sch\'emas relatifs. 
Ces conditions de finitudes, et la notion de faisceaux coh\'erents
qui en d\'ecoule sont importantes par exemple pour pouvoir
d\'efinir des groupes de K-th\'eorie de sch\'emas relatifs. 
Il serait par exemple int\'eressant de montrer que 
les calculs faits dans \cite{hut} donnent effectivement une
description de la K-th\'eorie de nos $\mathbb{S}$-vari\'et\'es 
toriques $X_{\mathbb{S}}(\Delta)$ (voir \S 5.3), 
comme il est sugg\'er\'e par l'auteur que cela est la cas
si l'on sait d\'efinir ce qu'est un mod\`ele sur $\mathbb{S}$
des vari\'et\'es toriques.

Notre notion de $\mathbb{N}$-sch\'emas n'est peut-\^etre 
pas sans relations avec la g\'eom\'etrie tropicale
(voir par exemple \cite{trop}), qui elle aussi utilise de fa\c{c}on essentielle
des semi-anneaux commutatifs. Il serait int\'eressant de trouver
des relations pr\'ecises entre la g\'eom\'etrie
des $\mathbb{N}$-sch\'emas de type fini sur le semi-anneau tropical
$\mathbb{R}_{trop}$ et les vari\'et\'es tropicales que l'on rencontre dans la litt\'erature.

\bigskip

\textbf{Remerciements:} Nous remercions G. Vezzosi pour plusieurs
conversations sur des sujets connexes au cours de ces derni\`eres
ann\'ees. Nous remercions aussi J. Kock pour de ses remarques nombreuses, 
pr\'ecises et utiles. Un grand merci enfin \`a M. Anel pour son int\'eret
constant aux dessous en g\'en\'eral et plus particuli\`erement \`a ceux
de $Spec\, \mathbb{Z}$. \\

\bigskip

\textbf{Convention:} Tous les mono\"\i des consid\'er\'es seront unitaires et 
associatifs. Tous les modules sur des mono\"\i des seront alors unitaires. 
De m\^eme toutes les cat\'egories mono\"\i dales seront munies
de contraintes d'unit\'e et d'associativit\'e. 

Nous n\'egligerons les prob\`emes ensemblistes li\'es aux choix d'univers. 
Le lecteur pourra consulter \cite{hagI,hagII} o\`u il trouvera une
fa\c{c}on de les r\'esoudre.

\bigskip

\section{G\'eom\'etrie alg\'ebrique relative}

Le but de cette premi\`ere partie est de pr\'esenter la notion de sch\'ema
relatif \`a une cat\'egorie sym\'etrique mono\"\i dale de base $C$. Nous commencerons
par un proc\'ed\'e g\'en\'eral de construction de topologies de Grothendieck 
\`a partir de pr\'echamps en cat\'egories v\'erifiant certaines conditions. Cela nous
permettra par la suite de d\'efinir la \emph{topologie fid\`element plate et quasi-compacte}, ainsi que
la \emph{topologie de Zariksi} dans des contextes tr\`es g\'en\'eraux. 
Nous d\'efinirons alors la notion de sch\'ema relatif en recollant 
des objets affines \`a l'aide de la topologie de Zariski. 

\subsection{Construction de topologies de Grothendieck}

Nous nous donnons une cat\'egorie $T$ qui poss\`ede des limites finies et 
un pseudo foncteur 
$$\begin{array}{cccc}
M : & T^{op} & \longrightarrow & Cat \\
& X & \longmapsto & M(X) 
\end{array}$$ 
qui v\'erifie les conditions suivantes. 
\begin{hyp}\label{h0}
\begin{enumerate} 
\item Pour tout $X$ dans $T$ la cat\'egorie $M(X)$ poss\`ede des limites 
et des colimites arbitraires. 

\item Pour tout $p:X'\to X$ dans $T$ le foncteur $M(p)= p^* : M(X) \to M(X')$ 
poss\`ede un adjoint \`a droite $p_* : M(X') \to M(X)$ qui est conservatif (i.e. 
un morphisme $u : x \rightarrow y$ dans $M(X')$ est un isomorphisme si et 
seulement si le morphisme induit $p_*(u) : p_*(x) \rightarrow p_*(y)$ est 
un isomorphisme dans $M(X)$).  
\item Pour tout diagramme cart\'esien dans $T$
$$\xymatrix{
Y' \ar[r]^{p'} \ar[d]^{q'} & Y \ar[d]^{q} \\
X' \ar[r]_{p} & X}$$ 
la transformation naturelle de changement de base
$$p^*q_* \Rightarrow q'_*p'^*$$
est un isomorphisme.  
\end{enumerate}  
\end{hyp}

Un mot sur la condition $(3)$ ci-dessus. La transformation naturelle
en question est construite de la fa\c{c}on suivante: on dispose  
d'isomorphismes naturels 
$$(q')^{*}p^{*}\simeq (pq')^{*}=(qp')^{*}\simeq (p')^{*}q^{*}$$ provenant de
la structure de pseudo-foncteur de $M$. Cela nous donne en composant 
\`a droite par $q_{*}$ un isomorphisme naturel
$(q')^{*}p^{*}q_{*} \simeq (p')^{*}q^{*}q_{*}$. En composant ce dernier
avec la co-unit\'e d'adjonction $q^{*}q_{*} \Rightarrow id$ on trouve une
transformation naturelle
$$(q')^{*}p^{*}q_{*} \Rightarrow (p')^{*},$$
qui par adjonction nous donne la transformation naturelle de changement de base
$$p^{*}q_{*} \Rightarrow q'_{*}(p')^{*}.$$

\begin{rmk}
\emph{L'exemple fondamental que nous avons en t\^ete est le suivant: $T$ est la
cat\'egorie des sch\'emas affines, et pour $X$ un tel sch\'ema $M(X)$ est la cat\'egorie 
des faisceaux quasi-coh\'erents sur $X$.}
\end{rmk}

Nous pouvons alors d\'efinir les notions de platitude et de fid\`ele 
platitude sur la cat\'egorie $T$ associ\'ees au pseudo foncteur $M$. 

\begin{df}\label{dmtop}
Soit $\{p_{i} : X_{i}\to X\}_{i\in I}$ une famille de morphismes dans $T$. 
\begin{enumerate}  
\item La famille $\{p_{i} : X_{i}\to X\}_{i\in I}$ est \emph{$M$-couvrante} s'il existe
un ensemble fini $J\subset I$ telle que la famille de foncteurs 
$$\{p_{i}^* : M(X) \to M(X_{i})\}_{i\in J}$$ 
soit conservative.  
\item La famille $\{p_{i} : X_{i}\to X\}_{i\in I}$ est \emph{$M$-plate} si tous les foncteurs 
$p_{i}^*:M(X)\to M(X_{i})$ 
sont exacts \`a gauche (i.e. commutent aux limites finies). 
\item La famille $\{p_{i} : X_{i}\to X\}_{i\in I}$ est \emph{$M$-fid\`element plate} si elle est \`a la fois 
$M$-couvrante et $M$-plate.  
\end{enumerate} 
\end{df}  

Dans la suite, nous utiliserons aussi la terminologie de \emph{morphisme
$M$-plat}, qui, comme le lecteur peut le deviner, est un morphisme
qui forme une famille (\`a un unique \'el\'ement) $M$-plate au sens
pr\'ec\'edent. En d'autre termes, $p : X'\to X$ est $M$-plat si le
foncteur $p^{*}$ est exact. Remarquons par ailleurs que pour un morphisme
$M$-plat $p : X'\to X$ le foncteur $p^{*}$ est en r\'ealit\'e exact. En effet
l'existence de l'adjoint \`a droite implique que $p^{*}$ commute aux colimites
arbitraires. \\

Comme la cat\'egorie $M(X)$ a des limites finies, un foncteur $p^*$ qui est \`a 
la fois exact \`a gauche et conservatif est fid\`ele, ce qui justifie la 
notion de \emph{$M$-fid\`element plat}.   

\begin{prop}\label{m-top}  
Les familles $M$-fid\`element plates d\'efinissent une 
pr\'etopologie sur $T$. 
\end{prop}  

\textit{Preuve:}
Il est imm\'ediat que tout isomorphisme est $M$-fid\`element plat. De plus, 
de part l'existence des isomorphismes naturels $(pq)^{*}\simeq q^{*}p^{*}$ on voit que
si $\{X_{i}\to X\}_{i}$ est une famille $M$-couvrante, et si 
pour tout $i\in I$, $\{Y_{i,j} \to X_{i}\}_{j}$ est une famille $M$-couvrante, alors 
la famille totale $\{Y_{i,j} \to X\}_{i,j}$ est encore $M$-couvrante.

Il nous reste \`a d\'emontrer la stabilit\'e des familles $M$-couvrantes
par changements de bases. Soit $\{p_{i} : X_{i}\to X\}_{i}$ une famille
$M$-couvrante et $f :Y\to X$ un morphisme. Notons
$\{q_{i} : Y_{i}:=Y\times_{X}X_{i} \to Y\}$ la famille obtenue par changement de base.
Pour montrer que cette famille est $M$-couvrante, nous consid\'erons le foncteur
$$\prod_{i}q_{i}^{*} : M(Y) \longrightarrow \prod_{i}M(Y_{i})$$
dont nous cherchons \`a prouver le caract\`ere conservatif. D'apr\`es l'hypoth\`ese
$(2)$ \ref{h0} faite sur $M$ il suffit pour cela de montrer que le foncteur
$$\prod_{i}(f_{i})_{*}q_{i}^{*} : M(Y) \longrightarrow \prod_{i}M(X_{i})$$
est conservatif (o\`u $f_{i} : Y_{i} \longrightarrow X_{i}$ est la seconde projection).
D'apr\`es \ref{h0} $(3)$ ce morphisme est isomorphe au morphisme
$$(\prod_{i}p_{i}^{*})f_{*} : M(Y) \longrightarrow \prod_{i}M(X_{i}).$$
Comme $f_ {*}$ est conservatif (d'apr\`es \ref{h0} $(2)$) et que
$\prod_{i}p_{i}^{*}$ est conservatif par hypoth\`ese, il en est de m\^eme du foncteur
compos\'e 
$(\prod_{i}p_{i}^{*})f_{*}$. Ce la termine la preuve de la proposition.
\hfill $\Box$ \\

La topologie d\'efinie par la pr\'etopologie pr\'ec\'edente est appel\'ee 
la topologie \emph{$M$-fid\`element plate} sur $T$.  \\

Le r\'esultat principal de cette section est que $M$ est toujours un 
champ pour la topologie 
$M$-fid\`element plate. Avant d'\'enoncer ce r\'esultat rappelons qu'ici 
$M$ est un pseudo foncteur en cat\'egories et non en groupo\"\i des
et que la notion de champ est la suivante. On peut consid\'erer
$M^{iso}$, le sous pseudo foncteur (non plein) de $M$ qui poss\`ede les m\^eme
objets et form\'e de tous les isomorphismes. Ce sous pseudo foncteur 
est un pseudo foncteur en groupo\"\i des. De plus, pour deux objets
$x$ et $y$ dans $M(X)$, on dispose d'un pr\'efaisceau 
de morphisme $\underline{Hom}(x,y)$: il s'agit du pr\'efaisceau sur $T/X$
qui \`a $u : Y\to X$ associe l'ensemble $Hom(u^{*}(x),u^{*}(y))$. 
Par d\'efinition, $M$ est un champ si les deux conditions suivantes sont 
satisfaites.

\begin{enumerate}
\item  Le pseudo foncteur en groupo\"\i des $M^{iso}$ est un champ 
au sens de \cite{lm}.
\item Pour tout $X\in T$ et tout $x,y$ des objets de $M(X)$, le
pr\'efaisceau $\underline{Hom}(x,y)$ est un faisceau sur $T/X$.  
\end{enumerate}

On peut aussi combiner ces deux conditions en une: pour un
recouvrement $\{U_{i}\to X\}_{i}$ dans $T$ on dispose d'une 
cat\'egorie de donn\'ees de descente $Desc(U/X,M)$. 
Rappelons que les objets de cette cat\'egorie sont les donn\'ees
de descentes $\{x_{i},\phi_{i,j}\}_{i,j}$, o\`u les $x_{i}$ sont des objets 
$x_{i}\in M(U_{i})$ et les $\phi_{i,j} : (x_{i})_{|U_{i,j}} \simeq (x_{j})_{|U_{i,j}}$
sont des isomorphismes dans $M(U_{i,j})$ satisfaisant \`a la condition
de cocycle usuelle $\phi_{j,k}\circ \phi_{i,j}=\phi_{i,k}$ 
dans  $M(U_{i,j,k})$. Un morphisme 
$\{x_{i},\phi_{i,j}\}_{i,j} \to \{y_{i},\psi_{i,j}\}_{i,j}$ entre de telles donn\'ees de descente
est une famille de morphismes $f_{i} : x_{i} \rightarrow y_{i}$ 
dans $M(U_{i})$, compatible  aux $\phi_{i,j}$ et aux
$\psi_{i,j}$ au sens o\`u $\psi_{i,j}f_{i}=f_{j}\phi_{i,j}$ dans $M(U_{i,j})$. 
Avec ces notations, les deux conditions pr\'ec\'edentes pour que $M$ soit un champ  
expriment alors que le foncteur
naturel
$$M(X) \longrightarrow Desc(U/X,M)$$
est une \'equivalence pour tout recouvrement $\{U_{i} \to X\}_{i}$. \\

\begin{thm}\label{mdesc}
Le pseudo foncteur $M$ est un champ pour la topologie
$M$-fid\`element plate. 
\end{thm}

\textit{Preuve:}
Pour commencer, le pseudo foncteur $M$ peut-\^etre strictifi\'e, par le
proc\'ed\'e standard qui consiste \`a lui associer le pr\'efaisceau
en cat\'egories qui envoie $X\in T$ sur la cat\'egorie des
pseudo transformations naturelles de $Hom(-,X)$ vers $M$ (voir par 
exemple \cite[\S 3.3]{hol}).
Les conditions \ref{h0} sont bien entendues pr\'eserv\'ees par 
ce proc\'ed\'e de strictification. Nous allons donc supposer que
$M$ est un pr\'efaisceau en cat\'egories.

On consid\`ere alors  
$Gpd(T)$ la cat\'egorie des pr\'efaisceaux en groupo\"\i des sur $T$, 
ainsi que $SPr(T)$ la cat\'egorie des pr\'efaisceaux simpliciaux sur $T$.
Rappelons qu'il existe une paire de foncteurs adjoints
$$\Pi_{1} : SEns \longrightarrow Gpd \qquad Gpd \longleftarrow SEns : N,$$
o\`u $N$ est le foncteur nerf et $\Pi_{1}$ est son adjoint 
\`a gauche qui envoie un ensemble simplicial sur son groupo\"\i de fondamental. Cette
adjonction identifie de plus $Gpd$ \`a une sous-cat\'egorie reflexive
de $SEns$. Pour un ensemble simplicial $K$ et une cat\'egorie
$C$, les foncteurs $\Pi_{1}(K) \to C$ sont en bijection naturelle
avec les donn\'ees suivantes:
\begin{enumerate}
\item une application $f : K_{0}\to Ob(C)$ de l'ensemble des $0$-simplexes
de $K$ vers l'ensemble des objets de $C$.  
\item une application $m : K_{1} \to Iso(C)$ de l'ensemble des $1$-simplexes
de $K$ vers l'ensemble des isomorphismes de $C$ telle que
la source de $m(k)$ soit $d_{0}(k)$ et le but de $m(k)$ soit $d_{1}(k)$ et avec 
$m(s_{0}(x))=id_{f(x)}$ pour $x\in K_{0}$.
\item on demande de plus que pour tout $l\in K_{2}$ on ait
$$m(d_{1}(l))=m(d_{0}(l)d_{2}(l)).$$ 
\end{enumerate}

En d'autres termes, $\Pi_{1}(K)$ est le groupo\"\i de librement
engendr\'e par le graphe unitaire $K_{1} \rightrightarrows K_{0}$
avec les relations induites par le morphisme
de bord $K_{2} \longrightarrow K_{1}\times K_{1}\times K_{1}$. 

En passant aux pr\'efaisceaux on trouve une paire de foncteurs adjoints
$$\Pi_{1} : SPr(T) \longrightarrow Gpd(T) \qquad Gpd(T) \longleftarrow SPr(T) : N.$$
Le foncteur $\Pi_{1}$ commute aux produits finis et envoie  l'ensemble simplicial 
$\Delta^{1}$ sur le groupo\"\i de ayant deux objets et un unique isomorphisme
entre eux. Cela implique en particulier que le foncteur $\Pi_{1}$ envoie
\'equivalences d'homotopie sur \'equivalences cat\'egoriques. En passant aux pr\'efaisceaux
on voit de plus que $\Pi_{1}$ envoie \'equivalence d'homotopie de pr\'efaisceaux
simpliciaux sur une \'equivalence de pr\'efaisceaux en groupo\"\i des
poss\'edant un quasi-inverse \emph{global} (ce d\'etail nous sera essentiel). \\

Soit maintenant $\{U_{i} \to X\}_{i}$ un 
recouvrement dans $T$. On consid\`ere le pr\'efaisceau
d'ensembles $U:=\coprod_{i}U_{i}$ (cette somme est prise
dans les pr\'efaisceaux et non dans $T$), qui est muni d'une augmentation 
naturelle $U \to X$. On d\'efinit $N(U/X) \in SPr(T)$ qui est le nerf du morphisme
$U \longrightarrow X$:
$$\begin{array}{cccl}
N(U/X)_* : & \Delta ^{op} & \longrightarrow & Pr(T) \\
& [n] & \longmapsto & \underbrace{U\times_X U\times_X \dots \times_X U}_{n+1\; fois},
\end{array}$$
o\`u les faces et d\'eg\'en\'erescences sont donn\'ees par les
diff\'erentes projections et diagonales. On pose
$$U/X:=\Pi_{1}N(U/X),$$
qui est encore muni d'une augmentation naturelle $p : U/X \to X$ dans $Gpd(T)$. 
Le pr\'efaisceau en groupo\"\i des $U/X$ est construit de sorte \`a ce que  
la cat\'egorie des morphismes de pr\'efaisceaux (au sens strict)
$\underline{Hom}(U/X,M)$, s'identifie naturellement \`a la cat\'egorie 
$Desc(U/X,M)$ des donn\'ees de descentes de $M$ pour le
recouvrement $\{U_{i} \to X\}_{i}$ (nous laissons le soin au lecteur
de v\'erifier cela en utilisant la description des morphismes
$\Pi_{1}(K)\to C$ donn\'ee plus haut).   Pour d\'emontrer le th\'eor\`eme
\ref{mdesc} il s'agit donc de montrer que le foncteur naturel
$$p^{*} : M(X)\simeq \underline{Hom}(X,M) \longrightarrow
\underline{Hom}(U/X,M)\simeq Desc(U/X,M)$$
est une \'equivalence. 

En utilisant l'hypoth\`ese \ref{h0} $(1)$ et $(2)$ on voit facilement que le foncteur
$$p^{*} : M(X) \longrightarrow Desc(U/X,M)$$
poss\`ede un adjoint \`a droite $p_{*}$.  Ce foncteur envoie une
donn\'ee de descente $\{x_{i},\phi_{i,j}\}$ sur l'objet
$$p_{*}(\{x_{i},\phi_{i,j}\})=Lim\left ( 
\prod_{i}(p_{i})_{*}(x_{i}) \rightrightarrows
\prod_{i,j}(p_{i,j})_{*}((x_{i})_{|U_{i,j}})
\right ).$$
Dans cette expression $p_{i,j} : U_{i,j}\to X$ d\'esigne la projection naturelle, et 
les deux morphismes proviennent d'une part des morphismes
naturels $(p_{i})_{*}(x_{i}) \to (p_{i,j})_{*}((x_{i})_{|U_{i,j}})\simeq (p_{i})_{*}((q_{i,j})_{*}q_{i,j}^{*}(x_{i}))$
(avec $q_{i,j} : U_{i,j} \longrightarrow U_{i}$ la premi\`ere projection), et d'autre
par le morphisme 
$$(p_{i})_{*}(x_{i}) \to (p_{i,j})_{*}((x_{j})_{|U_{i,j}})\simeq (p_{j})_{*}((q_{i,j})_{*}q_{i,j}^{*}(x_{i}))$$
obtenu en composant $(p_{i})_{*}(x_{i}) \to (p_{i,j})_{*}((x_{i})_{|U_{i,j}})$ avec 
l'isomorphisme $\phi_{i,j}$. 

Comme le foncteur
$p^*$ est conservatif (d'apr\`es la d\'efinition
des familles $M$-couvrantes) il suffit de montrer la co-unit\'e de l'adjonction
$p^*p_{*} \Rightarrow id$ est un isomorphisme. 
Pour cela, on consid\`ere le carr\'e cart\'esien dans $Gpd(T)$ suivant
$$\xymatrix{
Y \ar[d] \ar[r] & U/X \ar[d] \\
U \ar[r] & X.}$$
Comme $U=\coprod_{i}U_{i}$, on peut \'ecrire $Y=\coprod_{i}Y_{i}$ avec 
des diagrammes cart\'esiens
$$\xymatrix{
Y_{i} \ar[r] \ar[d] & Y \ar[d] \ar[r] & U/X \ar[d] \\
U_{i} \ar[r] & U \ar[r] & X.}$$
Par construction on voit que $Y_{i} \to U_{i}$ est lui-m\^eme isomorphe 
au morphisme
$q_{i} : V_{i}/U_{i} \to U_{i}$, o\`u $V_{i}$ d\'esigne le recouvrement
$\{U_{i,j} \to U_{i}\}_{i}$ obtenu \`a partir de $\{U_{i} \to X\}_{i}$ par changement 
de bases le long de $U_{i} \to X$. 
Ces carr\'es cart\'esiens induit des diagrammes commutatifs de cat\'egories
$$\xymatrix{
M(X) \ar[r] \ar[d]_-{p^{*}} & M(U_{i}) \ar[d]^-{q_{i}^{*}} \\
M(U/X) \ar[r] & M(V_{i}/U_{i}),}$$
(o\`u l'on note $M(F)$ pour $\underline{Hom}(F,M)$ pour tout
pr\'efaisceau en groupo\"\i des $F$).   Par d\'efinition de la topologie 
$M$-fid\`element plate les deux
familles de foncteurs 
$$M(X) \longrightarrow M(U_{i}) \qquad M(U/X) \longrightarrow M(V_{i}/U_{i})$$
sont conservatives. De plus, on dispose aussi 
de diagrammes commutatifs \`a isomorphisme pr\'es
$$\xymatrix{
M(X)   \ar[r] & M(U_{i})  \\
M(U/X)  \ar[u]^-{p_{*}} \ar[r] & \ar[u]_-{(q_{i})_{*}}M(V_{i}/U_{i})}$$
(on utilisera ici l'hypoth\`ese \ref{h0} $(3)$, la description
explicite des adjoints $p_{*}$ et $(q_{i})_ {*}$ en termes de limites finies, 
et la $M$-platitude des morphismes $U_{i}\to X$).
Ainsi, pour montrer que la co-unit\'e d'adjonction
$p^*p_{*} \Rightarrow id$ est un isomorphisme il suffit de montrer que pour tout
$i$ la co-unit\'e d'adjonction $q_{i}^*(q_{i})_{*} \Rightarrow id$ est un isomorphisme. 
En d'autres termes, on peut supposer qu'il existe un indice $i$ avec 
$X=U_{i}$ (et le morphisme $U_{i} \to X$ \'egal \`a l'identit\'e).

Nous supposons donc qu'il existe $i$ tel que $U_{i}=X$. Dans ce cas le morphisme
$U\to X$ poss\`ede une section, et il est bien connu que cette section induit
un inverse \`a homotopie pr\`es de la projection $N(U/X) \to X$.
Plus pr\'ecis\'ement, on construit une contraction de $N(U/X)$ sur $X$ 
$$c : \Delta^{1}\times N(U/X) \longrightarrow N(U/X)$$
de la fa\c{c}on suivante: pour $m\in \Delta^{op}$ le morphisme
$$c_{m} : \Delta^{1}(m)\times N(U/X)_{m} \longrightarrow N(U/X)_{m}$$
envoie symboliquement un morphisme $u : [m] \to [1]$ dans $\Delta$ et 
un point $(x_{0},\dots,x_{m})$ dans $N(U/X)_{m}$ sur le point
$(x_{0},\dots,x_{i-1},sp(x_{i}),\dots,sp(x_{m}))$, o\`u 
$i$ est tel que $u(j)=0 \Leftrightarrow j<i$. En appliquant le foncteur
$\Pi_{1}$ on voit que la projection $p : U/X \to X$ poss\`ede
un quasi-inverse global. Cela implique en particulier que pour tout
pr\'efaisceau en cat\'egories $M'$ le foncteur induit sur les cat\'egories
de morphismes
$$\underline{Hom}(X,M') \longrightarrow \underline{Hom}(U/X,M')$$
poss\`ede un quasi-inverse et est donc une \'equivalence. En prenant
$M'=M$ on trouve bien que le foncteur
$$M(X)\simeq \underline{Hom}(X,M') \longrightarrow 
\underline{Hom}(U/X,M)\simeq Desc(U/X,M)$$
est une \'equivalence de cat\'egories. Ceci termine la preuve du th\'eor\`eme
\ref{mdesc}. \hfill $\Box$ \\

\subsection{La topologie fid\`element plate}

Dans cette section nous allons utiliser le th\'eor\`eme
g\'en\'eral \ref{mdesc} pour construire la topologie 
fi\`element plate et quasi-compacte sur la cat\'egorie
oppos\'ee des mono\"\i des commutatifs dans une cat\'egorie
mono\"\i dale sym\'etrique quelconque. Cette topologie sera pour nous
auxiliaire et nous la rafinerons en d\'efinissant la topologie
de Zariski \`a la section suivante. C'est cette derni\`ere qui sera
utilis\'ee pour la d\'efinition de sch\'emas. \\

Tout au long de cette section $(C,\otimes,\mathbf{1})$ d\'esignera une
cat\'egorie mono\"\i dale
sym\'etrique avec $\mathbf{1}$ comme objet unit\'e. Nous supposerons de plus
que $C$ satisfait aux conditions suivantes. 

\begin{hyp}\label{h1}
\begin{enumerate}

\item La cat\'egorie $C$ poss\`ede des limites et des colimites. 

\item La structure mono\"\i dale $\otimes$ est ferm\'ee. En d'autres termes
pour toute paire d'objets $X$ et $Y$ dans $C$, le foncteur
$$\begin{array}{ccc}
C^{op} & \longrightarrow & Ens \\
Z & \mapsto & Hom(Z\otimes X,Y)
\end{array}$$
est repr\'esentable par un objet $\underline{Hom}(X,Y)\in C$.

\end{enumerate}
\end{hyp}

L'hypoth\`ese $(2)$ entraine que le produit tensoriel commute avec les colimites en 
chacune de ses variables. C'est ce que nous utiliserons implicitement tout au long de 
ce travail, et nous n'utiliserons pas les objets $\underline{Hom}(X,Y)$. \\

Nous noterons tout au long de ce chapitre
$Comm(C)$ la cat\'egorie des mono\"\i des associatifs, unitaires
et commutatifs dans $(C,\otimes,\mathbf{1})$. Nous noterons aussi
$Aff_{C}:=Comm(C)^{op}$ la cat\'egorie oppos\'ee de $Comm(C)$. 
Pour un objet $A\in Comm(C)$, nous noterons symboliquement
$Spec\, A$ l'objet correspondant dans $Aff_{C}$. 

\begin{df}\label{d1}
La \emph{cat\'egorie des sch\'emas affines sur $C$} est 
$Aff_{C}$. 
\end{df}

Comme $C$ poss\`ede tout type de limites et de colimites il en est 
de m\^eme de $Comm(C)$, et donc de $Aff_{C}$. De plus, le foncteur
d'oubli $Comm(C) \longrightarrow C$ poss\`ede un adjoint  \`a gauche
$$L : C \longrightarrow Comm(C)$$
qui envoit un objet $M$ de $C$ sur le mono\"\i de commutatif libre
$L(M)$ engendr\'e par $M$. \\

Pour $A\in Comm(C)$, on dispose d'une notion de $A$-module dans $C$, 
qui forme une cat\'egorie not\'ee $A-Mod$. Cette cat\'egorie poss\`ede aussi 
des limites et des colimites, et de plus le foncteur d'oubli $A-Mod \longrightarrow C$
y commute. Plus g\'en\'eralement, pour un morphisme $A \longrightarrow B$ dans
$Comm(C)$, on dispose d'une adjonction
$$-\otimes_{A}B : A-Mod \longrightarrow B-Mod \qquad 
A-Mod \longleftarrow B-Mod,$$
dont l'adjoint  $B-Mod \longrightarrow A-Mod$ est le foncteur
d'oubli \'evident. Ce foncteur poss\`ede aussi un adjoint \`a droite, 
et donc commute avec les limites et les colimites. Les morphismes
dans $A-Mod$ entre deux objets $M$ et $N$ 
seront not\'es $Hom_{A}(M,N)$.

Pour un diagramme de morphismes $\xymatrix{A' & \ar[l] A \ar[r]&  B}$
dans $Comm(C)$, il existe un isomorphisme naturel dans $B-Mod$
$$A'\coprod_{A}B \longrightarrow A'\otimes_{A}B,$$
o\`u $A'\coprod_{A}B$ est la somme calcul\'ee dans $Comm(C)$, et 
$A'\otimes_{A}B$ est le changement de base du $A$-module $A'$ 
par le morphisme $A \longrightarrow B$. Lorsque
$A=\mathbf{1}$ on trouve que la somme directe dans la cat\'egorie
$Comm(C)$ est donn\'ee par le produit tensoriel de mono\"\i des commutatifs. 

Enfin, pour $A\in Comm(C)$, la cat\'egorie $A-Mod$ est munie
d'une structure mono\"\i dale sym\'etrique $\otimes_{A}$, qui fait de
$A-Mod$ une cat\'egorie mono\"\i dale ferm\'ee. Les mono\"\i des
commutatifs dans $A-Mod$ forment une cat\'egorie \'equivalente
\`a la cat\'egorie $A/Comm(C)$, des objets de $Comm(C)$ en dessous
de $A$. Les objets de $A/Comm(C)$ seront appel\'es des
$A$-alg\`ebres commutatives. \\

Nous appliquons maintenant les r\'esultats de la section pr\'ec\'edente au cas
o\`u $T:=Aff_{C}$, 
et $M$ est le pseudo foncteur qui envoie $A\in Comm(C)$ 
sur la cat\'egorie $A-Mod$, et un morphisme $Spec\, A \to Spec\, B$ 
sur le foncteur $-\otimes_{A}B : A-Mod \to B-Mod$. Les conditions 
de l'hypoth\`ese \ref{h0} sont alors satisfaites. Noter que le point $(3)$ 
provient pr\'ecis\'ement du fait que la somme amalgamm\'ee
d'un diagramme $\xymatrix{B & A\ar[r] \ar[l] & C}$ dans $Comm(C)$
est $B\otimes_{A}C$. 

\begin{df}\label{dplat}
Avec les notations ci-dessus, la topologie $M$ -fid\`element plate
sur $Aff_{C}$ sera appel\'ee la \emph{topologie fid\`element plate
et quasi-compacte} (ou simplement \emph{fpqc}).  
\end{df}

On d\'eduite directement du th\`eor\`eme \ref{mdesc} que 
$Spec\, A \to A-Mod$ est un champ pour la topologie fpqc. 

\subsection{La topologie de Zariski}

Nous continuons avec une cat\'egorie mono\"\i dale sym\'etrique $(C,\otimes,\mathbf{1})$
qui satisfait aux conditions de la section pr\'ec\'edente.
Dans ce paragraphe nous allons d\'efinir la topologie de Zariski sur
la cat\'egorie $Aff_{C}$. Nous profiterons de l'occasion pour 
redonner de mani\`ere plus explicite la topologie fpqc d\'efinie dans la section 
pr\'ec\'edente.
Pour cela nous commen\c{c}ons par 
les d\'efinitions suivantes. 

\begin{df}\label{d2}
Soit $f : Y=Spec\, B \longrightarrow X=Spec\, A$ un morphisme
dans $Aff_{C}$.
\begin{enumerate}
\item Le morphisme $f$ 
est \emph{plat} si le foncteur
$$-\otimes_{A}B : A-Mod \longrightarrow B-Mod$$
est exact (i.e. commute aux limites finies).
\item Le morphisme $f$ est un \emph{\'epimorphisme} si 
pour tout $A' \in Comm(C)$, le morphisme
$$f^{*} : Hom(B,A') \longrightarrow Hom(A,A')$$
est injectif.
\item Le morphisme $f$ est \emph{de pr\'esentation finie}, 
si pour tout diagramme filtrant d'objets $A'_{i} \in A/Comm(C)$, 
le morphisme naturel
$$Colim_{i} Hom_{A/Comm(C)}(B,A'_{i}) \longrightarrow 
Hom_{A/Comm(C)}(B,Colim_{i}A'_{i})$$
est un isomorphisme. 
\item 
Le morphisme $f$
est \emph{un ouvert de Zariski} (ou encore
\emph{une immersion Zariski ouverte}) si le morphisme correspondant
$A \longrightarrow B$ dans $Comm(C)$
est un \'epimorphisme plat et de pr\'esentation finie. 

\end{enumerate}
\end{df}

A l'aide des d\'efinitions pr\'ec\'edentes nous d\'efinissons
les recouvrements fpqc et Zariski de la fa\c{c}on suivante. 

\begin{df}\label{topzar}
\begin{enumerate}
\item
Une famille de morphismes 
$$\{X_{i}=Spec\, A_{i} \longrightarrow X=Spec\, A\}_{i\in I}$$
dans $Aff_{C}$ est \emph{un recouvrement fpqc} (ou plus
simplement \emph{recouvrement plat}) si les deux 
conditions suivantes sont satisfaites. 
\begin{enumerate}
\item Pour tout $i\in I$ le morphisme $X_{i} \longrightarrow X$ est 
plat.
\item Il existe un sous ensemble fini $J\subset I$, tel que 
le foncteur
$$\prod_{j\in J} -\otimes_{A}A_{j} : A-Mod \longrightarrow 
\prod_{j\in J}A_{j}-Mod$$
est conservatif (i.e. un morphisme 
$u : M \rightarrow N$ de $A$-modules est 
un isomorphisme si et seulement si pour tout
$j\in J$ le morphisme induit $M\otimes_{A}A_{j} \rightarrow N\otimes_{A}A_{j}$
est un isomorphisme).
\end{enumerate}
\item Une famille de morphismes 
$$\{X_{i} \longrightarrow X\}_{i\in I}$$
dans $Aff_{C}$ est \emph{un recouvrement de Zariski} si
c'est un recouvrement plat et si tous les morphismes
$X_{i} \longrightarrow X$ sont des ouverts de Zariski.
\end{enumerate}
\end{df}

Il est facile de voir que 
les recouvrements fpqc et Zariski d\'efinissent deux 
pr\'etopologies sur la cat\'egorie $Aff_{C}$. 
Les topologies
de Grothendieck associ\'ees \`a ces deux pr\'etopologies seront 
appel\'ee respectivement la \emph{topologie fpqc} (ou encore
\emph{topologie plate}) et la
\emph{la topologie de Zariski}. 
Nous nous 
int\'eresserons principalement \`a la topologie de Zariski, et la 
topologie fpqc sera utilis\'ee que de mani\`ere auxiliaire (essentiellement
pour le corollaire \ref{p1} ci-dessous). 

On dispose donc
d'une cat\'egorie de pr\'efaisceaux (d'ensembles)
$Pr(Aff_{C})$, ainsi que deux sous-cat\'egories
de faisceaux
$$Sh^{fpqc}(Aff_{C}) \subset Sh^{Zar}(Aff_{C}) \subset Pr(Aff_{C}).$$
Comme la cat\'egorie qui nous int\'eressera principalement est 
$Sh^{Zar}(Aff_{C})$ nous noterons simplement
$$Sh(Aff_{C}):=Sh^{Zar}(Aff_{C}).$$
De m\^eme, l'expression \emph{faisceau}, sans plus de pr\'ecision, fera toujours
r\'ef\'erence \`a la notion de faisceau pour la topologie de Zariski. Les objets
de $Sh^{fpqc}(Aff_{C})$ seront eux appel\'es \emph{faisceaux fpqc}. \\

On dispose bien entendu du plongement de Yoneda
$$h_{-} : Aff_{C} \longrightarrow Pr(Aff_{C}),$$
et on identifiera toujours
$Aff_{C}$ avec la sous-cat\'egorie pleine
de $Pr(Aff_{C})$ form\'ee de l'image de $h_{-}$. 

\begin{cor}\label{p1}
\begin{enumerate}
\item Pour tout $X\in Aff_{C}$, le pr\'efaisceau $h_{X} \in Pr(Aff_{C})$
est un faisceau fpqc. Ce faisceau
sera simplement not\'e $X\in Sh^{fpqc}(Aff_{C})\subset Sh(Aff_{C})$. 
\item Le pr\'e-champ sur $Aff_{C}$, qui \`a un sch\'ema affine $X=Spec\, A$
associe la cat\'egorie $A-Mod$, et \`a un morphisme 
$Y=Spec\, B \longrightarrow X=Spec\, A$ associe le foncteur
$-\otimes_{A}B$ est un champ pour la topologie fpqc.
\end{enumerate}
\end{cor}

\textit{Preuve:} Le point $(2)$ est une cons\'equence du th\'eor\`eme \ref{mdesc}.
Pour $(1)$, il s'agit de montrer que pour tout
recouvrement fpqc $\{U_{i}=Spec\, B_{i} \longrightarrow X=Spec\, B\}_{I\in I}$, 
avec $I$ fini, et tout $B$-module $M\in B-Mod$, le diagramme 
$$M \longrightarrow \prod_{i}M\otimes_{B}B_{i} \rightrightarrows 
\prod_{i,j}M\otimes_{B}B_{i}\otimes_{B}B_{j}$$
est exact dans $C$. Cela se d\'eduit aussi du th\'eor\`eme
\ref{mdesc} en utilisant l'\'equivalence
$$M(X) \leftrightarrows Desc(U/X,M)$$
ainsi que la forme explicite de l'adjoint \`a droite $p_{*}$ 
donn\'ee dans la section $\S 2.1$.
\hfill $\Box$ \\

La proposition pr\'ec\'edente nous dit d'une part que la
topologie fqpc est sous-canonique, et d'autre part qu'elle
satisfait la condition descente pour les modules. Comme la
topologie de Zariski est moins fine que la topologie
plate nous en d\'eduisons que cela reste vrai pour
la topologie Zariski. 

A l'aide de la proposition pr\'ec\'edente nous identifierons 
la cat\'egorie $Aff_{C}$ avec son image dans $Sh^{fpqc}(Aff_{C})$, et donc
comme une sous-cat\'egorie pleine de $Sh(Aff_{C})$. Ainsi, 
un faisceau isomorphe \`a un objet de $Aff_{C}$ sera simplement
appel\'e un sch\'ema affine (au dessus de $C$ si l'on veut pr\'eciser). 

\subsection{Sch\'emas}

Dans ce paragraphe nous pr\'esentons la d\'efinition principale de ce
travail, \`a savoir celle de sch\'ema au dessus de $C$. Pour cela, nous
commencerons par introduire la notion d'ouvert et de recouvrement de Zariski 
dans $Sh(Aff_{C})$. Les sch\'emas seront d\'efinis comme les
faisceaux poss\'edant un recouvrement ouvert Zariski par 
des sch\'emas affines. Nous montrerons alors quelques propri\'et\'es de base
des sch\'emas (e.g. recollement, stabilit\'e par produits fibr\'es et r\'eunion disjointes). 

\begin{df}\label{d3}
\begin{enumerate}
\item 
Soit $X$ un sch\'ema affine et $F\subset X$ un sous-faisceau
de $X$. Nous dirons que $F$ est un \emph{ouvert de Zariski
de $X$} s'il existe une famille d'ouverts de Zariski 
$\{X_{i} \longrightarrow X\}_{i\in I}$ dans $Aff_{C}$ (au sens
de la d\'efinition \ref{d2}) tel que 
$F$ soit l'image du morphisme de faisceaux
$$\coprod_{i\in I} X_{i} \longrightarrow X.$$
\item Un morphisme $f : F \longrightarrow G$ dans $Sh(Aff_{C})$ est 
un \emph{ouvert de Zariski} (ou encore une \emph{immersion 
Zariski ouverte}) si pour tout sch\'ema affine $X$ et tout
morphisme $X \longrightarrow G$, le morphisme induit
$$F\times_{G}X \longrightarrow X$$
est un monomorphisme d'image un ouvert de Zariski de $X$. 
\end{enumerate}
\end{df}

On remarquera que l'on ne suppose pas l'ensemble d'indices
$I$ fini dans la d\'efinition pr\'ec\'edente. 
On v\'erifie ais\'emment que les ouverts de Zariski sont 
des monomorphismes dans $Sh(Aff_{C})$. 

\begin{lem}\label{l0'}
Les ouverts Zariski sont stables par changements de bases et 
composition dans $Sh(Aff_{C})$.
\end{lem}

\textit{Preuve:} Ceci est imm\'ediat. \hfill $\Box$ \\

Le lemme suivant montre
de plus que la notion d'ouvert Zariski pr\'ec\'edente est compatible
avec celle d\'efinie dans la d\'efinition \ref{d2}.

\begin{lem}\label{l0}
Soit $f : Z \longrightarrow Y$ un morphisme de sch\'emas affines. 
Le morphisme $f$ est un ouvert Zariski au sens de la d\'efinition 
\ref{d2} si et seulement s'il est un ouvert de Zariski au sens
de la d\'efinition \ref{d3} $(2)$. 
\end{lem}

\textit{Preuve:} Commen\c{c}ons par supposer que
$f$ soit un ouvert Zariski au sens de la d\'efinition \ref{d2}. Pour tout
sch\'ema affine $X$ et tout morphisme $X\longrightarrow Y$ on sait que
le morphisme induit $X\times_{Y}Z \longrightarrow X$ est encore 
un ouvert Zariski au sens de la d\'efinition \ref{d2}. Ceci montre que
l'on peut prendre la famille r\'eduite \`a un \'el\'ement 
$\{X\times_{Y}Z \longrightarrow X\}$ dans la d\'efinition \ref{d3} $(2)$, et donc
que $f$ est un ouvert de Zariski au sens de la d\'efinition 
\ref{d3} $(2)$. 

Inversement, si $f$ est un ouvert de Zariski au sens de la d\'efinition
\ref{d3} $(2)$. On prend $X=Y$ et $X \longrightarrow Y$ l'identit\'e. 
On voit alors que $f$ est un monomorphisme et que son image
est l'image d'un morphisme 
$$\coprod Y_{i} \longrightarrow Y$$
pour un recouvrement Zariski $\{Y_{i} \longrightarrow Y\}$ (que l'on 
peut supposer fini car $Z$ est affine et donc quasi-compact par d\'efinition
de la topologie de Zariski). Ainsi, 
$Y_{i}\simeq Y_{i}\times_{Y}X$, et $\{Y_{i} \longrightarrow X\}$
est donc un recouvrement Zariski tel que chaque morphisme induit
$Y_{i} \longrightarrow Y$ soit un ouvert Zariski. Tranduisant en termes
de monoides commutatifs, si $f : Z=Spec\, B \longrightarrow Y=Spec\, A$, 
il existe une famille de morphismes $\{B \longrightarrow B_{i}\}$ qui forme
un recouvrement Zariski et tel que chaque morphisme
$A \longrightarrow B_{i}$ soit un ouvert Zariski. On voit facilement \`a l'aide
des d\'efinitions que cela implique que $A \longrightarrow B$ est un morphisme
plat au sens de la d\'efinition \ref{d2}. Le morphisme $f$ \'etant un monomorphisme 
on trouve que le morphisme $A \longrightarrow B$
est donc un \'epimorphisme plat. Il nous reste \`a voir qu'il est 
de pr\'esentation finie. Pour cela, soit 
$B'_{\alpha} \in A/Comm(C)$ un diagramme filtrant d'objets dans $A/Comm(C)$, 
et montrons que
le morphisme naturel
$$Colim_{\alpha} Hom_{A/Comm(C)}(B,B'_{\alpha}) \longrightarrow 
Hom_{A/Comm(C)}(B,Colim_{\alpha}B'_{\alpha})$$
est bijectif. Comme $A \longrightarrow B$ est un \'epimorphisme on voit que
les ensembles source et but du morphisme pr\'ec\'edent sont
ou bien vides ou bien r\'eduits \`a un point. Il nous suffit donc de montrer que
si $Hom_{A/Comm(C)}(B,Colim_{\alpha}B'_{\alpha})$ est non vide alors il en 
est de m\^eme de $Colim_{\alpha} Hom_{A/Comm(C)}(B,B'_{\alpha})$, c'est \`a dire
qu'il existe un $\alpha_{0}$ tel que $Hom_{A/Comm(C)}(B,B'_{\alpha_{0}})$
soit non vide. 

Notons $B':=Colim_{\alpha}B'_{\alpha}$ et posons
$$B'_{i}:=B'\otimes_{B}B_{i} \qquad B'_{i,\alpha}:=B'_{\alpha}\otimes_{B}B_{i}.$$
On pose aussi pour deux indices $i$ et $j$
$$B_{ij}:=B_{i}\otimes_{B}B_{j} \qquad B'_{ij}:=B'\otimes_{B}B_{ij} \qquad
B'_{ij,\alpha}:=B'_{\alpha}\otimes_{B}B_{ij}.$$
On dispose de diagrammes dans $A/Comm(C)$ (index\'es par la cat\'egorie
$\rightrightarrows$)
$$\prod_{i} B_{i} \rightrightarrows \prod_{i,j}B_{ij} \qquad
\prod_{i} B'_{i} \rightrightarrows \prod_{i,j}B'_{ij},$$
et le morphisme $B \longrightarrow B'$ induit un morphisme du premier vers le second. 
Comme chaque $B_{i}$ et chaque $B_{ij}$ est de pr\'esentation finie
dans $A/Comm(C)$, ce morphisme de diagrammes se factorise par un morphisme
vers le diagramme 
$$\prod_{i} B'_{i,\alpha_{0}} \rightrightarrows \prod_{i,j}B'_{ij,\alpha_{0}}$$
pour un certain indice $\alpha_{0}$. En passant \`a la limite
le long de la cat\'egorie $\rightrightarrows$, et en appliquant le corollaire
\ref{p1} $(2)$ on trouve un 
morphisme $B \longrightarrow B_{\alpha_{0}}$ qui factorise
le morphisme $B \longrightarrow B'$ dans $A/Comm(C)$. Ceci termine la preuve du fait que
$f$ est un ouvert Zariski au sens de la d\'efinition \ref{d3} $(2)$. 
\hfill $\Box$ \\

Nous sommes maintenant pr\^ets pour d\'efinir la notion de sch\'ema
relatif. 

\begin{df}\label{sch}
Un faisceau $F\in Sh(Aff_{C})$ est un \emph{sch\'ema 
relatif \`a $C$} (ou simplement un \emph{sch\'ema} si le
contexte est clair) s'il existe une famille
de sch\'emas affines $X_{i}$ et un morphisme
$$p : \coprod_{i\in I} X_{i} \longrightarrow F$$
v\'erifiant les deux conditions suivantes. 
\begin{enumerate}
\item Le morphisme $p$ est un \'epimorphisme de faisceaux. 
\item Pour tout $i\in I$ le morphisme $X_{i} \longrightarrow F$
est une immersion Zariski ouverte. 
\end{enumerate}
Une famille de morphismes $\{X_{i} \longrightarrow F\}$ comme
ci-dessus est appel\'e un \emph{recouvrement Zariski affine de $F$}.

La \emph{cat\'egorie des sch\'emas relatifs \`a $C$}
est la sous-cat\'egorie pleine de $Sh(Aff_{C})$ form\'ee
des sch\'emas au sens ci-dessus. Nous la
noterons $Sch(C)$. 
\end{df}

Une propri\'et\'e fondamentale des sch\'emas relatifs est 
la propri\'et\'e de recollement suivante. 

\begin{prop}\label{p3}
L'application qui \`a $X\in Aff_{C}$ associe 
la cat\'egorie $Sch(C)/X$ d\'efinit un sous-champ  plein
du champ des faisceaux sur $Aff_{C}$ (qui \`a $X$ associe
$Sh(Aff_{C})/X$). 
\end{prop}

\textit{Preuve:}  Pour cela il nous faut v\'erifier les deux assertions
suivantes.
\begin{enumerate}
\item Soit 
$$\xymatrix{
F \ar[r] \ar[d] & G \ar[d] \\
Y \ar[r] & X}$$
un diagramme cart\'esien de faisceaux avec $X$ et $Y$ des sch\'emas affines. 
Si $G$ est un sch\'ema alors $F$ est un sch\'ema. 
\item Soit $X$ un sch\'ema affine et $F \longrightarrow X$
un morphisme de faisceaux. S'il existe un recouvrement Zariski affine
$\{X_{i} \longrightarrow X\}$ tel que $F\times_{X}X_{i}$
soit un sch\'ema pour tout $i$, alors $F$ est un sch\'ema. 
\end{enumerate}

Si 
$$\xymatrix{
F \ar[r] \ar[d] & G \ar[d] \\
Y \ar[r] & X}$$
est un carr\'e cart\'esien comme dans $(1)$ ci-dessus, 
et si $p : \coprod_{i\in I} X_{i} \longrightarrow G$ est 
un morphisme comme dans la d\'efinition \ref{sch}, 
les faisceaux $Y_{i}:=X_{i}\times_{G}F\simeq X_{i}\times_{X}Y$ sont des
sch\'emas affines, et le morphisme induit
$$\coprod_{i\in I} Y_{i} \longrightarrow F$$
est un \'epimorphisme tel que chaque morphisme
$Y_{i} \longrightarrow F$ soit un ouvert Zariski (voir lemme \ref{l0'}). 

Soit $F \longrightarrow X$ tel qu'au point $(2)$ ci-dessus. Pour tout
$i$, on choisit des sch\'emas affines $U_{ij}$ et un recouvrement Zariski affine
$$\coprod_{j} U_{ij} \longrightarrow F\times_{X}X_{i}.$$
On voit alors que 
le morphisme total
$$\coprod_{i,j} U_{ij} \longrightarrow F$$
est un \'epimorphisme. De plus, le lemme \ref{l0'} et la factorisation 
$$U_{ij} \longrightarrow  F\times_{X}X_{i} \longrightarrow F$$
montre que chaque $U_{ij} \longrightarrow F$ est un ouvert Zariski. 
Ceci montre que $F$ est un sch\'ema. \hfill $\Box$ \\

Notons aussi les deux faits suivants.

\begin{prop}\label{p3+}
\begin{enumerate}
\item Soit $F$ un sch\'ema et $F_{0} \subset F$ un ouvert de Zariski
au sens de la d\'efinition  \ref{d3}. Alors $F_{0}$ est un sch\'ema. 
\item Soit $f : F \longrightarrow G$ un morphisme entre sch\'emas. Alors $f$ 
est un ouvert de Zariski au sens de la d\'efinition \ref{d3} si et seulement si $f$ v\'erifie les
deux conditions suivantes.
\begin{enumerate}
\item Le morphisme $f$ est un monomorphisme.
\item  Il existe un recouvrement Zariski affine $\{X_{i} \longrightarrow F\}$ tel que
chaque morphisme $X_{i} \longrightarrow G$ soit un ouvert Zariski. 
\end{enumerate}
\end{enumerate}
\end{prop}

\textit{Preuve:} $(1)$ Soit $\{X_{i} \longrightarrow F\}$ un recouvrement 
Zariski affine. Pour chaque $i$, on pose $F_{0,i}:=F_{0}\times_{F}X_{i}$. 
Il existe donc une famille d'ouverts Zariski affine $\{U_{ij} \longrightarrow X_{i}\}$, tel que
$F_{0,i}$ soit l'image de $\coprod_{j}U_{ij} \longrightarrow X_{i}$. Les morphismes
$U_{ij} \longrightarrow F_{0}$ sont des ouverts Zariski, et de plus
$$\coprod_{i,j}U_{ij} \longrightarrow F_{0}$$
est un \'epimorphisme. Ceci implique que $F_{0}$ est un sch\'ema. \\

$(2)$ Supposons que $f$ soit un ouvert Zariski. Soit $\{Y_{i} \longrightarrow G\}$
un recouvrement Zariski affine de $G$. Pour tout $i$ soit 
$\{X_{ij} \longrightarrow Y_{i}\times_{G}F\}$ un recouvrement ouvert Zariski affine.
A l'aide du lemme \ref{l0'} on voit que 
la famille totale $\{X_{ij} \longrightarrow F\}$ est un recouvrement Zariski affine 
et chaque morphisme $X_{ij} \longrightarrow G$ est un ouvert Zariski.

Inversement, supposons que $f$ v\'erifie les deux conditions de la proposition. Soit $X$ un sch\'ema
affine et $X\longrightarrow G$ un morphisme.  On pose $F_{X}:=F\times_{G}X$, qui est 
un sous-faisceau de $X$. Chaqun des morphismes
$X_{i}\times_{G}X \longrightarrow X$ est un ouvert Zariski. Il existe donc 
pour tout $i$ un recouvrement Zariski affine $\{U_{ij} \longrightarrow X_{i}\times_{G}X\}$. 
Chaqun des morphismes $U_{ij} \longrightarrow X$ est un ouvert Zariski, et de plus
l'image de $\coprod_{i,j}U_{ij} \longrightarrow X$ est $F_{X}$. Ceci montre que
$f$ est un ouvert Zariski au sens de la d\'efinition \ref{d3}. \hfill $\Box$ \\

La proposition suivante donne certaines propri\'et\'es de stabilit\'e
des sch\'emas. 

\begin{prop}\label{p4-}
\begin{enumerate}
\item La sous-cat\'egorie $Sch(C)$ de
$Sh(Aff_{C})$ est stable par r\'eunions disjointes et
par produits fibr\'es. 
\item Un faisceau $F\in Sh(Aff_{C})$ est un sch\'ema si et seulement 
s'il existe une relation d'\'equivalence $R \subset X\times X$
dans $Sh(C)$ v\'erifiant les propri\'et\'es suivantes. 
\begin{enumerate}
\item On a 
$$X\simeq \coprod_{i\in I}U_{i}$$
avec $U_{i}$ des sch\'emas affines.  
\item Pour tout $(i,j)\in I^{2}$, consid\'erons le sous-faisceau
$R_{i,j}\subset U_{i}\times U_{j}$ d\'efini par le carr\'e cart\'esien
suivant
$$\xymatrix{
R_{i,j} \ar[r] \ar[d] & U_{i}\times U_{j} \ar[d] \\
R \ar[r] & X\times X,}$$
de sorte \`a ce que $R\simeq \coprod_{i,j}R_{i,j}$. 
Alors chacun des morphismes induits
$$R_{i, j} \longrightarrow U_{i}$$ 
est un ouvert Zariski.
\item Pour tout $i\in I$, le sous-objet $R_{i,i}\subset U_{i}\times U_{i}$ est \'egal \`a l'image
du morphisme diagonal $U_{i} \longrightarrow U_{i}\times U_{i}$. 
\item On a $F\simeq X/R$ (i.e. $F$ est isomorphe au faisceau quotient de $X$ par la relation
$R$). 
\end{enumerate}
\end{enumerate}
\end{prop}

\textit{Preuve:} $(1)$ Commen\c{c}ons par les sommes disjointes. 
Si $\{F_{i}\}$ est une famille de sch\'emas, et 
$$\coprod_{j} U_{ij} \longrightarrow F_{i}$$
un morphisme comme dans la d\'efinition \ref{sch}, 
le morphisme total 
$$\coprod_{i,j} U_{ij} \longrightarrow F$$
est un \'epimorphisme. Pour montrer que chaque 
$U_{ij} \longrightarrow F$ est un ouvert Zariski on utilise la factorisation
$$U_{ij} \longrightarrow F_{i} \longrightarrow F.$$
On est donc ramen\'es \`a montrer que $F_{i} \longrightarrow F$ 
est un ouvert Zariski. De fa\c{c}on plus g\'en\'erale on a le
lemme suivant. 

\begin{lem}
Soit $\{F_{i}\}_{i\in I}$ un famille de faisceaux. Pour tout $i$ le morphisme
$$F_{i} \longrightarrow F:=\coprod_{i\in I} F_{i}$$
est un ouvert Zariski. 
\end{lem} 

\textit{Preuve:} Soit $X$ un sch\'ema affine et $X \longrightarrow F$ un morphisme. 
Il existe un recouvrement Zariski $\{X_{j} \longrightarrow X\}_{j\in J}$, que l'on peut
supposer fini, et pour $j\in J$ un \'el\'ement $a(j)\in I$ et des diagrammes commutatifs
$$\xymatrix{
X_{j} \ar[r] \ar[d] & X \ar[d] \\
F_{a(j)} \ar[r] & F.}$$
Notons $J_{0}$ le sous-ensemble des $j\in J$ avec $a(j)=i$, alors
l'image du morphisme $F_{i} \times_{F}X \longrightarrow X$ est l'image
du morphisme 
$$\coprod_{j\in J_{0}} X_{j} \longrightarrow X.$$
Ceci montre que $F_{i} \longrightarrow F$ est un ouvert Zariski et finit la preuve du lemme. \hfill $\Box$ \\

Passons \`a la stabilit\'e de $Sch(C) \subset Sh(Aff_{C})$ par produits fibr\'es. 
Soit 
$$\xymatrix{ & F \ar[d] \\
G \ar[r] & H}$$
un diagramme dans $Sch(C)$, et consid\'erons $F\times_{H}G\in Sh(Aff_{C})$.

\begin{lem}
Soit $F_{1} \longrightarrow F_{0}$ un morphisme de faisceaux
avec $F_{0}$ un sch\'ema. S'il existe un recouvrement Zariski  affine
$\{X_{i} \longrightarrow F_{0}\}_{i\in I}$ tel que $F_{1}\times_{F_{0}}X_{i}$
soit un sch\'ema pour tout $i\in I$, alors $F_{1}$ est un sch\'ema. 
\end{lem}

\textit{Preuve:} C'est comme lors de la preuve de la proposition \ref{p3}. Si 
pour tout $i$ $\{X_{ij} \longrightarrow F_{1}\times_{F_{0}}X_{i}\}$
est un recouvrement Zariski affine, alors la famille totale
$\{X_{ij} \longrightarrow F_{1}\}$ est un recouvrement 
Zariski affine. \hfill $\Box$ \\

Le lemme pr\'ec\'edent appliqu\'e aux projections
$F\times_{H}G \longrightarrow F$ et $F\times_{H}G \longrightarrow G$ 
permet de se ramener au cas o\`u $F$ et $G$ sont des sch\'emas affines.
On peut alors trouver des recouvrements Zariski $\{X_{i} \longrightarrow X\}$
et $\{Y_{i} \longrightarrow Y\}$
et des factorisations
$$\xymatrix{
X_{i} \ar[r] \ar[d] & Z_{i} \ar[d] & Y_{i} \ar[l] \ar[d] \\ 
X \ar[r] & H& Y, \ar[l]}$$
o\`u chaque $Z_{i}$ est un sch\'ema affine et $Z_{i} \longrightarrow H$ est un
ouvert Zariski. Ainsi en appliquant une fois de plus le lemme aux projections
de $X\times_{H}Y$ sur $X$ et $Y$ on peut supposer qu'il existe
un sch\'ema affine $Z$, un ouvert Zariski $Z \longrightarrow H$ 
et des factorisations
$$\xymatrix{X \ar[rd] \ar[r] & Z \ar[d] & Y \ar[ld] \ar[l] \\
 & H. &}$$
Mais comme le morphisme $Z \longrightarrow H$ est un monomorphisme on 
a $X\times_{H}Y\simeq X\times_{Z}Y$ qui est donc un sch\'ema affine. \\

$(2)$ La n\'ecessit\'e se voit facilement en consid\'erant la relation d'\'equivalence 
induite par l'\'epimorphisme
$$\coprod_{i} X_{i} \longrightarrow F$$
o\`u $\{X_{i} \longrightarrow F\}$ est un recouvrement Zariski affine. 

Supposons maintenant qu'un faisceau $F$ s'\'ecrive comme $X/R$ avec 
$X$ et $R$ comme dans l'\'enonc\'e de la proposition et montrons que
$F$ est un sch\'ema. Le fait que chaque morphisme $U_{i} \longrightarrow F$ soit 
un monomorphisme se d\'eduit de la condition $(c)$ sur les $R_{i,i}$. 
Comme le morphisme 
$$X=\coprod U_{i} \longrightarrow F=X/R$$
est un \'epimorphisme il nous reste \`a voir que 
chaque morphisme $U_{i} \longrightarrow F$ est un ouvert Zariski. 

Soit $Y$ un sch\'ema affine, $Y \longrightarrow F$ un morphisme et 
posons $Y_{i}:=Y\times_{F}U_{i}$. On consid\`ere le morphisme
$Y_{i} \longrightarrow Y$, et on cherche \`a montrer que c'est
un ouvert Zariski. Il existe un recouvrement Zariski affine 
$\{Z_{j} \longrightarrow Y\}_{j\in J}$,  une application $a : J \longrightarrow I$ et 
des factorisations
$$\xymatrix{
Z_{j} \ar[r] \ar[d] & Y \ar[d] \\
U_{a(j)} \ar[r] & F.}$$
On a donc, pour tout $j\in J$ 
$$Y_{i}\times_{Y}Z_{j}\simeq R_{i,a(j)}\times_{U_{a(j)}}Z_{j}.$$ 
Comme $R_{i,a(j)} \longrightarrow U_{a(j)}$ est un ouvert Zariski, on voit que 
le morphisme
$$Y_{i}\times_{Y}Z_{j} \longrightarrow Y$$ 
est un ouvert Zariski pour tout $j\in J$. Il existe donc 
une famille d'ouverts Zariski affine $\{W_{ijk} \longrightarrow Y\}_{k\in K}$ telle que
l'image du morphisme
$$\coprod_{k} W_{ijk} \longrightarrow Y$$
soit \'egale \`a $Y_{i}\times_{Y}Z_{j}$. On voit alors que la famille
d'ouverts Zariski $\{W_{ijk} \longrightarrow Y\}_{(j,k)\in J\times K}$
est telle que l'image du morphisme
$$\coprod_{j,k}W_{ijk} \longrightarrow Y$$
est \'egale \`a $Y_{i}$. Ceci montre que $Y_{i} \longrightarrow Y$ est un
ouvert Zariski et donc que $U_{i} \longrightarrow F$ est un ouvert Zariski. 
Ceci termine la preuve du fait que $F$ est un sch\'ema. \hfill $\Box$ \\

Pour terminer ce paragraphe signalons un autre point de vue sur les sch\'emas
relatifs, que nous n'utiliserons pas, mais qui est plus proche 
de la notion usuelle de sch\'ema en tant qu'espace annel\'e. 

Soit $X$ un sch\'ema au dessus de $C$. On d\'efinit la  
cat\'egorie $Zar(X)$ des ouverts de Zariski de $X$ comme \'etant la sous cat\'egorie
pleine de $Sh(Aff_{C})/X$ form\'ee des $u : Y\rightarrow X$, avec $Y$ un sch\'ema et 
$u$ une immersion Zariski ouverte. La cat\'egorie $Zar(X)$ est un \emph{lieu} (\emph{locale} en 
anglais, voir \cite{mm}), c'est \`a dire que c'est la cat\'egorie sous-jacente 
\`a un ensemble partiellement ordonn\'e qui poss\`ede des sup ainsi que des
inf finis, et tel que les inf se distribuent sur les sup. En d'autres termes, 
$Zar(X)$ est une cat\'egorie qui se comporte comme la cat\'egorie des ouverts
d'un espace topologique. On dispose d'une topologie naturelle induite sur 
$Zar(X)$, qui n'est autre que la restriction de la topologie canonique de $Sh(Aff_{C})$. 
Ainsi, une famille de morphismes $\{Y_{i} \rightarrow Y\}$ est couvrante
dans $Zar(X)$ si le morphisme $\coprod Y_{i} \longrightarrow Y$ est un
\'epimorphisme de faisceaux. 
La sous-cat\'egorie pleine de $Zar(X)$ form\'ee des $Y\rightarrow X$ avec $Y$ un sch\'ema
affine sera not\'ee $ZarAff(X)$. Elle est aussi munie de la restriction de la topologie
canonique de $Sh(Aff_{C})$. De plus, le foncteur d'inclusion 
$ZarAff(X) \longrightarrow Zar(X)$ est continu, et induit une 
\'equivalence sur les cat\'egories
de faisceaux 
$$Sh(Zar(X)) \simeq Sh(ZarAff(X)).$$
Par cette \'equivalence, nous identifierons souvent implicitement 
ces deux cat\'egories. Nous la noterons simplement $Sh(X_{Zar})$.

Nous savons que le site $Zar(X)$ est un lieu, et que sa topologie
est engendr\'ee par une pr\'e-topologie quasi-compacte (i.e. les
familles couvrantes sont finies), \`a savoir
le site $ZarAff(X)$. Ceci implique que $Zar(X)$ est 
naturellement \'equivalent au lieu $Ouv(|X|)$ des ouverts d'un espace
topologique $|X|$ (voir \cite[Cor. IX.3.4, Cor IX.11.3]{mm}). 
L'espace topologique $|X|$ est \'evidemment
tel que le topos des faisceaux sur $|X|$ est \'equivalent au topos
des faisceaux sur $Zar(X)$. Ainsi, nous verrons tout faisceau
sur $Zar(X)$ comme un faisceau sur $|X|$. 
Ainsi, $Sh(X_{Zar})$ est  \'equivalente \`a $Sh(|X|)$, la cat\'egorie des faisceaux
sur l'espace topologique $|X|$. 

Soit $Y=Spec\, A \longrightarrow X$ un objet de $ZarAff(X)$. On lui associe l'objet
correspondant $A\in Comm(C)$. Lorsque $Y$ varie dans $ZarAff(X)$ ceci d\'efinit un
foncteur
$$\mathcal{O}_{X} : ZarAff(X)^{op} \longrightarrow Comm(C)$$
que l'on voit \^etre un faisceau par le corollaire \ref{p1}. On dispose ainsi d'un faisceau
$\mathcal{O}_{X}$ sur le topos $Sh(X_{Zar})$ \`a valeurs dans $Comm(C)$, ou 
de mani\`ere \'equivalente d'un espace $|X|$ muni d'un faisceau
$\mathcal{O}_{X}$ \`a valeurs dans $Comm(C)$. 
Le couple
$(|X|,\mathcal{O}_{X})$  joue dans notre situation le r\^ole 
des espaces annel\'es pour les sch\'emas au sens usuel. On pourrait d\'evelopper la 
th\'eorie de sch\'emas relatifs de ce point de vue, en d\'efinissant les sch\'emas au dessus
de $C$ somme des espaces munis de faisceaux \`a valeurs dans $Comm(C)$ et qui 
localement sont \'equivalents \`a un mod\`ele affine. Ceci donne une version un peu plus
g\'eom\'etrique des sch\'emas relatifs, mais ceci dit \'equivalente \`a celle que nous donnons
dans la d\'efinition \ref{sch}. Il nous semble que ce point de vue n'apporte en r\'ealit\'e pas
grand chose, d'autant plus que le point de vue fonctoriel nous serait indispensable par
exemple pour pouvoir consid\'erer les champs g\'eom\'etriques (notion que nous
n'aborderons cependant pas dans ce travail). 

\subsection{Changements de bases}

Supposons que $(C,\otimes,\mathbf{1})$ et $(D,\otimes,\mathbf{1})$ 
soient deux cat\'egories mono\"\i dales sym\'etriques qui v\'erifient les conditions
expos\'ees en d\'ebut de cette section. On se donne un foncteur mono\"\i dal sym\'etrique
(unitaire et associatif)
$$f : C \longrightarrow D$$
et l'on suppose qu'il poss\`ede un adjoint \`a droite
$$g : D \longrightarrow C.$$
Cet adjoint \`a droite n'est plus mono\"\i dal, mais on dispose 
de morphismes fonctoriels $g(X)\otimes g(Y) \longrightarrow g(X\otimes Y)$, qui sont
associatifs, commutatifs et unitaires. Cela suffit pour que $g$ induise
des foncteurs sur les cat\'egories de mono\"\i des commutatifs et 
de modules.  

Le foncteur $f$ induit un foncteur sur les cat\'egories des mono\"\i des commutatifs, et donc
sur les cat\'egories des sch\'emas affines
$$f : Aff_{C} \longrightarrow Aff_{D}.$$
Ce foncteur poss\`ede un adjoint \`a gauche
$$g : Aff_{D} \longrightarrow Aff_{C}$$
induit par la foncteur $g$ (nous garderons les m\^emes
notations pour les foncteurs induits sur les sch\'emas
affines).  

On dispose ainsi d'une adjonction sur les cat\'egories de pr\'efaisceaux
$$g_{!} : Pr(Aff_{D}) \longrightarrow Pr(Aff_{C}) \qquad
Pr(Aff_{D}) \longleftarrow Pr(Aff_{C}) : g^{*}=f_{!}.$$
Explicitement, pour $F\in Pr(Aff_{C})$, et 
$X \in Aff_{D}$ on a
$f_{!}(F)(X)=F(g(X))$. 

Rappelons qu'on dit que le foncteur $f : Aff_{C} \longrightarrow Aff_{D}$
est continu pour la topologie Zariski (resp. fpqc) si le foncteur
$f^{*} : Pr(Aff_{D}) \longrightarrow Pr(Aff_{C})$ pr\'eserve
les sous-cat\'egories de faisceaux Zariski  (resp. fpqc). Dans ce cas, 
l'adjonction $(f_{!},f^{*})$ descend en une adjonction sur les cat\'egories
de faisceaux
$$f^{\sim}_{!} : Sh(Aff_{C}) \longrightarrow Sh(Aff_{D}) \qquad
Sh(Aff_{C}) \longleftarrow Sh(Aff_{D}) : f^{*}.$$
Explicitement, le foncteur $f^{\sim}_{!}=a\circ f_{!} $ est le compos\'e
de $f_{!}$ sur les pr\'efaisceaux suivi du foncteur 
faisceau associ\'e. Comme le foncteur $f$ commute aux limites finies, 
il en est de m\^eme de $f^{\sim}_{!}$, et ainsi
l'adjonction 
$$f^{\sim}_{!} : Sh(Aff_{C}) \longrightarrow Sh(Aff_{D}) \qquad
Sh(Aff_{C}) \longleftarrow Sh(Aff_{D}) : f^{*}$$
d\'efinit un morphisme g\'eom\'etrique de topos. 

\begin{prop}\label{p4}
On garde les notations ci-dessus. 
\begin{enumerate}
\item Le foncteur $f : Aff_{C} \longrightarrow Aff_{D}$
pr\'eserve les monomorphismes et les limites finies. 
\item Si le foncteur $f : Aff_{C} \longrightarrow Aff_{D}$ 
est continu pour la topologie plate, et si de plus $g : D \longrightarrow C$ 
commute avec les colimites filtrantes, alors le foncteur
$f : Aff_{C} \longrightarrow Aff_{D}$ est aussi continue pour
la topologie de Zariski.
\end{enumerate}
\end{prop}

\textit{Preuve:} $(1)$ Le foncteur $f : C \longrightarrow D$ pr\'eserve
toutes les colimites car il poss\`ede un adjoint \`a droite. Ainsi, il en est 
de m\^eme pour $f : Comm(C) \longrightarrow Comm(D)$, et 
donc le foncteur induit sur les cat\'egories oppos\'ees pr\'eserve
les limites et en particulier les limites finies. De ceci nous d\'eduisons aussi
que $f : Aff_{C} \longrightarrow Aff_{D}$ pr\'eserve les monomorphismes, car 
$X \longrightarrow Y$ est un monomorphisme si et seulement si
le morphisme diagonal $X\longrightarrow X\times_{Y}X$ est un isomorphisme. \\

$(2)$ Par le point $(1)$ et l'hypoth\`ese on sait d\'ej\`a que 
$f : Comm(C) \longrightarrow Comm(D)$ pr\'eserve les \'epimorphismes
plats. Il nous reste donc \`a voir qu'il pr\'eserve aussi les
morphismes de pr\'esentation finie. Mais ceci se voit imm\'ediatement
par adjonction.  \hfill $\Box$ \\

\begin{cor}\label{c1}
On suppose que le foncteur $g : D \longrightarrow C$ est conservatif et qu'il
commute aux colimites filtrantes. On suppose de plus que
pour tout morphisme plat $A \longrightarrow B$ dans
$Comm(C)$, et tout $N\in f(A)-Mod$, le morphisme naturel
$$g(N)\otimes_{A}B \longrightarrow g(N\otimes_{f(A)}f(B))$$
est un isomorphisme dans $B-Mod$. 
Alors, 
$f : Aff_{C} \longrightarrow Aff_{D}$ est continu pour la topologie
de Zariski, et le foncteur
$$f^{\sim}_{!} : Sh(Aff_{C}) \longrightarrow Sh(Aff_{D})$$
pr\'eserve les sous-cat\'egories des sch\'emas et induit un foncteur
$$\begin{array}{ccc}
Sch(C) & \longrightarrow & Sch(D) \\
 X & \mapsto & X\times_{C}D:=f^{\sim}_{!}(X).
\end{array}$$
Pour tout sch\'ema affine $X$ au dessus de $C$  
$f^{\sim}_{!}(X)$ est isomorphe au sch\'ema affine $f(X) \in Aff_{D}$. 
\end{cor}

\textit{Preuve:} Commen\c{c}ons par montrer que
le foncteur $f : Aff_{C} \longrightarrow Aff_{D}$ est continu pour
la topologie plate. Pour cela, montrons tout d'abord que pour
un morphisme $A \longrightarrow B$ plat dans 
$Comm(C)$ le morphisme induit
$f(A) \longrightarrow f(B)$ est plat dans $Comm(D)$. Par hypoth\`ese, 
on dispose d'un diagramme commutatif (\`a isomorphisme pr\`es)
$$\xymatrix{
f(A)-Mod \ar[r]^-{g} \ar[d] & A-Mod \ar[d] \\
f(B)-Mod \ar[r]_-{g} & B-Mod,}$$
o\`u les foncteurs verticaux sont les changements
de bases. Comme le foncteur $g$ est conservatif et que
$A \longrightarrow B$ est plat, on voit facilement que cela
implique que $f(A) \longrightarrow f(B)$ est plat. Le m\^eme
argument montre que si 
$\{X_{i} \longrightarrow X\}$ est un recouvrement plat 
dans $Aff_{C}$ alors 
$\{f(X_{i}) \longrightarrow f(X)\}$ est un recouvrement plat 
dans $Aff_{D}$. Ceci fini de montrer que 
$f : Aff_{C} \longrightarrow Aff_{D}$ est continu
pour la topologie plate. Par la proposition \ref{p4} $(2)$ 
on voit que $f$ est continu pour la topologie de Zariski. 

Soit maintenant $X$ un sch\'ema au dessus de $C$. On sait 
d'apr\`es la proposition \ref{p4-} que l'on peut l'\'ecrire
comme un quotient $Y/R$, o\`u $Y$ est une r\'eunion disjointe
de sch\'emas affines et $R \subset Y\times Y$ une relation 
d'\'equivalence v\'erifiant les conditions de la proposition \ref{p4}. 
Comme
$f^{\sim}_{!}$ est un adjoint \`a gauche on trouve
que $f^{\sim}_{!}(X)$ est isomorphe au quotient 
$f^{\sim}_{!}(Y)/f^{\sim}_{!}(R)$. Comme le foncteur 
$f^{\sim}_{!}$ pr\'eserve clairement les sch\'emas affines, 
on a bien que $f^{\sim}_{!}(Y)$ est une r\'eunion disjointe 
de sch\'emas affines. De plus, comme $f ^{\sim}_{!}$ pr\'eserve les
\'epimorphismes (car il commute aux colimites et aux limites finies) et les
ouverts Zariski (voir proprosition \ref{p4}), on voit que la relation d\'equivalence 
$f^{\sim}_{!}(R) \longrightarrow f^{\sim}_{!}(Y)\times f^{\sim}_{!}(Y)$
v\'erifie les conditions de la proposition \ref{p4-} et donc
que $f^{\sim}_{!}(X)$ est un sch\'ema. 
\hfill $\Box$ \\

Pour un mono\"\i de commutatif  $A \in Comm(C)$, on a un isomorphisme naturel
$$f^{\sim}_{!}(Spec\, A)\simeq Spec\, f(A).$$
En effet, par d\'efinition le foncteur 
$f^{\sim}_{!}$ est le compos\'e de $f_{!}$, d\'efini sur les pr\'efaisceaux, et 
du foncteur faisceau associ\'e. Il est facile de voir que 
$f_{!}(Spec\, A)\simeq Spec\, f(A)$, en tant que pr\'efaisceau. Comme
les pr\'efaisceaux repr\'esentables par des sch\'emas affines sont 
des faisceaux, ceci montre bien que 
$f^{\sim}_{!}(Spec\, A)\simeq Spec\, f(A)$. De ceci on tire que
pour tout sch\'ema affine $X=Spec\, A\in Aff_{C}$, le foncteur
$f^{\sim}_{!}(X)$ est donn\'e par
$$\begin{array}{ccc}
Comm(C) & \longrightarrow & Ens \\
B & \mapsto & Hom(A,g(B)). 
\end{array}$$

\section{Trois exemples de g\'eom\'etries relatives}

Dans cette section nous pr\'esentons nos trois premiers exemples
de cat\'egories de sch\'emas relatifs. 

\subsection{$Spec\, \mathbb{Z}$}

Posons $(C,\otimes,\mathbf{1})=(\mathbb{Z}-Mod,\otimes,\mathbb{Z})$, 
la cat\'egorie mono\"\i dale sym\'etrique des groupes ab\'eliens. 
Elle v\'erifie les conditions de la section pr\'ec\'edente, donc
on dispose d'une cat\'egorie $Sch(\mathbb{Z}-Mod)$ des
sch\'emas au dessus de $\mathbb{Z}-Mod$. On voit facilement en 
d\'ecanulant les d\'efinitions que les notions de topologie
plate et d'ouverts de Zariski coincident avec les notions
usuelles (il suffit de remarquer qu'un morphisme
de sch\'emas est une immersion ouverte si et seulement si c'est un
monomorphisme plat et localement de pr\'esentation finie). Ainsi, 
on voit facilement que
la cat\'egorie $Sch(\mathbb{Z}-Mod)$ est naturellement \'equivalente
\`a la cat\'egorie des sch\'emas au sens usuel. 

\begin{df}\label{d6}
La \emph{cat\'egorie des $\mathbb{Z}$-sch\'emas} est 
$Sch(\mathbb{Z}-Mod)$. Elle sera not\'ee
$\mathbb{Z}-Sch$. 
\end{df}

\subsection{$Spec\, \mathbb{N}$}

Posons maintenant $(C,\otimes,\mathbf{1}):=(\mathbb{N}-Mod,\otimes,\mathbb{N})$, 
la cat\'egorie des mono\"\i des commutatifs (unitaires et associatifs), munie
de son produit tensoriel usuel. Nous appellerons les objets
de $\mathbb{N}-Mod$ aussi des $\mathbb{N}$-modules, pour des
raisons \'evidentes. Les conditions de la section pr\'ec\'edente sont
bien entendues v\'erifi\'ees, et on dispose donc d'une 
cat\'egorie $Sch(\mathbb{N}-Mod)$ des sch\'emas
relatifs \`a $\mathbb{N}-Mod$. 

\begin{df}\label{d7}
La \emph{cat\'egorie des $\mathbb{N}$-sch\'emas} est 
$Sch(\mathbb{N}-Mod)$. Elle sera not\'ee
$\mathbb{N}-Sch$. 
\end{df}

Par d\'efinition la cat\'egorie des $\mathbb{N}$-sch\'emas affines
est \'equivalente \`a la cat\'egorie oppos\'ee des 
semi-anneaux commutatifs. \\

On dispose d'un adjonction
$$f : \mathbb{N}-Mod \longrightarrow \mathbb{Z}-Mod \qquad 
\mathbb{N}-Mod \longleftarrow \mathbb{Z}-Mod : g,$$
o\`u $g$ est le foncteur d'inclusion des groupes
ab\'eliens dans les mono\"\i des ab\'eliens, et 
$f$ est le foncteur de compl\'etion des mono\"\i des
ab\'eliens vers les groupes ab\'eliens. Le foncteur
$f$ est mono\"\i dal sym\'etrique. 

\begin{prop}\label{p5}
Le foncteur
$$f : \mathbb{N}-Mod \longrightarrow \mathbb{Z}-Mod$$
v\'erifie les conditions du corollaire \ref{c1}. 
\end{prop}

\textit{Preuve:} Le fait que $g$ soit conservatif et commute aux colimites
filtrantes est clair. 
Pour un semi-anneau commutatif $A$, 
le foncteur $g : f(A)-Mod \longrightarrow A-Mod$
est pleinement fid\`ele, et son image consiste en les $A$-modules
$N$ tel que le mono\"\i de sous-jacent \`a $N$ est un groupe. Il est alors imm\'ediat 
de constater que pour un morphisme quelconque de 
semi-anneaux $A \longrightarrow B$, et tout 
$f(A)$-module $N$, le morphisme naturel
$$g(N)\otimes_{A}B \longrightarrow g(N\otimes_{f(A)}f(B))$$
est un isomorphisme (ceci est \'equivalent au fait que 
le mono\"\i de sous-jacent \`a $g(N)\otimes_{A}B$ est un groupe). 
\hfill $\Box$ \\

D'apr\`es le corollaire \ref{c1} on trouve donc un foncteur de changement 
de bases
$$-\otimes_{\mathbb{N}}\mathbb{Z} : 
\mathbb{N}-Sch \longrightarrow \mathbb{Z}-Sch.$$

\begin{prop}\label{p6}
Le foncteur
$$-\otimes_{\mathbb{N}}\mathbb{Z} : 
\mathbb{N}-Sch \longrightarrow \mathbb{Z}-Sch$$
poss\`ede un adjoint \`a gauche qui est pleinement fid\`ele. 
\end{prop}

\textit{Preuve:} Revenons au foncteur de compl\'etion
$f : \mathbb{N}-Mod \longrightarrow \mathbb {Z}-Mod$, et \`a son
adjoint \`a droite le foncteur d'inclusion
$g : \mathbb{Z}-Mod \longrightarrow \mathbb {N}-Mod$.
Ils induisent une adjonction sur les cat\'egories
de sch\'emas affines
$$g : Aff_{\mathbb{Z}} \longrightarrow Aff_{\mathbb{N}}\qquad
Aff_{\mathbb{Z}} \longleftarrow Aff_{\mathbb{N}} : f$$
o\`u $f$ est maintenant l'adjoint \`a droite. On sait que $f$ est continu pour
la topologie Zariski, et que le foncteur induit
$$f^{\sim}_{!} : Sh(Aff_{\mathbb{N}}) \longrightarrow Sh(Aff_{\mathbb{Z}})$$
pr\'eserve les sous-cat\'egories de sch\'emas. Il est facile de voir
que $g$ est aussi un foncteur continu pour la topologie Zariski, et que l'on a donc
$$f^{\sim}_{!}\simeq g^{*} : Sh(Aff_{\mathbb{N}}) \longrightarrow Sh(Aff_{\mathbb{Z}}).$$
Le foncteur 
$$g^{\sim}_{!} : Sh(Aff_{\mathbb{Z}}) \longrightarrow Sh(Aff_{\mathbb{N}})$$
est alors un adjoint \`a gauche de $f^{\sim}_{!}\simeq g^{*}$. De plus, 
$g$ \'etant pleinement fid\`ele on voit facilement que 
$g^{\sim}_{!}$ est aussi pleinement fid\`ele. Il nous reste donc \`a montrer que
$g^{\sim}_{!}$ pr\'eserve aussi les sch\'emas. Pour cela il suffit d'utiliser, 
comme pour la preuve du corollaire \ref{c1}, qu'il pr\'eserve
la notion d'ouverts Zariski.   \hfill $\Box$ \\

D'apr\`es la proposition \ref{p6} nous verrons
la cat\'egorie des $\mathbb{Z}$-sch\'emas plong\'ee
dans celle des $\mathbb{N}$-sch\'emas. Il faut comprendre la proposition 
\ref{p6} \`a travers le fait que le morphisme
$$Spec\, \mathbb{Z} \longrightarrow Spec\, \mathbb{N}$$
est un monomorphisme.

\subsection{$Spec\, \mathbb{F}_{1}$}

Nous posons maintenant $(C,\otimes,\mathbf{1})=(Ens,\times,*)$, la cat\'egorie
des ensembles munie de sa structure mono\"\i dale sym\'etrique
donn\'ee par le produit direct. Cette cat\'egorie v\'erifie les
conditions de notre premier paragraphe, et on dispose donc
d'une cat\'egorie $Sch(Ens)$ de sch\'emas au dessus de $Ens$. 
Nous appellerons les objets de $Sch(Ens)$ des $\mathbb{F}_{1}$-sch\'emas. 

\begin{df}\label{d8}
La \emph{cat\'egorie des $\mathbb{F}_{1}$-sch\'emas} est 
$Sch(Ens)$. Elle sera not\'ee
$\mathbb{F}_{1}-Sch$. 
\end{df}

Par construction, la cat\'egorie $Aff_{\mathbb{F}_{1}}$ des 
$\mathbb{F}_{1}$-sch\'emas affines est \'equivalente \`a la cat\'egorie
oppos\'ee de la cat\'egorie des mono\"\i des commutatifs (associatifs et unitaires). \\

On dispose d'une adjonction 
$$f : \mathbb{F}_{1}-Mod=Ens \longrightarrow \mathbb{N}-Mod \qquad
\mathbb{F}_{1}-Mod=Ens \longleftarrow \mathbb{N}-Mod : g,$$
o\`u $g$ est le foncteur d'oubli qui envoit un mono\"\i de commutatif
sur son ensemble sous-jacent. Le foncteur $f$ envoit un ensemble
$X$ sur $X\otimes \mathbb{N}$ le mono\"\i de commutatif libre engendr\'e
par $X$. Le foncteur $f$ est mono\"\i dal sym\'etrique et induit donc
un foncteur sur les sch\'emas affines
$$Aff_{\mathbb{F}_{1}} \longrightarrow Aff_{\mathbb{N}}.$$
Ce foncteur est celui qui envoit un mono\"\i de commutatif $M$ sur
le semi-anneau $\mathbb{N}[M]$, analogue \emph{semi} 
des anneaux en groupes. 

\begin{prop}\label{p7}
Le foncteur
$$f : \mathbb{F}_{1}-Mod \longrightarrow \mathbb{N}-Mod$$
v\'erifie les conditions du corollaire \ref{c1}. 
\end{prop}

\textit{Preuve:} Le fait que $g$ soit conservatif et commute aux colimites
filtrantes est clair. Soit $A \longrightarrow B$ un morphisme plat 
dans $Comm(\mathbb{F}_{1}-Mod)$. Consid\'erons le morphisme
de semi-anneaux induit $\mathbb{N}[A] \longrightarrow \mathbb{N}[B]$. 
On consid\`ere le diagramme de cat\'egories et foncteurs
$$\xymatrix{
\mathbb{N}[A]-Mod \ar[r]^-{g} \ar[d] & A-Mod \ar[d] \\
\mathbb{N}[B]-Mod \ar[r]^-{g}  & B-Mod,}$$
dont on doit montrer qu'il commute \`a isomorphisme pr\`es 
(plus pr\'ecis\'ement que la transformation naturelle entre les deux
compositions possibles est un isomorphisme). On peut remarquer
que le cat\'egorie $\mathbb{N}[A]-Mod$ s'identifie \`a la cat\'egorie
des monoides ab\'eliens dans $A-Mod$. De m\^eme, 
le cat\'egorie $\mathbb{N}[B]-Mod$ s'identifie \`a la cat\'egorie
des monoides ab\'eliens dans $B-Mod$. Le fait que 
le diagramme ci-dessus commute est alors un cons\'equence
du fait que le foncteur de changement de bases
$A-Mod \longrightarrow B-Mod$ commute aux produits directs (qui lui-m\^eme
est un cons\'equence du fait que $A \longrightarrow B$ soit plat). \hfill $\Box$ \\

D'apr\`es le corollaire \ref{c1} on dispose donc d'un foncteur de changement
de bases
$$-\otimes_{\mathbb{F}_{1}}\mathbb{N} : 
\mathbb{F}_{1}-Sch \longrightarrow \mathbb{N}-Sch.$$
En composant avec le changement de bases des 
$\mathbb{N}$-sch\'emas vers les $\mathbb{Z}$-sch\'emas on trouve
un changement de bases
$$-\otimes_{\mathbb{F}_{1}}\mathbb{Z} : 
\mathbb{F}_{1}-Sch \longrightarrow \mathbb{Z}-Sch.$$

\section{Quelques exemples de sch\'emas au-dessous de $Spec\, \mathbb{Z}$}

\subsection{Quelques sch\'emas en groupes}

Soit $M$ un groupe ab\'elien, que l'on peut regarder 
comme un mono\"\i de commutatif. On consid\`ere
$\mathbb{D}_{\mathbb{F}_{1}}(M):=Spec\, M\in \mathbb{F}_{1}-Sch$. 
On remarque imm\'ediatement que 
$$\mathbb{D}_{\mathbb{F}_{1}}(M)\otimes_{\mathbb{F}_{1}}\mathbb{Z} \simeq
Spec\, \mathbb{Z}[M]\simeq \mathbb{D}(M),$$
est le $\mathbb{Z}$-sch\'ema en groupes diagonalisable associ\'e \`a $M$. 
On voit ainsi que tous les sch\'emas en groupes diagonalisables sont 
d\'efinis sur $\mathbb{F}_{1}$. En particulier, on 
pose $Spec\, \mathbb{F}_{1^{n}}:=\mathbb{D}_{\mathbb{F}_{1}}(\mathbb{Z}/n)$.
Ceci nous donne
$$Spec\, \mathbb{F}_{1^{n}}\otimes_{\mathbb{F}_{1}}\mathbb{Z}\simeq
\mu_{n},$$
comme on doit s'y attendre. \\

Examinons maintenant le cas du groupe lin\'eaire
$Gl_{n}$. Revenons un moment au cas g\'en\'eral 
d'une cat\'egorie mono\"\i dale sym\'etrique
$(C,\otimes,\mathbf{1})$ qui v\'erifie les conditions
\ref{h1}. Posons $E_{n}=\{1,\dots,n\}$, l'ensemble \`a n \'el\'ements. On d\'efinit un foncteur
$$\begin{array}{cccc}
Gl_{n,C} : & Comm(C) & \longrightarrow & Ens \\
 & A & \mapsto & Aut(\coprod_{E_{n}}A),
\end{array}$$
o\`u $\coprod_{E_{n}}A$ est le $A$-module libre sur $E_{n}$ (i.e. le 
$A$-module libre de rang $n$), et 
$Aut(\coprod_{E_{n}}A)$ est le groupe des automorphismes
de $\coprod_{E_{n}}A$ dans $A-Mod$. Ce foncteur est un faisceau 
et sera consid\'er\'e comme objet de
$Sh(Aff_{C})$. 
Nous noterons
$$Gl_{n,\mathbb{Z}} \qquad Gl_{n,\mathbb{N}} \qquad Gl_{n,\mathbb{F}_{1}}$$
pour les trois cas de $(C,\otimes,\mathbf{1})$ qui nous int\'eressent.

\begin{prop}\label{p8}
\begin{enumerate}
\item 
Le faisceau $Gl_{n,\mathbb{Z}}$ (resp. $Gl_{n,\mathbb{N}}$, resp. 
$Gl_{n,\mathbb{F}_{1}}$) est repr\'esentable par
un $\mathbb{Z}$-sch\'ema affine (resp. un $\mathbb{N}$-sch\'ema affine,
resp. un $\mathbb{F}_{1}$-sch\'ema). 
\item Le morphisme naturel
$$Gl_{n,\mathbb{N}}\otimes_{\mathbb{N}}\mathbb{Z} \longrightarrow
Gl_{n,\mathbb{Z}}$$
est un isomorphisme. 
\item Le $\mathbb{Z}$-sch\'ema
$Gl_{n,\mathbb{F}_{1}}\otimes_{\mathbb{F}_{1}}\mathbb{Z}$ est isomorphe
au produit semi direct de $\Sigma_{n}$ par $\mathbb{G}_{m}^{n}$. 
\item On a 
$$Gl_{n,\mathbb{Z}}(\mathbb{Z})=Gl_{n}(\mathbb{Z}) \qquad
Gl_{n,\mathbb{N}}(\mathbb{N})=Gl_{n,\mathbb{F}_{1}}(\mathbb{F}_{1})=\Sigma_{n}.$$

\end{enumerate}
\end{prop}

\textit{Preuve:} $(1)$ Commen\c{c}ons par noter $M_{n,C}$ le faisceau des endomorphismes
de $\coprod_{E_{n}}A$
$$\begin{array}{cccc}
M_{n,C} : & Comm(C) & \longrightarrow & Ens \\
 & A & \mapsto & End(\coprod_{E_{n}}A).
\end{array}$$
Supposons pour commencer que dans $C$ les sommes directes sont aussi des produits directs, 
c'est \`a dire que pour deux objets $X$ et $Y$ de $C$, le morphisme naturel
$X\coprod Y \longrightarrow X\times Y$ est un isomorphisme. Ceci est le cas par exemple
pour $C=\mathbb{N}-Mod$ et pour $C=\mathbb{Z}-Mod$. Dans ce cas, on voit facilement que 
$M_{n,C}$ est isomorphe au foncteur
$$\begin{array}{cccc}
M_{n,C} : & Comm(C) & \longrightarrow & Ens \\
 & A & \mapsto & A^{n^{2}}.
\end{array}$$
En d'autres termes, si l'on note $B:=L(\mathbf{1}^{n^{2}})\in Comm(C)$, le mono\"\i de commutatif libre
engendr\'e par l'objet $\mathbf{1}^{n^{2}}$, alors on a $M_{n,C}\simeq Spec\, B$. Ainsi, 
le faisceau $M_{n,C}$ est repr\'esentable par un sch\'ema affine. De plus, il existe des 
carr\'es cart\'esiens de faisceaux
$$\xymatrix{
Gl_{n,C} \ar[r] \ar[d] & M_{n,C} \times M_{n,C} \ar[d]^-{m} \\
\bullet \ar[r]_-{e} & M_{n,C}\times M_{n,C},
}$$
o\`u le morphisme $m$ envoie symboliquement $(u,v)$ sur $(uv,vu)$, et o\`u le morphisme
$e$ est le morphisme unit\'e. Ceci montre donc que 
$Gl_{n,C}$ est repr\'esentable par un sch\'ema affine, toujours sous l'hypoth\`ese que
les sommes directes sont aussi des produits directs dans $C$. Ceci montre donc que 
$Gl_{n,\mathbb{N}}$ et $Gl_{n,\mathbb{Z}}$ sont tous deux des sch\'emas affines.

Il nous reste \`a traiter le cas de $Gl_{n,\mathbb{F}_{1}}$. Mais, pour $A\in Comm(Ens)$, on voit facilement que
le groupe $Aut(\coprod_{E_{n}}A)$ est isomorphe au produit semi-direct 
du groupe sym\'etrique $\Sigma_{n}$ par le groupe $A^{\times}$ des \'el\'ements inversibles dans
le mono\"\i de $A$. En d'autres termes, si l'on note $\mathbb{G}_{m,\mathbb{F}_{1}}:=Gl_{1,\mathbb{F}_{1}}$, on 
trouve que le faisceau en groupes $Gl_{n,\mathbb{F}_{1}}$ est isomorphe au produit semi-direct du faisceau constant
$\Sigma_{n}$ par le faisceau en groupes $\mathbb{G}_{m,\mathbb{F}_{1}}^{n}$. Ainsi, en tant que faisceau d'ensembles
$Gl_{n,\mathbb{F}_{1}}$ est isomorphe \`a une r\'eunion disjointe de $\mathbb{G}_{m,\mathbb{F}_{1}}^{n}$. 
Il nous suffit donc de montrer que $\mathbb{G}_{m,\mathbb{F}_{1}}$ est repr\'esentable par un sch\'ema
affine. Mais on voit facilement que 
$$\mathbb{G}_{m,\mathbb{F}_{1}}\simeq Spec\, \mathbb{Z}\simeq \mathbb{D}(\mathbb{Z}),$$
o\`u
$\mathbb{Z}$ ici est consid\'er\'e comme un mono\"\i de pour la loi additive, et donc comme
un objet de $Comm(Ens)$ (ce $Spec\, \mathbb{Z}$ est un $\mathbb{F}_{1}$-sch\'ema et n'a \'evidemment rien \`a voir avec le 
$Spec\, \mathbb{Z}$ qui est un $\mathbb{Z}$-sch\'ema). 
Ceci termine la preuve du point $(1)$. \\

$(2)$ Ceci est \'evident car pour $A\in Comm(\mathbb{Z}-Mod)$ un anneau commutatif, on a
$$Gl_{n,\mathbb{N}}\otimes_{\mathbb{N}}\mathbb{Z}(A)\simeq 
Gl_{n,\mathbb{N}}(A),$$
o\`u $A$ est aussi condid\'er\'e comme un semi-anneau (ceci est vrai car
$Gl_{n,\mathbb{N}}$ est affine). Comme 
les $A$-modules dans $\mathbb{Z}-Mod$ forment une sous-cat\'egorie 
pleine des $A$-modules dans $\mathbb{N}-Mod$, on voit bien que 
$$Gl_{n,\mathbb{N}}\otimes_{\mathbb{N}}\mathbb{Z}(A)\simeq 
Gl_{n,\mathbb{N}}(A)\simeq Gl_{n}(A)\simeq Gl_{n,\mathbb{Z}}(A).$$

$(3)$ Nous avons vu lors de la preuve de $(1)$ que le faisceau en groupes
$Gl_{n,\mathbb{F}_{1}}$ est le produit semi-direct de
$\Sigma_{n}$ par $\mathbb{G}_{m,\mathbb{F}_{1}}^{n}$. Ainsi, 
$Gl_{n,\mathbb{F}_{1}}\otimes_{\mathbb{F}_{1}}\mathbb{Z}$ est le produit semi-direct 
du faisceau constant $\Sigma_{n}$ par 
$(\mathbb{G}_{m,\mathbb{F}_{1}}\otimes_{\mathbb{F}_{1}}\mathbb{Z})^{n}$.
Mais nous avons vu lors de la discussion sur les groupes de type multiplicatif que 
$$\mathbb{G}_{m,\mathbb{F}_{1}}\otimes_{\mathbb{F}_{1}}\mathbb{Z}\simeq 
\mathbb{D}_{\mathbb{F}_{1}}(\mathbb{Z})\otimes_{\mathbb{F}_{1}}\mathbb{Z}\simeq \mathbb{G}_{m}.$$

$(4)$ Les valeurs de $Gl_{n,\mathbb{Z}}(\mathbb{Z})$ et
$Gl_{n,\mathbb{F}_{1}}(\mathbb{F}_{1})$ ont d\'ej\`a \'et\'e vues. Le groupe
$Gl_{n,\mathbb{N}}(\mathbb{N})$ est par d\'efinition le groupe des automorphismes
de $\mathbb{N}^{n}$ en tant que $\mathbb{N}$-module. Il s'identifie donc au 
sous-groupe  des matrices $n\times n$ \`a coefficients dans $\mathbb{N}$ 
qui poss\`edent un inverse qui soit aussi \`a coefficients dans $\mathbb{N}$. Il s'agit donc
du sous-groupe des matrices de permutation et donc de $\Sigma_{n}$. 
\hfill $\Box$ \\

La conclusion de la proposition \ref{p8} est que $Gl_{n,\mathbb{Z}}$ est d\'efini
sur $\mathbb{N}$, mais pas sur $\mathbb{F}_{1}$ pour $n>1$. 

\subsection{Vari\'et\'es toriques}

Soit $\Delta$ un \'eventail au sens de \cite{od}. Pour un 
\'el\'ement $\sigma \in \Delta$, on dispose d'un 
cone dual $\sigma^{*}$ dont les points entiers forment 
un mono\"\i de commutatif not\'e $M_{\sigma}$. Pour $\sigma' \subset \sigma$ une
face de $\sigma$, on dispose d'un morphisme de mono\"\i des
$M_{\sigma} \longrightarrow M_{\sigma'}$. Si 
l'on pose $U_{\sigma}:=Spec\, M_{\sigma}\in \mathbb{F}_{1}-Sch$, 
alors le morphisme
$$U_{\sigma'} \longrightarrow U_{\sigma}$$
est un ouvert de Zariski. 
En effet, le morphisme de mono\"\i des $M_{\sigma} \longrightarrow M_{\sigma'}$
est obtenu par inversion d'un nombre fini d'\'el\'ements de $M_{\sigma}$ 
(voir Proposition 1.3 \cite{od}). 

On d\'efinit des faisceaux sur $Aff_{\mathbb{F}_{1}}$
$$X:=\coprod_{\sigma\in \Delta}U_{\sigma} \qquad R:=
\coprod_{\sigma,\sigma'\in \Delta}U_{\sigma\cap \sigma'}.$$
Les deux morphismes naturels de $R$ dans $X$ d\'efinissent une relation d'\'equivalence
$$R\subset X\times X.$$
Cette relation d'\'equivalence v\'erifie les conditions de la proposition \ref{p4-} $(2)$, et on obtient ainsi un 
$\mathbb{F}_{1}$-sch\'ema
$$X_{\mathbb{F}_{1}}(\Delta):=X/R \in \mathbb{F}_{1}-Sch.$$

On dispose aussi d'un $\mathbb{Z}$-sch\'ema $X_{\mathbb{Z}}(\Delta)$, 
tel que par exemple pr\'esent\'e dans \cite{od}. 
On voit imm\'ediatement par construction que l'on a un isomorphisme
$$X_{\mathbb{F}_{1}}(\Delta)\otimes_{\mathbb{F}_{1}}\mathbb{Z}\simeq X_{\mathbb{Z}}(\Delta).$$

Ainsi, les vari\'et\'es toriques sont bien d\'efinies sur $\mathbb{F}_{1}$. Elles sont donc aussi d\'efinies
sur $\mathbb{N}$, en posant $X_{\mathbb{N}}(\Delta):=
X_{\mathbb{F}_{1}}(\Delta)\otimes_{\mathbb{F}_{1}}\mathbb{N}$, 
ou encore par une construction directe. 

\section{Les mod\`eles homotopiques}

Dans cette derni\`ere section nous supposerons que
$(C,\otimes,\mathbf{1})$ est une cat\'egorie de mod\`eles 
mono\"\i dale sym\'etrique au sens de \cite{ho}. Remarquer
que par d\'efinition la cat\'egorie mono\"\i dale sym\'etrique sous-jacente
v\'erifie l'hypoth\`ese \ref{h1}. Nous supposerons de plus que
les conditions suivantes sont satisfaites. 

\begin{hyp}\label{h2}
\item La cat\'egorie de mod\`eles $C$ est simpliciale et combinatoire
au sens de \cite[Appendix\, A.2]{hagI}. 
\item La cat\'egorie de mod\`eles mono\"\i dale sym\'etrique
$(C,\otimes,\mathbf{1})$ satisfait \`a l'axiome du mono\"\i de
de \cite{ss}. 
\item L'objet $\mathbf{1}$ est cofibrant dans $C$.
\item Les sources et buts des cofibrations et cofibrations triviales
g\'en\'eratrices sont des objets cofibrants dans $C$. 
\end{hyp}

Sous les hypoth\`ese ci-dessus, nous pouvons consid\'erer la
cat\'egorie $Comm(C)$, des $\mathcal{L}$-alg\`ebres dans $C$, 
o\`u  $\mathcal{L}$ est (l'image dans $C$ de) l'op\'erade des isom\'etries lin\'eaires (voir 
par exemple \cite{sp}). La cat\'egorie $Comm(C)$ est donc un mod\`ele
pour les $E_{\infty}$-mono\"\i des dans $C$, et ses objets
seront donc simplement appel\'es des \emph{$E_{\infty}$-mono\"\i des dans
$C$}. 

D'apr\`es l'hypoth\`ese \ref{h2} et
\cite[Thm. 4.3,Cor. 9.7]{sp} on sait qu'elle est munie d'une structure de cat\'egorie 
de mod\`eles dont les fibrations et les \'equivalences
sont d\'efinies sur les objets sous-jacents dans $C$. De plus, 
pour $A\in Comm(C)$, on peut d\'efinir $A-Mod$, une cat\'egorie
des $A$-modules dans $C$ (voir \cite{sp}). Elle peut aussi 
se voir comme une cat\'egorie de modules (\`a gauche) sur
un mono\"\i de $U(A)$ dans $C$ (voir \cite{sp}).
Ainsi, d'apr\`es nos
hypoth\`eses \ref{h2}, et d'apr\`es 
\cite{ss}, la cat\'egorie $A-Mod$ est aussi munie d'une 
structure de cat\'egorie de mod\`eles o\`u les fibrations et les \'equivalences
sont d\'efinies dans $C$. Enfin, pour un morphisme
$f : A \longrightarrow B$ entre objets dans $Comm(C)$, on dispose d'une adjonction
de Quillen
$$-\otimes_{A}B : A-Mod \longrightarrow B-Mod \qquad
A-Mod \longleftarrow B-Mod : F,$$
o\`u $F$ est le foncteur d'oubli naturel. 
Nous utiliserons la notation usuelle
$$-\otimes^{\mathbb{L}}_{A}B : Ho(A-Mod) \longrightarrow Ho(B-Mod)$$
pour d\'esigner le foncteur d\'eriv\'e \`a gauche
de $-\otimes_{A}B$. Il est d\'emontr\'e
dans \cite[Prop. 9.10]{sp} que lorsque $A$ et $B$ sont cofibrants et que 
$f : A \longrightarrow B$ est une
\'equivalence alors l'adjonction de Quillen pr\'ec\'edente
est une \'equivalence de Quillen. Comme ceci ne semble plus vrai 
lorsque $A$ et $B$ ne sont plus tous deux cofibrants nous
nous interdirons de consid\'erer les cat\'egories
de modules sur des $E_{\infty}$-mono\"\i des qui ne sont pas
cofibrants dans $Comm(C)$. 

Pour finir ces rappels, nous noterons $Map_{C}$, ou simplement 
$Map$ si $C$ est claire, les 
espaces de morphismes (\emph{mapping spaces} en anglais)
de la cat\'egorie de mod\`eles $Comm(C)$, et $Map_{A}$ ceux
de $A-Mod$. Rappelons que comme les cat\'egories de mod\`eles
$Comm(C)$ et $A-Mod$ sont simpliciales, nous pourrons prendre
comme mod\`eles pour ces espaces de morphismes les
ensembles simpliciaux de morphismes entre objets cofibrants et fibrants
(voir \cite{ho}). 

\begin{df}\label{d9}
La \emph{cat\'egorie de mod\`eles des sch\'emas affines sur $C$}
est la cat\'egorie
$$Aff_{C}:=Comm(C)^{op}$$
oppos\'ee de la cat\'egorie de mod\`eles des $E_{\infty}$-mono\"\i des cofibrants. 
\end{df}

On remarquera que lorsque la structure de mod\`ele sur $C$ est
triviale (i.e. les \'equivalences sont les isomorphismes), 
et que l'on prend alors l'enrichissement simplicial trivial (i.e. 
pour $K$ un ensemble simplicial et $X\in C$ on pose
$K\otimes X:=\coprod_{\pi_{0}(K)}X$), alors 
la cat\'egorie $Comm(C)$ est naturellement \'equivalente 
\`a la cat\'egorie des mono\"\i des commutatifs dans $C$. 
Ainsi, $Aff_{C}$ d\'efinie ci-dessus coincide avec 
$Aff_{C}$ pr\'ec\'edemment d\'efini dans \ref{d1} (strictement 
parlant il faut supposer que $C$ soit une cat\'egorie
localement pr\'esentable pour que la structure de mod\`eles
triviale soit combinatoire et satisfasse donc l'hypoth\`ese
\ref{h2}, et cela sera toujours le cas pour les exemples
consid\'er\'es). 

\subsection{Rappel sur la g\'eom\'etrie alg\'ebrique homotopique}

Nous commencerons par d\'efinir des topologies fpqc et de Zariski sur 
$Aff_{C}$, on sens des topologies sur les cat\'egories de mod\`eles
de \cite{hagI}. Il s'agit, bien entendu, d'une g\'en\'eralisation 
de la d\'efinition \ref{d2} que l'on retrouvera en prenant la structure
de mod\`eles triviale sur $C$.

Nous noterons comme auparavant $Spec\, A\in Aff_{C}$  l'objet
correspondant \`a $A\in Comm(C)$. Noter que $Spec\, B \longrightarrow Spec\, A$ est 
une fibration dans $Aff_{C}$ si et seulement si $A \longrightarrow B$
est une cofibration dans $Comm(C)$. 

\begin{df}\label{d10}
Soit $f : Y=Spec\, B \longrightarrow X=Spec\, A$ un morphisme
dans $Aff_{C}$.
\begin{enumerate}
\item Nous dirons que le morphisme $f$
est \emph{plat} si le foncteur correspondant
$$-\otimes^{\mathbb{L}}_{Q(A)}Q(B) : Ho(Q(A)-Mod) \longrightarrow Ho(Q(B)-Mod)$$
commute aux produits fibr\'es homotopiques finis (ici
$Q(A)$ et $Q(B)$ sont des mod\`eles cofibrants dans $Comm(C)$ pour $A$ et $B$).
\item Le morphisme $f$ est un \emph{\'epimorphisme} si 
pour tout $A' \in Comm(C)$, le morphisme
$$f^{*} : Map_{Comm(C)}(B,A') \longrightarrow Map_{Comm(C)}(A,A')$$
est \'equivalent \`a une inclusion de composantes connexes
(i.e induit un morphisme injectif
$$\pi_{0}(Map(B,A')) \hookrightarrow 
\pi_{0}(Map(A,A')),$$
et pour tout $i>0$ et tout $u\in\pi_{0}(Map(B,A'))$ des isomorphismes
$\pi_{i}(Map(B,A'),u) \simeq 
\pi_{i}(Map(A,A'),f^{*}(u)).$). 
\item Le morphisme $f$ est \emph{de pr\'esentation finie}, 
si pour tout diagramme filtrant d'objets $A'_{i} \in Q(A)/Comm(C)$, 
le morphisme naturel
$$Colim_{i} Map_{Q(A)/Comm(C)}(B,A'_{i}) \longrightarrow 
Map_{Q(A)/Comm(C)}(B,Hocolim_{i}A'_{i})$$
est un isomorphisme dans $Ho(SEns)$ (ici $Q(A)$ est un mod\`ele
cofibrant de $A$ dans $Comm(C)$). 
\item 
Le morphisme $f : X=Spec\, B \longrightarrow Y=Spec\, A$ dans $Aff_{C}$
est \emph{un ouvert de Zariski} (ou encore
\emph{une immersion Zariski ouverte}) si le morphisme correspondant
$A \longrightarrow B$ dans $Comm(C)$
est un \'epimorphisme plat et de pr\'esentation finie. 
\end{enumerate}
\end{df}

A l'aide des d\'efinitions pr\'ec\'edentes nous d\'efinissons
les recouvrements fpqc et Zariski de la fa\c{c}on suivante. 

\begin{df}\label{topzarhom}
\begin{enumerate}
\item
Une famille de morphismes 
$$\{X_{i}=Spec\, A_{i} \longrightarrow X=Spec\, A\}_{i\in I}$$
dans $Aff_{C}$ est \emph{un recouvrement fpqc} (ou plus
simplement \emph{recouvrement plat}) si les deux 
conditions suivantes sont satisfaites. 
\begin{enumerate}
\item Pour tout $i\in I$ le morphisme $X_{i} \longrightarrow X$ est 
plat.
\item Il existe un sous ensemble fini $J\subset I$, tel que 
le foncteur
$$\prod_{j\in J} -\otimes_{Q(A)}Q(A_{j}) : Ho(Q(A)-Mod) \longrightarrow 
\prod_{j\in J}Ho(Q(A_{j})-Mod)$$
est conservatif.
\end{enumerate}
\item Une famille de morphismes 
$$\{X_{i} \longrightarrow X\}_{i\in I}$$
dans $Aff_{C}$ est \emph{un recouvrement de Zariski} si
c'est un recouvrement plat et si tous les morphismes
$X_{i} \longrightarrow X$ sont des ouverts de Zariski.
\end{enumerate}
\end{df}

Il est facile de voir que 
les recouvrements fpqc et Zariski d\'efinissent deux
pr\'etopologies sur la cat\'egorie de mod\`eles $Aff_{C}$, 
au sens de \cite[\S 4.3]{hagI} (cela demande 
d'utiliser la formule du changement de base
\cite[Prop. 9.12]{sp}). Les topologies
de Grothendieck associ\'ees \`a ces pr\'etopologies seront 
appel\'ees la \emph{topologie fpqc}, ou plus simplement la \emph{topologie
plate}, et la \emph{topologie de Zariski}. Comme il est d\'emontr\'e dans \cite{hagI} il s'agit 
de  donn\'ees de deux topologies de Grothendieck, au sens usuel, sur
la cat\'egorie $Ho(Aff_{C})$. \\

Comme il est expliqu\'e dans \cite[\S 4.3]{hagI}, la cat\'egorie de mod\`eles
$Aff_{C}$ munie de la topologie fpqc (resp. Zariski) forme un \emph{site de mod\`eles}, 
et on peut donc consid\'erer la cat\'egorie de mod\`eles
$Aff_{C}^{\sim,fpqc}$ (resp. $Aff_{C}^{\sim,Zar}$) des champs sur $Aff_{C}$. Rappelons que
$Aff_{C}^{\sim,fpqc}$ (resp. $Aff_{C}^{\sim,Zar}$) est une localisation de Bousfield \`a gauche 
de la cat\'egorie de mod\`eles des pr\'efaisceaux simpliciaux
sur $Aff_{C}$ dont les objets locaux sont les foncteurs
$$F : Aff_{C}^{op}=Comm(C) \longrightarrow SEns$$
v\'erifiant les deux conditions suivantes. 

\begin{enumerate}
\item Pour toute \'equivalence $X \longrightarrow Y$ dans $Aff_{C}$, 
le morphisme induit 
$$F(Y) \longrightarrow F(X)$$
est un \'equivalence faible d'ensembles simpliciaux. 

\item Pour tout hyper-recouvrement plat (resp. de Zariski) $X_{*} \longrightarrow X$
(voir \cite[\S 4.4]{hagI}), le morphisme naturel
$$F(X)\simeq Map(X,F) \longrightarrow Holim_{[n]\in \Delta}Map(X_{n},F)$$
est un \'equivalence faible d'ensembles simpliciaux. 
\end{enumerate}

\begin{df}\label{d11}
La \emph{cat\'egorie des champs au-dessus de $C$} pour la topologie plate (resp. 
de Zariski)
est la cat\'egorie homotopique 
$Ho(Aff_{C}^{\sim,fpqc})$ (resp. $Ho(Aff_{C}^{\sim,Zar})$). Comme pr\'ec\'edemment
nous adopterons les notations
$$Ch^{fpqc}(Aff_{C}):=Ho(Aff_{C}^{\sim,fpqc}) \qquad 
Ch(Aff_{C}):=Ho(Aff_{C}^{\sim,Zar}).$$
Les objets de $Ch(Aff_{C})$ seront simplement appel\'es des
\emph{champs au-dessus de $C$} (ou simplement des \emph{champs}
si le contexte est clair).  
\end{df}

Comme pour le cas non homotopique pr\'esent\'e pr\'ec\'edemment nous nous int\'eresserons
principalement aux champs pour la topologie de Zariski et les champs fpqc ne seront
utilis\'es que de mani\`ere auxiliaire. \\

Il faut noter que lorsque la structure de mod\`eles est triviale sur $C$ alors
$Aff^{\sim,fpqc}_{C}$ est la cat\'egorie de mod\`eles des pr\'efaisceaux
simpliciaux au sens de \cite{ja} (ou plus pr\'ecis\'ement sa version
projective d\'ecrite dans \cite{bl}). Ainsi, 
la cat\'egorie $Ch(Aff_{C})$ est la cat\'egorie homotopique des
pr\'efaisceaux simpliciaux sur le site $Aff_{C}$, qui s'identifie aussi
\`a la cat\'egorie homotopique des objets simpliciaux dans
le topos $Sh(Aff_{C})$ des faisceaux sur $Aff_{C}$. La cat\'egorie
$Sh(Aff_{C})$ s'identifie alors \`a la sous-cat\'egorie pleine 
de $Ch(Aff_{C})$ des pr\'efaisceaux simpliciaux discrets (i.e. 
dont les valeurs sont des ensembles vus comme ensembles
simpliciaux constants). Ceci nous permettra par la suite
de voir les faisceaux sur $Aff_{C}$ comme des champs
sur $Aff_{C}$, et d'identifier la cat\'egorie $Sh(Aff_{C})$
\`a son image dans $Ch(Aff_{C})$. \\

Venons-en maintenant \`a la d\'efinition de sch\'emas 
relatifs \`a $C$. Pour cela nous aurons besoin
de la g\'en\'eralisation suivante de la d\'efinition \ref{d3}. 
Rappelons que l'on peut d\'efinir une notion 
de \emph{h-monomorphisme} dans $Ch(Aff_{C})$, qui 
sont les morphismes de champs $F\longrightarrow G$ tel que le morphisme
diagonal $F\longrightarrow F\times_{G}^{h}F$ soit un isomorphisme 
dans $Ch(Aff_{C})$. Il s'agit d'une g\'en\'eralisation 
de nature homotopique de la notion usuelle de monomorphisme. Les 
h-monomorphismes dans $Ch(Aff_{C})$ sont aussi les morphismes
qui induisent un monomorphisme sur les faisceaux $\pi_{0}$ et des 
isomorphismes sur tous les faisceaux $\pi_{i}$ pour $i>0$ (voir 
\cite{hagI} pour la d\'efinition des faisceaux d'homotopie). Cette notion 
de h-monomorphisme permet aussi de donner un sens \`a  la notion 
de sous-objets d'un objet $F$ de $Ch(Aff_{C})$, appel\'es
\emph{sous-champs}, comme classe d'isomorphisme de h-monomorphismes
de but $F$. 

\begin{df}\label{d12}
\begin{enumerate}
\item 
Soit $X$ un sch\'ema affine et $F\subset X$ un sous-champ
de $X$. Nous dirons que $F$ est un \emph{ouvert de Zariski
de $X$} s'il existe une famille d'ouverts de Zariski 
$\{X_{i} \longrightarrow X\}_{i\in I}$ dans $Aff_{C}$ (au sens
de la d\'efinition \ref{d10}) tel que 
$F$ soit l'image du morphisme de champs
$$\coprod_{i\in I} X_{i} \longrightarrow X.$$
\item Un morphisme $f : F \longrightarrow G$ dans $Ch(Aff_{C})$ est 
un \emph{ouvert de Zariski} (ou encore une \emph{immersion 
Zariski ouverte}) si pour tout sch\'ema affine $X$ et tout
morphisme $X \longrightarrow G$, le morphisme induit
$$F\times_{G}X \longrightarrow X$$
est un h-monomorphisme d'image un ouvert de Zariski de $X$. 
\end{enumerate}
\end{df}

On v\'erifie sans peine que les ouverts de Zariski 
se comportent comme on s'y attend (stabilit\'e par 
\'equivalences, compositions et produits fibr\'es homotopiques
dans $Aff_{C}$.)   

Comme il est expliqu\'e dans \cite[\S 4.2]{hagI} on dispose d'un plongement 
de Yoneda
$$\mathbb{R}\underline{h}_{-} : Ho(Aff_{C}) \longrightarrow 
Ho(Aff_{C}^{\wedge}),$$
qui permet d'identifier la cat\'egorie $Ho(Aff_{C})$ \`a une 
sous-cat\'egorie pleine de la cat\'egorie 
homotopique des pr\'e-champs $Ho(Aff_{C}^{\wedge})$. 
Malheureusement, nous ne savons pas montrer que 
la topologie plate est sous-canonique en g\'en\'eral, 
et ainsi nous ne savons pas si l'analogue du corollaire
\ref{p1} reste valable (en r\'ealit\'e nous n'avons  pas trouv\'e de
version homotopique raisonable du th\'eor\`eme \ref{mdesc}).
En particulier nous ne pouvons pas
d\'eduire que le plongement de Yoneda ci-dessus se factorise
par la sous-cat\'egorie pleine $Ch(Aff_{C}) \subset Ho(Aff_{C}^{\wedge})$. 
Nous sommes donc oblig\'es de rajouter l'hypoth\`ese suivante.

\begin{hyp}\label{h3}
La topologie plate sur $Aff_{C}$ est sous-canonique. En d'autres termes, 
pour tout $X\in Aff_{C}$ le pr\'e-champ
$\mathbb{R}\underline{h}_{X}$ est un champ. 
\end{hyp}

L'hypoth\`ese pr\'ec\'edente implique que le plongement de Yoneda induit 
un foncteur pleinement fid\`ele
$$\mathbb{R}\underline{h}_{-} : Ho(Aff_{C}) \longrightarrow Ch(C).$$
Ceci permet d'identifier la cat\'egorie homotopique
$Ho(Aff_{C})$ comme une sous-cat\'egorie pleine de $Ch(Aff_{C})$. 
Les objets dans l'image de $\mathbb{R}\underline{h}_{-}$ 
seront simplement appel\'es des \emph{sch\'emas affines au-dessus de $C$}
(ou simplement \emph{sch\'emas affines} si $C$ est claire).

Enfin, comme $Aff_{C}^{\sim,fpqc}$ est une cat\'egorie de mod\`eles on peut
parler de limites et colimites homotopiques, et en particulier de
produits fibr\'es homotopiques. Ces derniers seront not\'es
$-\times_{-}^{h}-$. Remarquer que le produit direct dans
$Aff_{C}^{\sim,Zar}$ ne pr\'eserve pas les \'equivalences en g\'en\'eral, et donc
que l'on a aussi une notion de produit homotopique
not\'e $-\times^{h}-$. Il s'agit bien entendu du produit direct dans 
la cat\'egorie $Ch(Aff_{C})$. 

On d\'efinit alors la notion de sch\'ema
relatif dans $Aff_{C}$ de la fa\c{c}on suivante.

\begin{df}\label{d13}
Un champ $F$ est un \emph{sch\'ema relatif \`a $C$} (ou simplement 
\emph{sch\'ema} si le contexte est clair) 
existe une famille de sch\'emas affines $\{X_{i}\}_{i}$ et un \'epimorphisme
$$\coprod_{i\in I} X_{i} \longrightarrow F$$
tel que chaque morphisme $X_{i} \longrightarrow F$ soit  un ouvert
Zariski. 
\end{df}

\begin{df}\label{d14}
La \emph{cat\'egorie des sch\'emas relatifs \`a $C$}
est la sous-cat\'egorie pleine de $Ch(Aff_{C})$ form\'ee
des sch\'emas au sens de la d\'efinition \ref{d13}. Nous la
noterons $Sch(C)$. 
\end{df}

Tout comme pour le cas non homotopique on montre les 
propri\'et\'es de stabilit\'e suivantes. 

\begin{prop}\label{p9}
La sous-cat\'egorie $Sch(C)$ de
$Ch(Aff_{C})$ est stable par r\'eunions disjointes et
par produits fibr\'es homotopiques. 
\end{prop}

\textit{Preuve:} C'est la m\^eme que pour 
la proposition \ref{p4-}. \hfill $\Box$ \\

On montre de m\^eme la proposition suivante, analogue homotopique 
de la proposition \ref{p4-} $(2)$. Pour cela on renvoit \`a \cite{hagI,hagII} pour la notion 
de groupoides de Segal dans une cat\'egorie de mod\`eles, qui est un analogue homotopique
de la notion usuelle de relation d'\'equivalence. 

\begin{prop}\label{p9+}
Un champ $F$ est un sch\'ema si et seulement s'il existe un 
groupoide de Segal $X_{*}$ dans $Aff_{C}^{\sim,Zar}$ 
v\'erifiant les conditions suivantes. 
\begin{enumerate}
\item On a 
$$X_{0}\simeq \coprod_{i\in I} U_{i}$$
avec $U_{i}$ des sch\'emas affines. 
\item Pour tout $(i,j)\in I^{2}$, consid\'erons le sous-champ
$R_{i,j}\subset U_{i}\times^{h} U_{j}$ d\'efini par la carr\'e homotopiquement 
cart\'esien suivant
$$\xymatrix{
R_{i,j} \ar[r] \ar[d] & U_{i}\times^{h} U_{j} \ar[d] \\
X_{1} \ar[r] & X_{0}\times^{h}X_{0}.
}$$
Alors chacun des morphismes 
$$R_{i,j} \longrightarrow U_{i}$$
est un ouvert Zariski.
\item Pour tout $i\in I$ le morphisme $R_{i,i} \longrightarrow U_{i}$ est 
un isomorphisme. 
\item On a 
$$F\simeq |X_{*}|=Hocolim_{n\in \Delta^{op}}X_{n}.$$
\end{enumerate}
\end{prop}

Tout comme dans le cas que nous avons trait\'e pr\'ec\'edemment, 
on peut pour tout sch\'ema $X$ sur $C$ au sens ci-dessus, d\'efinir
un site de Grothendieck $Zar(X)$ des ouverts de Zariski de $X$. 
Ce site poss\`ede des produits fibr\'es, il est engendr\'e par une 
pr\'etopologie quasi-compacte, 
et sa cat\'egorie sous-jacente est un lieu. On en d\'eduit donc que 
le topos des champs $Ch(Zar(X))$ est \'equivalent  aux champs  sur
un espace topologique bien d\'efini $|X|$. L'espace topologique
$|X|$ est appel\'e comme il se doit \emph{l'espace sous-jacent \`a $X$}. 
Par construction, le site $Ouv(X)$ des ouverts de $X$ 
est \'equivalent au site $Zar(X)$. 

Sans entrer dans les d\'etails signalons aussi que 
l'on peut construire un pr\'efaisceau $\mathcal{O}_ {X}$ sur $|X|$ \`a 
valeurs dans $Comm(C)$. Notre hypoth\`ese \ref{h3} implique que
ce pr\'efaisceau est un champ, au sens o\`u il satisfait \`a une 
condition de descente pour les hyper-recouvrements. Ainsi, 
le couple $(|X|,\mathcal{O}_{X})$ joue le r\^ole de l'espace
annel\'e sous-jacent au sch\'ema $X$, en un sens relatif
\`a $C$. Nous n'irons pas plus loin dans cette direction. \\

Pour terminer, signalons l'existence de foncteurs de changements
de bases g\'en\'eralisant ceux du paragraphe \S 2.3. On se donne
une adjonction de Quillen
$$f : (C,\otimes,\mathbf{1}) \longrightarrow (D,\otimes,\mathbf{1}) \qquad
(C,\otimes,\mathbf{1}) \longleftarrow (D,\otimes,\mathbf{1}) : g$$
avec $f$ de Quillen \`a gauche et sym\'etrique mono\"\i dal. On en d\'eduit
une adjonction de Quillen sur les cat\'egorie de $E_{\infty}$-mono\"\i des
$$f : Comm(C) \longrightarrow Comm(D) \qquad
Comm(C) \longleftarrow Comm(D) : g.$$
Comme il est expliqu\'e dans \cite[\S 4.8]{hagI} cette adjonction de Quillen
induit une adjonction sur les cat\'egorie de pr\'e-champs
$$\mathbb{L}f_{!} : Ho(Aff_{C}^{\wedge}) \longrightarrow 
Ho(Aff_{D}^{\wedge}) \qquad 
Ho(Aff_{C}^{\wedge}) \longleftarrow 
Ho(Aff_{D}^{\wedge}) : \mathbb{R}f^{*}.$$
Lorsque le foncteur $f$ est continu (i.e. $\mathbb{R}f^{*}$ pr\'eserve les
sous-cat\'egorie de champs) alors on obtient une adjonction sur les
cat\'egories de champs
$$\mathbb{L}f_{!}^{\sim} : Ch(C) \longrightarrow 
Ch(D) \qquad 
Ch(C) \longleftarrow 
Ch(D) : \mathbb{R}f^{*}.$$
Enfin, pour $A\in Comm(C)$ cofibrant, on dispose d'une adjonction
de Quillen
$$f : A-Mod \longrightarrow f(A)-Mod \qquad
A-Mod \longleftarrow f(A)-Mod : g.$$
Cette derni\`ere adjonction est de plus compatible aux changements 
de bases (voir \cite{sp} pour les d\'etails). \\

On dispose alors de l'\'enonc\'e analogue au corollaire \ref{c1}. 

\begin{prop}\label{p10}
On suppose que le foncteur d\'eriv\'e 
$$\mathbb{R}g : Ho(D) \longrightarrow Ho(C)$$ 
est conservatif et qu'il
commute aux colimites homotopiques filtrantes. On suppose de plus que
pour tout morphisme plat $A \longrightarrow B$ dans
$Comm(C)$, et tout $N\in f(A)-Mod$, le morphisme naturel
$$\mathbb{R}g(N)\otimes^{\mathbb{L}}_{A}B \longrightarrow 
\mathbb{R}g(N\otimes^{\mathbb{L}}_{f(A)}f(B))$$
est un isomorphisme dans $Ho(B-Mod)$. 
Alors, le foncteur 
$f : Aff_{C} \longrightarrow Aff_{D}$ est continu pour la topologie
de Zariski, et le foncteur induit sur les cat\'egories de champs
$$\mathbb{L}f^{\sim}_{!} : Ch(Aff_{C}) \longrightarrow Ch(Aff_{D})$$
pr\'eserve les sous-cat\'egories des sch\'emas et induit un foncteur
$$\begin{array}{ccc}
Sch(C) & \longrightarrow & Sch(D) \\
 X & \mapsto & X\times_{C}D:=f^{\sim}_{!}(X).
\end{array}$$
\end{prop}

\textit{Preuve:} C'est essentiellement la m\^eme que pour le corollaire
\ref{c1}. \hfill $\Box$ \\

\subsection{$Spec\, \mathbb{S}$ : La nouvelle g\'eom\'etrie alg\'ebrique courageuse}

Dans ce paragraphe nous posons $(C,\otimes,\mathbf{1})=(\mathcal{GS},\wedge,\mathbb{S})$, 
la cat\'egorie de mod\`eles mono\"\i dales des $\Gamma$-espaces 
telle que pr\'esent\'ee dans \cite{sch}.  
Rappelons que le cat\'egorie sous-jacente \`a 
$\mathcal{GS}$ est la cat\'egorie des foncteurs
$$F : \Gamma \longrightarrow SEns,$$
o\`u $\Gamma$ est la cat\'egorie des ensembles finis point\'es, 
v\'erifiant $F(*)=*$. Nous noterons $n^{+}$ l'ensemble
$\{0,\dots,n\}$ point\'es en $0$. 
La structure de mod\`eles sur 
$\mathcal{GS}$ est une localisation de la structure 
de mod\`eles niveaux par niveaux (o\`u fibrations et 
\'equivalences sont d\'efinies termes \`a termes), et dont les objets locaux
sont les $\Gamma$-espaces \emph{tr\`es sp\'eciaux}, c'est \`a dire
qui v\'erifient les deux conditions suivantes
\begin{enumerate}
\item Pour $I$ et $J$ deux ensembles finis point\'es, le morphisme naturel
$$F(I\vee J) \longrightarrow F(I)\times F(J)$$
est une \'equivalence faible d'ensembles simpliciaux (ici
$I\vee J$ d\'esigne le coproduit de $I$ et $J$ dans $\Gamma$). 

\item D'apr\`es la condition ci-dessus, le digramme dans $Ho(SEns)$
$$\xymatrix{F(1^{+})\times F(1^{+}) & \ar[r] F(2^{+}) \ar[r] & F(1^{+})}$$
d\'efinit une structure de mono\"\i de sur $F(1^{+})$ dans 
la cat\'egorie $Ho(SEns)$. On demande alors que 
le mono\"\i de $\pi_{0}(F(1^{+}))$ soit une loi de groupes. 
\end{enumerate}

La cat\'egorie $\Gamma$ poss\`ede une structure mono\"\i dale sym\'etrique
induite par le produits directs d'ensembles point\'es. Par convolution ceci 
d\'efinit une structure mono\"\i dale sym\'etrique sur la cat\'egorie 
$\mathcal{GS}$, not\'ee $-\wedge-$. L'objet neutre pour cette
structure mono\"\i dale est le $\Gamma$-espace en sph\`ere
$\mathbb{S}$ corepr\'esent\'e par $1^{+}\in \Gamma$ et donn\'e par
$$\mathbb{S}(n^{+})=Hom_{\Gamma}(1^{+},n^{+})=\{0,\dots,n\}.$$
On montre que $(\mathcal{GS},\wedge,\mathbb{S})$ est une cat\'egorie
de mod\`eles mono\"\i dale sym\'etrique qui v\'erifie
les conditions \ref{h2} (voir \cite{sch,mmss}). 

De fa\c{c}on g\'en\'erale il faut penser aux 
$\Gamma$-espaces tr\'es sp\'eciaux comme \`a des 
groupes simpliciaux commutatifs \`a homotopie de coh\'erences
pr\`es. La structure mono\"\i dale $\wedge$ est alors l'analogue
du produit tensoriel de groupes ab\'eliens. Ceci se justifie en montrant que la cat\'egorie
mono\"\i dal sym\'etrique
$(Ho(\mathcal{GS}),\wedge^{\mathbb{L}})$ est naturellement \'equivalente
\`a la cat\'egorie homotopique des spectres sym\'etriques connectifs
(i.e. sans
homotopie n\'egative) munie du 
smash produit. On montre m\^eme que les th\'eories homotopiques
de mono\"\i des, $E_{\infty}$-mono\"\i des et de modules associ\'ees sont \'equivalentes 
(voir \cite{mmss} pour plus de d\'etails). En conclusion, la th\'eorie
homotopique des $\Gamma$-espaces tr\`es sp\'eciaux est 
\'equivalente \`a celle des spectres connectifs, et ce de fa\c{c}on compatible
avec les constructions tensorielles. Ainsi, la cat\'egorie
de mod\`eles $Comm(\mathcal{GS})$ est un mod\`ele pour la th\'eorie
homotopique des \emph{nouveaux anneaux commutatifs courageux} (traduction, 
simpliste mais amusante,
de \emph{brave new commutative rings}), d'o\`u notre terminologie
de \emph{nouvelle g\'eom\'etrie alg\'ebrique courageuse}. \\

Pour un objet $X \in \mathcal{GS}$ on d\'efinit ses groupes d'homotopie stables
par la formule
$$\pi_{i}(X):=\pi_{i}(RX(1^{+})),$$
o\`u $RX$ est un mod\`ele fibrant de $X$ dans
$\mathcal{GS}$, et $RX(1^{+})$ est son ensemble simplicial 
sous-jacent. Pour tout $i\geq 0$, $\pi_{i}(X)$ est toujours
un groupe ab\'elien. 

Pour un objet $A\in Comm(\mathcal{GS})$, les groupes
d'homotopie $\pi_{i}(A)$ forment un anneau commutatif gradu\'e
que l'on notera $\pi_{*}(A)$. De m\^eme, pour $A$ cofibrant 
dans $A\in Comm(\mathcal{GS})$ et 
$M\in A-Mod$, $\pi_{*}(M)$ est naturellement un 
$\pi_{*}(A)$-module gradu\'e. 

\begin{lem}\label{l1}
La cat\'egorie de mod\`eles mono\"\i dale sym\'etrique
$(\mathcal{GS},\wedge,\mathbb{S})$ v\'erifie l'hypoth\`ese \ref{h3}. 
\end{lem}

\textit{Preuve:} C'est essentiellement la m\^eme que pour le cas
des anneaux simpliciaux commutatifs et des anneaux en spectres qui sont 
trait\'es dans \cite[\S 2.2, \S 2.4]{hagII} (on pourrait aussi 
utiliser essentiellement la m\^eme preuve que
celle de la proposition \ref{pconj1}). En deux mots, on commence
par v\'erifier qu'un morphisme $A \longrightarrow B$ 
dans $Comm(\mathcal{CS})$ est 
est plat si et seulement s'il v\'erifie les deux conditions suivantes
\begin{enumerate}
\item Le morphisme induit $\pi_{0}(A) \longrightarrow \pi_{0}(B)$
est un morphisme plat d'anneaux commutatifs. 
\item Pour tout $i>0$, le morphisme naturel
$$\pi_{i}(A)\otimes_{\pi_{0}(A)}\pi_{0}(B) \longrightarrow 
\pi_{i}(B)$$
est un isomorphisme.
\end{enumerate}
Ensuite, si $A \longrightarrow B_{*}$ est un hyper-recouvrement plat, 
le morphisme induit $\pi_{*}(A) \longrightarrow \pi_{*}(B)_{*}$ est un 
hyper-recouvrement plat d'anneaux commutatifs gradu\'es. Un argument 
de suite spectrale et la descente plate usuelle des anneaux implique alors
que le morphisme
$$A \longrightarrow Holim_{[n]\in \Delta^{op}}B_{n}$$
est une \'equivalence dans $Comm(\mathcal{GS})$. 
\hfill $\Box$ \\

Le lemme \ref{l1} et les consid\'erations g\'en\'erales
du paragraphe pr\'ec\'edent nous donne une notion
de sch\'emas relatifs au-dessus de $\mathcal{GS}$, dont 
la cat\'egorie sera not\'ee $\mathbb{S}-Sch$. 

\begin{df}\label{d15}
La cat\'egorie des sch\'emas relatifs
\`a $\mathcal{GS}$ sera not\'ee
$\mathbb{S}-Sch$. Ses objets seront appel\'es
des \emph{$\mathbb{S}$-sch\'emas}, ou si l'on veut
rire des \emph{nouveaux sch\'emas courageux}.
\end{df}

Noter 
que $Ho(Comm(\mathcal{GS}))$, qui est \'equivalente \`a la cat\'egorie
des $\mathbb{S}$-sch\'emas affines,  est naturellement \'equivalente
\`a la cat\'egorie des $\mathbb{S}$-alg\`ebres commutatives
connectives de \cite{ekmm}. Ceci justifie notre terminologie de 
$\mathbb{S}$-sch\'emas, ou encore de nouveaux sch\'emas courageux. \\

Il existe une adjonction de Quillen 
$$\pi_{0} : \mathcal{GS} \longrightarrow \mathbb{Z}-Mod \qquad
\mathcal{GS} \longleftarrow \mathbb{Z}-Mod : j,$$
qui est telle que $\pi_{0}$ soit de plus un foncteur mono\"\i dal sym\'etrique.
On voit facilement (par exemple \`a l'aide de la description 
des morphismes plats donn\'ee lors de la preuve du lemme
\ref{l1}) que le foncteur induit
$$\pi_{0} : Comm(C) \longrightarrow Comm(\mathbb{Z}-Mod)$$
est continu pour la topologie Zariski et que $j$ est conservatif (il est 
m\^eme pleinement fid\`ele) et commute avec 
les colimites filtrantes. La proposition \ref{p10} nous donne alors un 
foncteur de changement de bases
$$-\otimes_{\mathbb{S}}\mathbb{Z} : \mathbb{S}-Sch \longrightarrow \mathbb{Z}-Sch$$
des $\mathbb{S}$-sch\'emas vers les $\mathbb{Z}$-sch\'emas. Ce foncteur
poss\`ede en r\'ealit\'e un adjoint \`a droite
$$i :  \mathbb{Z}-Sch \longrightarrow  \mathbb{S}-Sch$$
qui est pleinement fid\`ele, et la situation est donc tr\`es proche
de celle du changement de bases de $Spec\, \mathbb{Z} \longrightarrow Spec\, \mathbb{N}$. 
Le morphisme naturel $Spec\, \mathbb{Z} \longrightarrow Spec\, \mathbb{S}$ doit 
aussi \^etre pens\'e comme un monomorphisme, et m\^eme comme
une immersion ferm\'ee. 

Nous d\'efinissons un $\mathbb{S}$-sch\'ema $Gl_{n,\mathbb{S}}$ 
de la fa\c{c}on suivante. Pour $A \in Comm(\mathcal{GS})$
on consid\`ere $QA$ un remplacement cofibrant de $A$, ainsi 
que $QA^{n} \in QA-Mod$ le $QA$-module libre de rang
$n$ sur $QA$. L'objet $QA^{n} \in QA-Mod$ est carat\'eris\'e
par le fait qu'un morphisme de $QA$-modules
$QA^{n} \longrightarrow M$ est la m\^eme chose qu'un morphisme
d'ensemble simpliciaux $E_{n} \longrightarrow M(1^{+})$, o\`u 
$E_{n}$ est l'ensemble fini \`a $n$ \'el\'ements vu comme
un ensemble simplicial discret, et 
$M(1^{+})$ est l'ensemble simplicial sous-jacent \`a l'objet 
$M\in \mathcal{GS}$. On d\'efinit alors un pr\'efaisceau simplicial
$$\begin{array}{cccc}
Gl_{n,\mathbb{S}} : & Comm(\mathcal{GS}) & \longrightarrow & SEns \\
 & A & \mapsto & Aut_{QA-Mod}(QA^{n}),
\end{array}$$
o\`u $Aut_{QA-Mod}(QA^{n})$ est le sous-ensemble simplicial
de $Map_{QA-Mod}(QA^{n},QA^{n})$ consistant en les
morphismes qui sont des \'equivalences (nous laissons le soin 
au lecteur de rendre cette construction 
rigoureusement fonctorielle, voir aussi \cite[\S 1.3.7]{hagII} pour plus de d\'etails).
On voit que ce pr\'efaisceau est un champ et on le consid\`ere
alors comme un objet dans $Ch(Aff_{\mathcal{GS}})$. 

\begin{prop}\label{p16}
\begin{enumerate}
\item 
Le champ $Gl_{n,\mathbb{S}}$ est repr\'esentable par un $\mathbb{S}$-sch\'ema
affine. 
\item On a un isomorphisme
$$Gl_{n,\mathbb{S}}\otimes_{\mathbb{S}}\mathbb{Z} \longrightarrow
Gl_{n,\mathbb{Z}}.$$
\end{enumerate}
\end{prop}

\textit{Preuve:} $(1)$ C'est essentiellement la m\^eme que pour 
la proposition \ref{p8} $(1)$. \\

$(2)$ Il suffit de constater que pour tout $A\in Comm(\mathbb{Z}-Mod)$,  
le foncteur de Quillen \`a droite
$$i : A-Mod \longrightarrow i(A)-Mod$$
induit un foncteur pleinement fid\`ele
$$\mathbb{R}i : Ho(A-Mod) \longrightarrow Ho(i(A)-Mod).$$
\hfill $\Box$ \\

On remarque qu'il n'est pas vrai que 
$i(Gl_{n,\mathbb{Z}})$ soit isomorphe \`a $Gl_{n,\mathbb{S}}$, o\`u 
$i : \mathbb{Z}-Sch \longrightarrow \mathbb{S}-Sch$ est le foncteur
d'inclusion. En effet, $Gl_{n,\mathbb{S}}(\mathbb{S})$ est un $\mathbb{S}$-sch\'ema affine
de la forme $Spec\, A$ pour un certain $A\in Comm(\mathcal{GS})$, qui v\'erifie que
pour tout $B\in Comm(\mathcal{GS})$ il existe un isomorphisme naturel
$$\pi_{0}(Map_{Comm(\mathcal{GS})}(A,B))\simeq Gl_{n}(\pi_{0}(B)).$$
On peut voir que l'anneau 
$\mathbb{Z}[T_{ij},Det(T_{ij})^{-1}]$, vu comme objet de $Comm(\mathcal{GS})$
ne v\'erifie pas cette propri\'et\'e universelle, car par exemple
il n'existe pas de morphisme $\mathbb{Z}[T_{ij},Det(T_{ij})^{-1}] \longrightarrow \mathbb{S}$
dans $Ho(Comm(\mathcal{GS}))$ (si l'on compose avec 
$\mathbb{Z} \longrightarrow \mathbb{Z}[T_{ij},Det(T_{ij})^{-1}]$ on trouverait une
section de la projection $\mathbb{S} \longrightarrow \mathbb{Z}=\pi_{0}(\mathbb{S})$
que l'on sait ne pas exister).  \\

\subsection{L'anneaux en spectres \`a un \'el\'ement $\mathbb{S}_{1}$}

Dans ce paragraphe nous posons $(C,\otimes,\mathbf{1})=(SEns,\times,*)$ 
la cat\'egorie de mod\`eles mono\"\i dale sym\'etrique des
ensembles simpliciaux munis du produit direct. Elle v\'erifie \'evidemment 
les hypoth\`eses de \ref{h2}. 

\begin{prop}\label{pconj1}
La cat\'egorie de mod\`eles mono\"\i dale sym\'etrique
$(SEns,\times,*)$ v\'erifie l'hypoth\`ese \ref{h3}. 
\end{prop}

\textit{Preuve:} Soit 
$F_{*} \longrightarrow X=Spec\, A$ un hyper-recouvrement 
plat et pseudo-repr\'esentable dans $Aff_{C}^{\sim,fpqc}$, au sens 
de \cite[\S 4.4]{hagI}. Par d\'efinition on peut \'ecrire
$F_{n}=\coprod_{\alpha} Spec\, B_{n,\alpha}$, avec 
tous les morphismes $A \longrightarrow B_{n,\alpha}$ plats. 
Posons $B_{n}:=\prod_{\alpha}B_{n,\alpha}$ pour tout
$n$. 
Il nous faut montrer que le morphisme naturel
$$A\longrightarrow Holim_{[n]\in \Delta}B_{n}$$
est un isomorphisme dans $Ho(Comm(SEns))$.

Pour d\'emontrer cela nous allons utiliser les foncteurs
de troncations $B \mapsto B_{\leq k}$ de
$Comm(SEns)$ dans elle-m\^eme. Rappelons qu'un 
objet $B\in Comm(SEns)$ est $k$-tronqu\'e si 
pour tout point $x\in \pi_{0}(B)$ on a $\pi_{i}(B,x)=0$ pour
$i>k$ (i.e. l'ensemble simplicial sous-jacent \`a $B$ est $k$-tronqu\'e).
On montre facilement (par exemple \`a l'aide des techniques
de localisation de cat\'egories de mod\`eles) que le foncteur 
d'inclusion des objets $k$-tronqu\'es
$$Ho(Comm(SEns)_{\leq k}) \longrightarrow Ho(Comm(SEns))$$
poss\'ede un adjoint \`a gauche 
$$(-)_{\leq k} : Ho(Comm(SEns)) \longrightarrow Ho(Comm(SEns)_{\leq k}).$$
On peut m\^eme r\'ealiser cette adjonction comme 
l'adjonction d\'eriv\'ee d'une adjonction de Quillen
$$Comm(SEns) \longrightarrow Comm(SEns)_{\leq k} \qquad
Comm(SEns) \longleftarrow Comm(SEns)_{\leq k},$$
o\`u la cat\'egorie sous-jacente $Comm(SEns)_{\leq k}$ est la m\^eme que
celle de $Comm(SEns)$, mais ses \'equivalences
sont les $k$-\'equivalences (i.e. les morphismes induisant 
des isomorphismes sur les groupes d'homotopie $\pi_{i}$
pour $i\leq k$). 

On commence alors par remarquer le lemme suivant.

\begin{lem}\label{lconj1}
Soit $A \longrightarrow B$ un morphisme plat 
dans $Comm(SEns)$, alors 
pour tout entier $k\geq 0$ le morphisme induit 
$$A_{\leq k} \longrightarrow B_{\leq k}$$
est un morphisme plat dans $Comm(SEns)$. 
\end{lem}

\textit{Preuve du lemme:} Tout d'abord, par platitude on voit facilement que
le foncteur de changement de bases
$$-\times^{\mathbb{L}}_{A}B : Ho(A-Mod) \longrightarrow Ho(B-Mod)$$
pr\'eserve les objets $k$-tronqu\'es. Ainsi, si l'on forme le carr\'e
homotopiquement cocart\'esien dans $Ho(Comm(SEns))$ suivant
$$\xymatrix{
A \ar[r] \ar[d] & B \ar[d] \\
A_{\leq k} \ar[r] & B',}$$
on voit que $B'$ est $k$-tronqu\'e. Ceci implique facilement 
que le morphisme naturel $B_{\leq k} \longrightarrow B'$ est un
\'equivalence, et donc que le diagramme suivant
$$\xymatrix{
A \ar[r] \ar[d] & B \ar[d] \\
A_{\leq k} \ar[r] & B_{\leq k}}$$
est homotopiquement cocart\'esien. En particulier on voit que
$A_{\leq k} \longrightarrow B_{\leq k}$ est un morphisme plat. 
\hfill $\Box$ \\

Pour tout objet $B\in Comm(SEns)$, on dispose d'un 
isomorphisme naturel dans $Ho(Comm(SEns))$
$B \simeq Holim_{k}B_{\leq k}$. Ainsi, le lemme \ref{lconj1}
permet de tronquer termes \`a termes l'hyper-recouvrement 
$F_{*} \longrightarrow X$, et donc de se ramener au cas
o\`u il existe un entier $k$ tel que $A$ soit $k$-tronqu\'e, ainsi que tous les
$B_{n,\alpha}$, et donc tous les $B_{n}$. 

Dans ce cas, la limite homotopique 
$Holim_{[n]\in \Delta}B_{n}$ se r\'eduit en une limite homotopique finie, 
ce qui implique que l'on peut supposer que l'hyper-recouvrement
$F_{*}$ est $m$-born\'e pour un certain entier $m$ (au sens
de \cite{hagI}), c'est \`a dire \'equivalent \`a son $m$-\`eme 
co-squelette. Un argument standard permet alors de se ramener
au cas o\`u l'hyper-recouvrement  est $0$-born\'e, c'est \`a dire o\`u il
s'agit du nerf d'un recouvrement plat 
$$\{X_{i}=Spec\, B_{0,i} \longrightarrow X\}$$
(voir la preuve du th\'eor\`eme \cite[Lem. 3.4.2]{hagI}). 

Pour montrer que $A \longrightarrow Holim_{[n]\in \Delta}B_{n}$
est alors un isomorphisme dans $Ho(Comm(SEns))$, on 
utilise le fait que $\{X_{i} \longrightarrow X\}$ soit un recouvrement plat, 
et donc le fait qu'il suffise de montrer que pour tout $i$, le morphisme
induit par le changement de bases $-\times_{A}^{\mathbb{L}}B_{0,i}$ 
$$B_{0,i} \longrightarrow Holim_{[n]\in \Delta}B_{n}\times_{A}^{\mathbb{L}}B_{0,i}$$
est un isomorphisme dans $Ho(Comm(SEns))$ (noter que
la limite homotopique est finie et donc qu'elle commute au changement
de bases par platitude). On se ram\`ene ainsi 
au cas o\`u le recouvrement plat $\{X_{i} \longrightarrow X\}$
poss\`ede une section $X \longrightarrow X_{i_{0}}$. Mais dans ce
cas l'hyper-recouvrement $F_{*} \longrightarrow X$ est le nerf
d'un morphisme qui poss\`ede une section, et donc induit une
\'equivalence de pr\'efaisceaux simpliciaux
$$Hocolim_{[n] \in \Delta} F_{n} \longrightarrow X$$
est une \'equivalence dans la cat\'egorie de mod\`eles
des pr\'e-champs $Aff_{C}^{\wedge}$. Ceci implique que pour
tout $Y\in Aff_{C}$, on a 
$$Map_{Aff_{C}^{\wedge}}(X,Y)\simeq
Holim_{[n]\in \Delta}Map_{Aff_{C}^{\wedge}}(F_{n},Y),$$
ce que l'on voit, par le lemme de Yoneda, \^etre \'equivalent au fait que le morphisme
$$A \longrightarrow Holim_{[n]\in \Delta}B_{n}$$
soit un isomorphisme dans $Ho(Comm(SEns))$. \hfill $\Box$ \\

La formalisme d\'ecrit au paragraphe \S 5.1 nous donne alors une
notion de sch\'ema au-dessus de $SEns$. 

\begin{df}\label{d16}
La cat\'egorie des sch\'emas relatifs \`a $SEns$ sera not\'ee
$\mathbb{S}_{1}-Sch$. Ses objets seront appel\'es
des \emph{$\mathbb{S}_{1}$-sch\'emas}. 
\end{df}

Le foncteur 
$$\pi_{0} : SEns \longrightarrow Ens$$
est de Quillen \`a gauche et mono\"\i dal. Il satisfait aux conditions
de la proposition \ref{p10}, et donc induit un foncteur de changements
de bases
$$-\otimes_{\mathbb{S}_{1}}\mathbb{F}_{1} : \mathbb{S}_{1}-Sch
\longrightarrow \mathbb{F}_{1}-Sch.$$
On remarque que ce foncteur poss\`ede un adjoint \`a gauche pleinement fid\`ele
$$i : \mathbb{F}_{1}-Sch \longrightarrow \mathbb{S}_{1}-Sch.$$
On consid\`ere aussi le foncteur de Quillen \`a droite
$$\mathcal{GS} \longrightarrow SEns$$
qui envoit un objet $F\in \mathcal{GS}$ sur son ensemble simplicial 
sous-jacent $F(1^{+})\in SEns$ (voir le paragraphe pr\'ec\'edent pour
les notations). Ce foncteur poss\`ede un adjoint \`a gauche
$$SEns \longrightarrow \mathcal{GS}$$
qui v\'erifie les conditions de la proposition \ref{p10}. Ainsi, il d\'efinit un 
foncteur de changement de bases
$$-\otimes_{\mathbb{S}_{1}}\mathbb{S} : \mathbb{S}_{1}-Sch
\longrightarrow \mathbb{S}-Sch.$$
On remarquera que le diagramme suivant
$$\xymatrix{
\mathbb{S}_{1}-Sch \ar[rr]^-{-\otimes_{\mathbb{S}_{1}}\mathbb{S}} \ar[dd]_-{-\otimes_{\mathbb{S}_{1}}\mathbb{F}_{1}} & & \mathbb{S}-Sch 
\ar[dd]^-{-\otimes_{\mathbb{S}}\mathbb{Z}} \\
 &  & \\
\mathbb{F}_{1}-Sch \ar[rr]^-{-\otimes_{\mathbb{F}_{1}}\mathbb{Z}} & & 
\mathbb{Z}-Sch}$$
est commutatif \`a isomorphisme naturel pr\`es. Ceci permet 
en quelque sorte de comparer les $\mathbb{F}_{1}$-sch\'emas et les 
$\mathbb{S}$-sch\'emas. Comme les $\mathbb{F}_{1}$-sch\'emas peuvent
\^etre naturellement consid\'er\'es comme des
$\mathbb{S}_{1}$-sch\'emas, ceci donne un foncteur de changement
de bases des $\mathbb{F}_{1}$-sch\'emas vers les $\mathbb{S}$-sch\'emas
$$i(-)\otimes_{\mathbb{S}_{1}}\mathbb{S} : \mathbb{F}_{1}-Sch \longrightarrow
\mathbb{S}-Sch.$$
Par exemple, les vari\'et\'es toriques $X_{\mathbb{F}_{1}}(\Delta)$ 
d\'efinie en \S 4.2 donnent lieux \`a des $\mathbb{S}$-sch\'emas
$$X_{\mathbb{S}}(\Delta):=i(X_{\mathbb{F}_{1}}(\Delta))\otimes_{\mathbb{S}_{1}}\mathbb{S}.$$
Bien entendu, on a
$$X_{\mathbb{S}}(\Delta)\otimes_{\mathbb{S}}\mathbb{Z}\simeq
X_{\mathbb{Z}}(\Delta),$$
ce qui donne des vari\'et\'es toriques d\'efinies sur 
$\mathbb{S}$. 

Prenons par exemple $\Delta$ l'\'eventail ne contenant qu'un unique
cone $\mathbb{R}^{+} \subset \mathbb{R}$. Dans ce cas, 
$X_{\mathbb{S}}(\Delta)$ est $\mathbb{S}[\mathbb{N}]$ la 
$\mathbb{S}$-alg\`ebre en groupes sur le mono\"\i de $\mathbb{N}$. Elle n'est 
pas isomorphe \`a $\mathbb{Z}[\mathbb{N}]$, et ainsi
les deux $\mathbb{S}$-sch\'emas 
$i(X_{\mathbb{Z}}(\Delta))$ et $X_{\mathbb{S}}(\Delta)$ 
ne sont pas isomorphes. En r\`egle g\'en\'erale
$X_{\mathbb{S}}(\Delta)$ est une extension non triviale
de $X_{\mathbb{Z}}(\Delta)$ sur $\mathbb{S}$. 

\subsection{Le semi-anneau en spectres des entiers positifs 
$\mathbb{S}_{+}$}

Dans ce paragraphe nous posons $(C,\otimes,\mathbf{1})=(\mathcal{MS},\wedge,\mathbb{S}_{+})$, 
la cat\'egorie de mod\`eles mono\"\i dale sym\'etrique
des $\Gamma$-espaces sp\'eciaux (et non tr\`es sp\'eciaux 
comme dans \S 5.2) que nous allons maintenant d\'ecrire. 

La cat\'egorie $\mathcal{MS}$ est la cat\'egorie $\mathcal{GS}$ des
foncteurs $F : \Gamma \longrightarrow SEns$ v\'erifiant $F(*)=*$. 
La structure de mod\`eles sur $\mathcal{MS}$ est une localisation 
de la structure de mod\`eles niveaux par niveaux (fibrations et 
\'equivalences d\'efinies niveaux par niveaux) telle que ses objets locaux
soient les foncteurs $F : \Gamma \longrightarrow SEns$ tels que 
pour tout $I$ et $J$ dans $\Gamma$ le morphisme naturel
$$F(I\vee J) \longrightarrow F(I)\times F(J)$$
soit une \'equivalence faible d'ensembles simpliciaux. Il s'agit donc
des $\Gamma$-espaces sp\'eciaux et non tr\`es sp\'eciaux. 
La cat\'egorie de mod\`eles $\mathcal{GS}$ est bien entendu une
localisation de Bousfield \`a gauche de $\mathcal{MS}$. 
Tout comme pour $\mathcal{GS}$ le smash produit fait de 
$\mathcal{MS}$ une cat\'egorie de mod\`eles mono\"\i dale sym\'etrique
qui v\'erifie les conditions \ref{h2}. Pour des raisons qui deviendrons claires 
nous noterons $\mathbb{S}_{+}$ l'unit\'e de
la structure mono\"\i dale $\wedge$ sur $\mathcal{MS}$. 

Enfin, comme pour le cas de $\mathcal{MS}$ nous 
noterons pour $X\in \mathcal{MS}$
$$\pi_{i}(X):=\pi_{i}(RX(1^{+}),e),$$
o\`u $RX$ est un mod\`ele fibrant de $X$, et $e \in RX(1^{+})$ est le point de base
naturel. Plus g\'en\'eralement, pour $x\in RX(1^{+})$ on pose
$$\pi_{0}(X):=\pi_{0}(RX(1^{+})) \qquad \pi_{i}(X,x):=\pi_{i}(RX(1^{+}),x).$$

\begin{prop}\label{pconj2}
La cat\'egorie de mod\`eles mono\"\i dale sym\'etrique
$(\mathcal{MS},\wedge,\mathbb{S}_{+})$ v\'erifie l'hypoth\`ese \ref{h3}. 
\end{prop}

\textit{Preuve:} C'est essentiellement la m\^eme que pour
la proposition \ref{pconj1}. \hfill $\Box$ \\

Par le formalisme g\'en\'eral de \S 5.1 nous obtenons ainsi une 
notion de sch\'ema  relatif \`a $(\mathcal{MS},\wedge,\mathbb{S}_{+})$. 

\begin{df}\label{d17}
La cat\'egorie des sch\'emas relatifs \`a $\mathcal{MS}$ sera not\'ee
$\mathbb{S}_{+}-Sch$. Ses objets seront appel\'es
des \emph{$\mathbb{S}_{+}$-sch\'emas}. 
\end{df}

Les objets de $Comm(\mathcal{MS})$ peuvent \^etre consid\'er\'es
comme des \emph{semi-anneaux en spectres}. L'objet unit\'e
$\mathbb{S}_{+}$ doit alors \^etre pens\'e comme le semi-anneau en spectres
des entiers positifs. Noter qu'on a
$\mathbb{S}_{+}=\mathbb{S}$, mais que leurs mod\`eles fibrants
pris dans $\mathcal{MS}$ et dans $\mathcal{GS}$ ne sont 
pas \'equivalents. Nous avons d\'ej\`a mention\'e que le mod\`ele
fibrant de $\mathbb{S}$ dans $\mathcal{GS}$ est un mod\`ele
du spectre en sph\`ere, et en particulier que
les groupes $\pi_{i}(\mathbb{S})$  sont les groupes d'homotopie stables
des sph\`eres. 

En contre partie, on peut montrer qu'un mod\`ele
fibrant de $\mathbb{S}_{+}$ dans $\mathbb{MS}$ est donn\'e par 
le nerf du groupo\"\i de des ensembles finis. Plus pr\'ecis\'ement, 
le groupo\"\i de des ensembles finis $ES$ poss\`ede une structure
mono\"\i dale sym\'etrique donn\'ee par la somme disjointe. On peut alors
construire un pseudo-foncteur
$$ES : \Gamma \longrightarrow Cat$$
qui envoit $n^{+}$ sur $ES^{n}$ et o\`u les foncteurs de transitions sont
donn\'es par la structure mono\"\i dale. En appliquant le proc\'ed\'e de strictification
usuelle on en d\'eduit un foncteur (au sens strict)
$$ES : \Gamma \longrightarrow Cat$$
qui compos\'e avec le foncteur nerf $N : Cat \longrightarrow SEns$ donne
un foncteur
$$N(ES) : \Gamma \longrightarrow SEns,$$
qui est un $\Gamma$-espace sp\'ecial. On montre alors que ce $\Gamma$-espace
est un mod\`ele fibrant de $\mathbb{S}_{+}$ dans $\mathcal{MS}$. En particulier, 
les groupes d'homotopie de $N(ES)$, sont donn\'es par
$$\pi_{0}(\mathbb{S}_{+}):=\pi_{0}(N(ES)(1^{+}))\simeq \mathbb{N}
\qquad
\pi_{1}(\mathbb{S}_{+},n)\simeq \Sigma_{n} \qquad
\pi_{i}(\mathbb{S}_{+},n)=0 \; \forall \; i>1.$$
On voit ainsi que $\mathbb{S}_{+}$ est un semi-anneau en spectres dont un mod\`ele
est le groupo\"\i de des ensembles finis avec ses structures mono\"\i dales 
donn\'ees par le somme disjointe (l'addition) et le produit direct (la multiplication), et qui est 
donc une version homotopique des entiers naturels. 
L'anneau en spectres $\mathbb{S}$ est par construction l'anneau 
obtenu par compl\'etion en groupes du semi-anneau $\mathbb{S}_{+}$. On retrouve
ici le th\'eor\`eme de Quillen-Priddy, affirmant que la compl\'etion en groupes
du groupo\"\i des des ensembles finis est \'equivalent au spectre en sph\`ere. Il est bienvenu de comparer
le morphisme naturel $\mathbb{S}_{+} \longrightarrow \mathbb{S}$ avec le morphisme naturel
$\mathbb{N} \longrightarrow \mathbb{Z}$.

On dispose d'un foncteur d'inclusion
$$\mathbb{N}-Mod \longrightarrow \mathcal{MS}$$
des mono\"\i des commutatifs dans $\mathcal{MS}$, qui est de Quillen \`a droite. 
Son adjoint \`a gauche est le foncteur de troncation
$$\pi_{0} : \mathcal{MS} \longrightarrow \mathbb{N}-Mod$$
qui envoit un objet $F$ sur $\pi_{0}(RF(1^{+})$. On v\'erifie que les
hypoth\`eses de la proposition \ref{p10} sont satisfaites et donc que l'on 
obtient un foncteur de changement de bases
$$-\otimes_{\mathbb{S}_{+}}\mathbb{N} : \mathbb{S}_{+}-Sch \longrightarrow
\mathbb{N}-Sch.$$
Ce foncteur poss\`ede de plus un adjoint \`a gauche pleinement fid\`ele
$$i : \mathbb{N}-Sch \longrightarrow \mathbb{S}_{+}-Sch.$$

De plus, le foncteur identit\'e
$$\mathcal{MS} \longrightarrow \mathcal{GS}$$
est mono\"\i dal sym\'etrique et de Quillen \`a gauche. Son adoint \`a droite 
induit un foncteur pleinement fid\`ele
$$Ho(\mathcal{GS}) \longrightarrow Ho(\mathcal{MS})$$
dont l'image est form\'e des $\Gamma$-espaces sp\'eciaux. 
On peut v\'erifier que les conditions de 
la proposition \ref{p10} sont satisfaites et que l'on obtient ainsi un foncteur
de changement de bases
$$\mathbb{S}_{+}-Sch \longrightarrow \mathbb{S}-Sch$$
qui poss\`ede encore un adjoint \`a gauche pleinement fid\`ele
$$\mathbb{S}-Sch \longrightarrow \mathbb{S}_{+}-Sch.$$
Bien entendu on dispose d'un diagramme commutatif \`a isomorphisme 
naturel pr\`es
$$\xymatrix{
\mathbb{S}_{+}-Sch \ar[rr]^-{-\otimes_{\mathbb{S}_{+}}\mathbb{S}} 
\ar[dd]_-{-\otimes_{\mathbb{S}_{+}}\mathbb{N}} & & \mathbb{S}-Sch 
\ar[dd]^-{-\otimes_{\mathbb{S}}\mathbb{Z}} \\
 &  & \\
\mathbb{N}-Sch \ar[rr]_-{-\otimes_{\mathbb{N}}\mathbb{Z}} & & 
\mathbb{Z}-Sch,}$$
ce qui permet de comparer les $\mathbb{N}$-sch\'emas et 
les $\mathbb{S}$-sch\'emas. Il faut noter que pour
$X\in \mathbb{N}-Sch$ le morphisme naturel
$$i(X)\otimes_{\mathbb{S}_{+}}\mathbb{S} \longrightarrow
i(X\otimes_{\mathbb{N}}\mathbb{Z})$$
n'est pas un isomorphisme en g\'en\'eral. En effet, 
le foncteur $-\otimes_{\mathbb{S}_{+}}\mathbb{S}$
est donn\'e par une version homotopique de la construction 
compl\'etion en groupes, des semi-anneaux en spectres vers
les anneaux en spectres, et l'on sait bien qu'appliqu\'e \`a un 
semi-anneaux discret on peut obtenir un anneau en 
spectres non discret. 

Tout comme pour le cas de $\mathbb{S}$, on peut 
d\'efinir un $\mathbb{S}_{+}$-sch\'ema affine, classifiant les 
auto-\'equivalence du $\mathbb{S}_{+}$-module libre de rang $n$. Nous le noterons 
$Gl_{n,\mathbb{S}_{+}}$. On montre facilement que l'on a les deux
propri\'et\'es suivantes

\begin{itemize}
\item L'ensemble simplicial $Gl_{n,\mathbb{S}_{+}}(\mathbb{S}_{+})$
est discret et  \'equivalent \`a $\Sigma_{n}$. 
\item Le morphisme naturel
$$Gl_{n,\mathbb{S}_{+}}\otimes_{\mathbb{S}_{+}}\mathbb{S} \longrightarrow
Gl_{n,\mathbb{S}}$$
est un isomorphisme de $\mathbb{S}$-sch\'emas affines. 
\end{itemize}

On montre aussi, tout comme pour le cas de $\mathbb{S}$ que le morphisme naturel
$$Gl_{n,\mathbb{S}_{+}} \longrightarrow i(Gl_{n,\mathbb{N}})$$
n'est pas un isomorphisme dans $\mathbb{S}_{+}$-Sch. 

Enfin, le foncteur de Quillen \`a droite
$$\mathcal{MS} \longrightarrow SEns$$
qui envoit $F$ sur $F(1^{+})$, 
dont l'adjoint \`a gauche v\'erifie les conditions
de la proposition \ref{p10} donne un foncteur de changements de base
$$\mathbb{S}_{1}-Sch \longrightarrow \mathbb{S}_{+}-Sch.$$
On dispose alors d'un second diagramme commutatif \`a homotopie pr\`es
$$\xymatrix{
\mathbb{S}_{1}-Sch \ar[rr]^-{-\otimes_{\mathbb{S}_{1}}\mathbb{S}_{+}} 
\ar[dd]_-{-\otimes_{\mathbb{S}_{1}}\mathbb{F}_{1}} & & \mathbb{S}_{+}-Sch 
\ar[dd]^-{-\otimes_{\mathbb{S}_{+}}\mathbb{N}} \\
 &  & \\
\mathbb{F}_{1}-Sch \ar[rr]_-{-\otimes_{\mathbb{F}_{1}}\mathbb{N}} & & 
\mathbb{N}-Sch.}$$

En conclusion, on dispose d'un diagramme commutatif \`a isomorphisme naturel pr\`es
$$\xymatrix{
\mathbb{S}_{1}-Sch \ar[rr]^-{-\otimes_{\mathbb{S}_{1}}\mathbb{S}_{+}} \ar[dd]_-{-\otimes_{\mathbb{S}_{1}}\mathbb{F}_{1}} & & \mathbb{S}_{+}-Sch 
\ar[dd]^-{-\otimes_{\mathbb{S}_{+}}\mathbb{N}} \ar[rr]^-{-\otimes_{\mathbb{S}_{+}}\mathbb{S}} & & \mathbb{S}-Sch
\ar[dd]^-{-\otimes_{\mathbb{S}}\mathbb{Z}} \\
 &  & & & \\
\mathbb{F}_{1}-Sch \ar[rr]_-{-\otimes_{\mathbb{F}_{1}}\mathbb{N}} & & 
\mathbb{N}-Sch \ar[rr]_-{-\otimes_{\mathbb{N}}\mathbb{Z}} & & \mathbb{Z}-Sch.}$$
Sch\'ematiquement on repr\'esente ce diagramme par le diagramme suivant
$$\xymatrix{
Spec\, \mathbb{Z} \ar[r] \ar[d] & Spec\, \mathbb{N} \ar[r] \ar[d] & Spec\, \mathbb{F}_{1} \ar[d] \\
Spec\, \mathbb{S} \ar[r] & Spec\, \mathbb{S}_{+} \ar[r] & Spec\, \mathbb{S}_{1}.}$$


\begin{thebibliography}{99}


\bibitem[Bl]{bl} B. Blander, \textit{Local projective model structure on simpicial presheaves},
$K$-theory \textbf{24} (2001) No. $3$, 283-301.

\bibitem[De]{de} Deitmar, \textit{Schemes over $F_1$}, 
dans \textit{Number fields and function fields---two parallel worlds},  
87--100, Progr. Math. \textbf{239}, Birkhäuser Boston, Boston, MA, 2005. 


\bibitem[D]{del} P. Deligne, \textit{Cat\'egories tannakiennes}, 
The Grothendieck Festschrift, Vol. II, 111--195, Progr. Math., 87, Birkhäuser Boston, Boston, MA, 1990.

\bibitem[EKMM]{ekmm} A.D. Elmendorf, I. Kriz, M.A. Mandell, J.P. May, \textit{Rings, modules, and  
algebras in stable homotopy theory}, Mathematical Surveys and Monographs, vol. $47$,  
American Mathematical Society, Providence, $RI$, $1997$.  
 
\bibitem[Ha]{ha} M. Hakim, \textit{Topos annel\'es et sch\'emas relatifs},
Ergebnisse der Mathematik und ihrer Grenzgebiete, Band $64$. Springer-Verlag
Berlin-New York, $1972$.
 
\bibitem[Hol]{hol} S. Hollander, \textit{A homotopy theory for stacks}, pr\'e-publication
math.AT/0110247.

\bibitem[Ho]{ho} M. Hovey, \textit{Model categories}, Mathematical surveys and monographs, Vol. $\mathbf{63}$,
Amer. Math. Soc., Providence 1998.

\bibitem[Hu]{hut}  T. Huettemann, \textit{Algebraic K-Theory of non-linear projective toric varieties},
J. Pure Appl. Algebra  \textbf{170}  (2002),  no. 2-3, 185--242. 

\bibitem[Ja]{ja} J. F. Jardine, \textit{Simplicial presheaves}, J. Pure and Appl. Algebra $\mathbf{47}$ (1987),
35-87.

\bibitem[La-Mo]{lm} G. Laumon, L. Moret-Bailly, 
\textit{Champs alg\'ebriques}, 
Ergebnisse der Mathematik und ihrer Grenzgebiete, Band $39$. Springer-Verlag
Berlin, $2000$.

\bibitem[Ma-Mo]{mm}
S. Mac Lane, I. Moerdijk, 
\textit{Sheaves in geometry and logic. A first introduction to topos theory}, 
Universitext. Springer-Verlag, New York, 1994.

\bibitem[M-M-S-S]{mmss} M. Mandel, P. May, S. Schwede, B. Shipley, 
\textit{Model categories of diagram spectra},  
Proc. London Math. Soc. (3)  $\mathbf{82}$  (2001),  no. 2, 441--512.

\bibitem[Mi]{mi} Milne, \textit{Etale cohomology}, 
Princeton Mathematical Series, 33. 
Princeton University Press, Princeton, N.J., 1980. xiii+323. 

\bibitem[O]{od} T. Oda, \textit{Convex bodies and algebraic geometry}, 
Ergebnisse der Mathematik und ihrer Grenzgebiete, Band $15$. Springer-Verlag
Berlin-New York, $1972$.

\bibitem[R-S-T]{trop} J. Richter-Gebert, B. Sturmfels, T. Theobald, 
\textit{First steps in tropical geometry}, 
Proc. Conference on Idempotent Mathematics and Mathematical Physics, 
Vienna 2003 (G.L. Litvinov and V.P. Maslov, eds.), Contemporary Mathematics, AMS.

\bibitem[Sa]{sav}  N. Saavedra, \textit{Cat\'egories Tannakiennes},  
Lecture Notes in Mathematics, Vol. $\mathbf{265}$ 
Springer-Verlag, Berlin-New York, 1972. ii+418 pp.

\bibitem[S]{sch} S. Schwede, \textit{Stable homotopical algebra and $\Gamma$-spaces},  
Math. Proc. Cambridge Philos. Soc.  $\mathbf{126}$  (1999),  no. 2, 329--356.

\bibitem[S-S]{ss} S. Schwede, B. Shipley, \textit{Algebras and modules
in monoidal model categories}, Proc. London Math. Soc. $(3)$
$\mathbf{80}$ ($2000$), $491-511$.


\bibitem[So]{so} C. Soul\'e, \textit{Les vari\'et\'es sur le corps \`a un \'el\'ement},
Mosc. Math. J.  \textbf{4}  (2004),  no. 1, 217--244.

\bibitem[Sp]{sp} M. Spitzweck, \textit{Operads, algebras and modules in model categories and motives},
Ph.D. Thesis, Mathematisches Instit\"ut, Friedrich-Wilhelms-Universit\"at Bonn (2001), accessible \`a
\textsf{http://www.uni-math.gwdg.de/spitz/}.


\bibitem[To-Va]{tv} B. To\"en, M. Vaqui\'e, \textit{Moduli of objects in dg-categories}, 
\`a parraitre aux Annales Sci. de l'ENS,
pr\'e-publication math.AG/0503269.

\bibitem[HAGDAG]{hagdag} B. To\"en, G. Vezzosi, 
\textit{From HAG to DAG: derived moduli stacks}, in Axiomatic, enriched and motivic homotopy theory, 173--216, NATO Sci. Ser. II Math. Phys. Chem., 131, Kluwer Acad. Publ., Dordrecht, 2004.

\bibitem[HAGI]{hagI} B. To\"en, G. Vezzosi, \textit{Homotopical algebraic geometry I: Topos theory},
Adv. Math. $\mathbf{193}$ (2005), no. 2, 257--372.

\bibitem[HAGII]{hagII} B. To\"en, G. Vezzosi,
\textit{Homotopical algebraic geometry II: Geometric stacks and applications}, \`a parraitre
dans Memoires of the AMS, 
pr\'e-publication math.AG/0404373.

\bibitem[Ve]{ve} G. Vezzosi, \textit{A sketchy note on enriched homotopical 
topologies and enriched homotopical stacks}, 
pr\'e-publication math.CT/0507447.

\end{thebibliography}
\end{document}